\documentclass[10pt,oneside]{article}
\usepackage{amssymb}
\usepackage[cp1250]{inputenc}
\usepackage[T1]{fontenc}
\usepackage{color}

\begin{document}
\newcommand{\qed}{\rule{1.5mm}{1.5mm}}
\newcommand{\proof}{\textit{Proof. }}
\newcommand{\ccon}{\rightarrowtail}
\newtheorem{theorem}{Theorem}[section]
\newtheorem{lemma}[theorem]{Lemma}
\newtheorem{remark}[theorem]{Remark}
\newtheorem{example}[theorem]{Example}
\newtheorem{corollary}[theorem]{Corollary}
\newtheorem{proposition}[theorem]{Proposition}
\newtheorem{claim}[theorem]{Claim}
%\numberwithin{equation}{section}
\begin{center}
{\LARGE\textbf{Effective approximation of the solutions of algebraic equations}\vspace*{3mm}}\\
{\large Marcin Bilski\footnote{Department of Mathematics and Computer Science, Jagiellonian University, \L
ojasiewicza 6, 30-348 Krak\'ow,
Poland. E-mail address: marcin.bilski@im.uj.edu.pl
}}\hspace*{10mm}
{\large Peter Scheiblechner\footnote{Lucerne University of Applied Sciences and Arts
    School of Engineering and Architecture
           6048 Horw, Switzerland,
E-mail address: peter.scheiblechner@hslu.ch\vspace*{1mm}\\
M. Bilski was partially supported by the NCN grant 2014/13/B/ST1/00543
}}\vspace*{2mm}\\
\end{center}
{\small\textbf{Abstract.} Let $F$ be a holomorphic map whose
components satisfy some polynomial relations. We present an algorithm for
constructing Nash maps locally approximating $F,$ whose components
satisfy the same relations.\vspace*{2mm}\vspace*{0mm}}
%+++++++++++++++++++++++++++++++++++++++++++++++++++++++++++++++++

\section{Introduction}\label{introduction}

Let $Q_1,\ldots,Q_q$ be polynomials in $\hat{m}$ complex variables
and let $F:U\rightarrow\mathbf{C}^{\hat{m}}$ be a holomorphic map
such that
$$Q_1(F(x))=\cdots=Q_q(F(x))=0$$
for every $x\in U,$ where $U$ is a domain in $\mathbf{C}^n.$ Let us fix any point $x_0\in U.$ The aim of this
paper is to present an algorithm for constructing a complex Nash map $F^{\nu}$ uniformly approximating  $F,$ in
some neighborhood $U_0$ of $x_0,$ such that
$$Q_1(F^{\nu}(x))=\cdots=Q_q(F^{\nu}(x))=0$$
for every $x\in U_0.$ The correctness of the
algorithm will follow from the proof of the fact that such approximations
always exist, presented in Section \ref{mainIprof} (see Theorem
\ref{main}).

The existence of local approximation of the solutions of algebraic or analytic equations was investigated in
\cite{Ar1}, \cite{Ar2} and \cite{Ar3}, and Theorem~\ref{main} can be derived from the results of these papers. A
more general global version of Theorem~\ref{main} was proved in \cite{Lem}, pp 338-339, by means of general
N\'eron desingularization (which makes it not suitable for effective constructions; for approximating over the
field of real numbers see also \cite{CRS}). In order to obtain an effective procedure of approximation, in the
present paper we base the proof on a combination of some of the ideas of \cite{Ar1}, \cite{vD} and \cite{B6}. This
allows us to reduce the problem to the case where $q=1$ and $Q_1$ is of the form
$$Q_1(y_1,\ldots,y_{\hat{m}})=y_{{\hat{m}}}^s+y_{\hat{m}}^{s-1}p_1(y_1,\ldots,y_{\hat{m}-1})+\cdots+p_s(y_1,\ldots,y_{\hat{m}-1}),$$ where
$s\in\mathbf{N}$ and $p_1,\ldots,p_s\in\mathbf{C}[y_1,\ldots,y_{\hat{m}-1}].$ Then we can use the following claim:
in order to find complex Nash approximations $F^{\nu}:U_0\rightarrow\mathbf{C}^{\hat{m}}$ of $F|_{U_0}$ such that
$Q_1(F^{\nu}(x))=0$ it is sufficient to find complex Nash approximations $\bar{F}^{\nu}$ of $F|_{U_0}$ such that
$Q_1(\bar{F}^{\nu}(x))=d(x)(\frac{\partial Q_1}{\partial{y_{\hat{m}}}}(\bar{F}^{\nu}(x)))^2,$ where $d(x)$ is a
complex Nash function bounded by a small constant.

The claim was proved in \cite{vD}, p. 393. Since it will be used for several times in the present paper, we sketch
the idea of its proof. First, using the Taylor expansion of {$Q_1$ about $\bar{F}^{\nu}(x)$ with respect to
$y_{\hat{m}},$ we get
$$Q_1(\bar{F}^{\nu}(x)+(0,\ldots,0,h))=$$
$$Q_1(\bar{F}^{\nu}(x))+\frac{\partial Q_1}{\partial{y_{\hat{m}}}}(\bar{F}^{\nu}(x))h+\sum_{i=2}^s\frac{1}{i!}\frac{\partial^i
Q_1}{\partial{y_{\hat{m}}^i}}(\bar{F}^{\nu}(x))h^i=$$

$$d(x)(\frac{\partial Q_1}{\partial{y_{\hat{m}}}}(\bar{F}^{\nu}(x)))^2+\frac{\partial Q_1}{\partial{y_{\hat{m}}}}(\bar{F}^{\nu}(x))h+\sum_{i=2}^s\frac{1}{i!}\frac{\partial^i
Q_1}{\partial{y_{\hat{m}}^i}}(\bar{F}^{\nu}(x))h^i.$$ Next, taking $h=d(x)\frac{\partial
Q_1}{\partial{y_{\hat{m}}}}(\bar{F}^{\nu}(x))Y,$ where $Y$ is a new variable, we obtain

$$Q_1(\bar{F}^{\nu}(x)+(0,\ldots,0,d(x)\frac{\partial Q_1}{\partial{y_{\hat{m}}}}(\bar{F}^{\nu}(x))Y))=$$
$$d(x)\big{(}\frac{\partial
Q_1}{\partial{y_{\hat{m}}}}(\bar{F}^{\nu}(x))\big{)}^2\cdot\big{(}1+Y+\sum_{i=2}^s\frac{1}{i!}\frac{\partial^i
Q_1}{\partial{y_{\hat{m}}^i}}(\bar{F}^{\nu}(x))\big{(}\frac{\partial
Q_1}{\partial{y_{\hat{m}}}}(\bar{F}^{\nu}(x))\big{)}^{i-2}(d(x))^{i-1}Y^i\big{)}.$$ The last factor of the
right-hand side of the equation above is a polynomial in $Y$ such that the coefficient  of $Y^i,$ for $i\geq 2,$
is a Nash function bounded by a small constant (since $d(x)$ is bounded by a small constant). One can show that in
such a situation there is a Nash function $Y=Y(x)$ for which the last factor vanishes identically (cf. \cite{vD},
p. 392). Then $$F^{\nu}=\bar{F}^{\nu}(x)+(0,\ldots,0,d(x)\frac{\partial
Q_1}{\partial{y_{\hat{m}}}}(\bar{F}^{\nu}(x))Y(x))$$ is the Nash approximation we look for.}

The claim discussed above was applied in a simi\-lar context to (global) maps depending on one variable in
\cite{vD}. If the number of the variables $F$ depends on is greater than $1,$ then we will decrease it using the
Weierstrass Division Theorem in the spirit of \cite{Ar1}, see also \cite{JP} pp. 295-298 or \cite{R} pp. 97-98.

To perform an algorithm we need
to represent holomorphic functions in such a way that they can constitute its
input. This will be done similarly to \cite{vH2}, \cite{vH4}.
More precisely, every component $f$ of $F$ will be defined by the following data:\vspace*{2mm}\\
(a) a procedure computing the coefficients of the Taylor expansion
of $f$ at $x_0,$\\
(b) the polyradius of a polydisc $U_f$ centered at $x_0$ on which the Taylor expansion\linebreak
\hspace*{5.5mm}of $f$ is convergent,\\
(c) a constant $M_f$ such that $\sup_{U_f}|f|<M_f.$\vspace*{2mm}

The class of functions represented in this way
is large and it contains important transcendental functions (such as $\mathrm{exp},$ $\mathrm{log},$ $\mathrm{sin},$ $\mathrm{cos}$ and many others). Moreover, it is closed under algebraic operations and differentiation which
can be performed effectively.
Furthermore, one can use computer algebra
procedures \cite{GP} to obtain effective versions of the Weierstrass Preparation and
Division Theorems. Using basic facts from complex analysis (Cauchy's integral formula), it is also possible to control accuracy of polynomial approximations of such functions
(and, consequently, accuracy of Nash approximations of $F$;
see Section \ref{examp} below for details).

In Section \ref{examp} we present algorithms for computing approximating maps. The algorithms rely on effective
versions of the Weierstrass Preparation and Division Theorems and on some techniques from computational algebraic
geometry (equidimensional decomposition \cite{GP}, \cite{JKSS}, \cite{JS}, \cite{Lec},  elimination theory
\cite{CLOS}, \cite{DL}, \cite{GP}, description of factors of iterated discriminants \cite{LMC}).

This paper is related to \cite{B6} where another method for approximation of holomorphic solutions of polynomial
equations is presented. However, the method of \cite{B6} contains steps which are not computable in the data
structure described above. More precisely, it relies on the zero-test and on factorization of monic polynomials
with holomorphic coefficients into the product of powers of reduced monic and pairwise relatively prime
polynomials with holomorphic coefficients. (The zero-test is only semi-decidable, and computability of the
factorization would imply decidability of the zero-test, cf. Remark \ref{comparis} below.) The approach proposed
in the present paper allows us to avoid these steps. First, in Section \ref{examp}, we give a partial algorithm
which does not rely on the factorization of polynomials with holomorphic coefficients, but still uses the
zero-test. Next after refining the partial algorithm, we obtain a complete one.

Our work is motivated by the natural question whether every purely dimensional analytic set can be (locally)
approximated by algebraic or Nash ones of the same pure dimension. This question is related to the classical
problem of characterizing those analytic sets which are analytically equivalent to algebraic ones (see e.g.
\cite{Ar1}, \cite{Ar2}, \cite{BoK}, \cite{BuK}, \cite{Lem}). As for algebraic approximation, the problem is
interesting especially for non-complete intersections, i.e. for analytic sets for which the number of defining
equations is greater than their codimensions. Then one cannot simply replace the defining functions by
approximating polynomials as this gives sets of strictly smaller dimensions. Nevertheless, algebraic
approximations do exist for a large subclass of the class of non-complete intersections (see \cite{B6}, \cite{B7},
\cite{B9}). In particular, every analytic set admits local algebraic approximations which in many cases can be
effectively constructed. One of the main tools used in the construction is Theorem \ref{main} discussed in the
present paper (cf. \cite{B7}, p. 284) and the algorithm following from its proof.

The organization of our paper is as follows. Sections \ref{mainIprof} and \ref{examp} are devoted to
Theorem~\ref{main} and the algorithms, respectively. Preliminary material concerning Nash maps and sets as well as
analytic sets with proper projection is gathered in Section~\ref{preli} below.

\section{Preliminaries}\label{preli}
\subsection{Nash mappings and sets}\label{prelnash}
Let $\Omega$ be an open subset of $\mathbf{C}^n$ and let $f$ be a
holomorphic function on $\Omega.$ We say that $f$ is a Nash
function at $x_0\in\Omega$ if there exist an open neighborhood $U$
of $x_0$ and a polynomial
$P:\mathbf{C}^n\times\mathbf{C}\rightarrow\mathbf{C},$ $P\neq 0,$
such that $P(x,f(x))=0$ for $x\in U.$ A holomorphic function
defined on $\Omega$ is said to be a Nash function if it is a Nash
function at every point of $\Omega.$ A holomorphic mapping defined
on $\Omega$ with values in $\mathbf{C}^N$ is said to be a Nash
mapping if each of its components is a Nash function.

{The following lemma is well known. Here we recall the proof
to emphasize its algorithmic nature.}
\begin{lemma}\label{algNashfunct}{Let $f,g:\Omega\rightarrow\mathbf{C}$ be Nash functions
and let $P\in\mathbf{C}[x,y], Q\in\mathbf{C}[x,z]$ be monic in $y, z$, respectively, such that
$$P(x,f(x))=0, Q(x,g(x))=0\mbox{ for }x\in\Omega.$$ Then there are $S, R\in \mathbf{C}[x,u]$, monic in $u$, such that
$$S(x,(f+g)(x))=0, R(x,(f\cdot g)(x))=0 \mbox{ for }x\in\Omega.$$ Given $P, Q,$ one can effectively compute $S, R.$
Moreover,
the degrees of $S, R$ in $u$ are bounded by a constant depending only on the degrees of $P, Q$ in $y, z.$}
\end{lemma}
\textit{Proof. }{ Let $\diamond$ denote either multiplication or addition in $\mathbf{C}$.
Let $\pi: \mathbf{C}_x^n\times\mathbf{C}_u\times\mathbf{C}_y\times\mathbf{C}_z\rightarrow\mathbf{C}_x^n\times\mathbf{C}_u $
denote the natural projection. Consider the algebraic set
$$V=\{(x, u, y, z)\in\mathbf{C}^n\times\mathbf{C}\times\mathbf{C}\times\mathbf{C}: P(x,y)=Q(x,z)=0, u=y\diamond z\}.$$
Since $P, Q$ are monic in $y,z$, and $u-(y\diamond z)$ is monic in $u$ we see that the map
$\pi|_V:V\rightarrow\mathbf{C}^n_x\times\mathbf{C}_u$ is proper and also the projection of $\pi(V)$ onto
$\mathbf{C}^n_x$ is proper. Consequently $\pi(V)$ is an algebraic hypersurface described by a single polynomial
denoted by {$S(x,u)$} (if $\diamond=+$) or {$R(x,u)$} (if $\diamond=\cdot$). Clearly, $S(x,(f+g)(x))=0$ and $R(x,
(f\cdot g)(x))=0.$ Observe that, by properness of the projection of $\pi(V)$ onto $\mathbf{C}^n_x,$ the
polynomials $S, R$ are monic in $u.$}

{Note that $S, R$ can be effectively computed and the degrees in $u$ of $S, R$ are bounded by
a constant depending only on the degrees of $P, Q$ in $y, z,$ respectively. This follows by
the fact that $S(x,u), R(x,u)$ are obtained by taking the resultants of pairs of polynomials
in $y$ with coefficients in $\mathbf{C}[x,u].$ Namely, $S(x,u)=\mathrm{res}_y(P(x,y),
Q(x,u-y))$ and $R(x,u)=\mathrm{res}_y(P(x,y), y^dQ(x,\frac{u}{y})),$ where $d$ is the degree
of $Q$ in $z.$}

{In other words, the projection of $V$ to $\mathbf{C}^n_x\times\mathbf{C}_u$ corresponds to the elimination of
$y,z$ from $P, Q, u-(y\diamond z),$ which is algorithmic (see also \cite{BeWa}, \cite{CLOS}, Chapter 3, or
\cite{GP}, pp. 69-73) .\qed}\vspace*{2mm}

A subset $Y$ of an open set $\Omega\subset\mathbf{C}^n$ is said to be a Nash subset of $\Omega$ if and only if for
every $y_0\in\Omega$ there exists a neighborhood $U$ of $y_0$ in $\Omega$ and there exist Nash functions
$f_1,\ldots,f_s$ on $U$ such that $$Y\cap U=\{x\in U: f_1(x)=\cdots=f_s(x)=0\}.$$ We will use the following fact
from \cite{Tw}, p. 239. Let $\pi:\Omega\times\mathbf{C}^k\rightarrow\Omega$ denote the natural projection.
\begin{theorem}
\label{projnash} Let $X$ be a Nash subset of
$\Omega\times\mathbf{C}^k$ such that $\pi|_{X}:X\rightarrow\Omega$
is a proper mapping. Then $\pi(X)$ is a Nash subset of $\Omega$
and $\mathrm{dim}(X)=\mathrm{dim}(\pi(X)).$
\end{theorem}
The fact from \cite{Tw} stated below explains the relation between
Nash and algebraic sets.
\begin{theorem}\label{charactNash}Let $X$ be a Nash subset of an open set
$\Omega\subset\mathbf{C}^n.$ Then every analytic irreducible
component of $X$ is an irreducible Nash subset of $\Omega.$
Moreover, if $X$ is irreducible then there exists an algebraic
subset $Y$ of $\mathbf{C}^n$ such that $X$ is an analytic
irreducible component of $Y\cap\Omega.$
\end{theorem}
\begin{corollary}{\label{functHypNash}Let $U$ be an open subset of $\mathbf{C}^n.$ Let $f:U\rightarrow\mathbf{C}$ be a holomorphic function and
let $N\subset U\times\mathbf{C}$ be a Nash hypersurface such that $\mathrm{graph}(f)\subset N.$
Then $f$ is a Nash function.}
\end{corollary}
{\textit{Proof.} One may assume that $U$ is connected (because each of its connected components
can be considered separately). If $N$ is reducible, then replace it by its analytic irreducible component containing
$\mathrm{graph}(f)$ which, by Theorem~\ref{charactNash}, is a Nash set. Again by Theorem~\ref{charactNash}, there is an algebraic hypersurface
$M\subset\mathbf{C}^n\times\mathbf{C}$ such that $N\subset M.$ Let $P\in\mathbf{C}[x_1,\ldots,x_n,z]$ be a (non-zero) polynomial such that
$\{P=0\}=M.$ Then $P(x,f(x))=0$ for every $x\in U.$ \qed}
%**************************************************************************************************
\subsection{Analytic sets}\label{setswithprop} {For any map $f,$ by $\mathrm{graph}(f)$
we denote the set of all $(x,y)$ such that $f(x)=y$.}

Let $U, U'$ be domains in $\mathbf{C}^n,\mathbf{C}^k,$ respectively,
and let
$\pi:\mathbf{C}^n\times\mathbf{C}^k\rightarrow\mathbf{C}^n$ denote
the natural projection. For any purely $n$-dimensional analytic
subset $Y$ of $U\times U'$ such that $\pi|_Y:Y\rightarrow U$ is a
proper mapping, we denote by $\mathcal{S}(Y,\pi)$ the set of
singular points of $\pi|_{Y}:$
$$\mathcal{S}(Y,\pi)=\mathrm{Sing}(Y)\cup\{x\in \mathrm{Reg}(Y):(\pi|_Y)'(x) \mbox{ is not an isomorphism}\}.$$
We often write $\mathcal{S}(Y)$ instead of $\mathcal{S}(Y,\pi)$
when it is clear which projection is taken into consideration.

It is well known that $\mathcal{S}(Y)$ is an analytic subset of $U\times U',$
$\mathrm{dim}(\mathcal{S}(Y))<n$ (cf. \cite{Ch}, p. 50), hence, by the Remmert Theorem
$\pi(\mathcal{S}(Y))$ is also analytic. Moreover, the following holds. The mapping $\pi|_{Y}$
is surjective
and open and there exists an integer $s=s(\pi|_{Y})$ such that:\vspace*{2mm}\\
(1) $\sharp(\pi|_{Y})^{-1}(\{a\})<s$ for
$a\in\pi(\mathcal{S}(Y)),$\\
(2) $\sharp(\pi|_{Y})^{-1}(\{a\})=s$ for $a\in U\setminus
\pi(\mathcal{S}(Y)),$\\
\\
(3) for every $a\in U\setminus \pi(\mathcal{S}(Y))$ there exist a neighborhood $W$ of $a$ and
holomor-\linebreak \hspace*{5.4mm}phic mappings $f_1,\ldots,f_s:W\rightarrow U'$ such that {
$\mathrm{graph}(f_i)\cap \mathrm{graph}(f_j)=\emptyset$ for\linebreak \hspace*{5.4mm}$i\neq
j$ and $\mathrm{graph}(f_1)\cup\ldots\cup \mathrm{graph}(f_s)=$$(W\times U')\cap
Y.$}\vspace*{2mm}

Let $E$ be a purely $n$-dimensional analytic subset of $U\times
U'$ with proper projection onto a domain $U\subset\mathbf{C}^n,$
where $U'$ is a domain in $\mathbf{C}.$ Then there is a monic
polynomial $p\in\mathcal{O}(U)[z]$ (i.e. a polynomial in $z$ whose leading coefficient is $1$)
such that $E=\{(x,z)\in
U\times\mathbf{C}:p(x,z)=0\}$ and the discriminant $\Delta_p$ of
$p$ is not identically zero. $p$ will be called the optimal
polynomial for $E.$ We have $\tilde{\pi}(\mathcal{S}(E))=\{x\in
U:\Delta_p(x)=0\},$ where
$\tilde{\pi}:U\times\mathbf{C}\rightarrow U$ is the natural
projection. If $E$ is algebraic and $U=\mathbf{C}^n, U'=\mathbf{C},$ then the coefficients
of the optimal polynomial $p$ are polynomials.

Let $V$ be an analytic subset of $\mathbf{C}^n\times\mathbf{C}^k$ with proper projection onto
$\mathbf{C}^n,$ let
$\Phi_L:\mathbf{C}^n\times\mathbf{C}^k\rightarrow\mathbf{C}^n\times\mathbf{C}$ be given by
$\Phi_L(u,v)=(u,L(v)),$ where $L$ is $\mathbf{C}$-linear. Then
$\Phi_L|_V:V\rightarrow\mathbf{C}^n\times\mathbf{C}$ is a proper map. Indeed, fix any compact
subset $K$ of $\mathbf{C}^n\times\mathbf{C}.$ Let $K'$ be the image of the projection of $K$
onto $\mathbf{C}^n.$ Then $(\Phi_L|_V)^{-1}(K)$ is clearly closed (in the Euclidean
topology). It is also bounded because $(\Phi_L|_V)^{-1}(K)\subset (K'\times\mathbf{C}^k)\cap
V,$ so $\Phi_L|_V$ is proper. Consequently, $\Phi_L(V)$ is analytic (or algebraic if $V$ is
such). {Also observe that if generic fibers in $\Phi_L(V)$ and in $V$ over $\mathbf{C}^n$
have the same number of elements, then
$\Phi_L(\mathcal{S}(V))\subset\mathcal{S}(\Phi_L(V)).$} {Indeed, pick any
$(a,b)\in\mathcal{S}(V)\subset \mathbf{C}^n\times\mathbf{C}^k$ and a polydisc $P_1\times
P_2\subset\mathbf{C}^n\times\mathbf{C}^k$ such that $V\cap(P_1\times P_2)$ has proper
projection onto $P_1$ and $V\cap(\{a\}\times P_2)$ is a one element set. Then, by definition
of $\mathcal{S}(V),$ generic fibers in $V\cap(P_1\times P_2)$ have at least $2$ elements.
Since $\Phi_L$ preserves the cardinality of generic fibers in $V,$ then it also preserves the
cardinality of generic fibers in $V\cap(P_1\times P_2).$ This implies that the fiber in
$\Phi_L(V\cap(P_1\times P_2))$ over $a$ is one-element, and generic fibers in
$\Phi_L(V\cap(P_1\times P_2))$ over $P_1$ have at least $2$ elements. Hence,
$\Phi_L(a,b)\in\mathcal{S}(\Phi_L(V\cap(P_1\times P_2)))\subset\mathcal{S}(\Phi_L(V)).$}

Finally, for any analytic subset $X$ of an open set
$\tilde{U}\subset\mathbf{C}^m$ let $X_{(k)}\subset\tilde{U}$
denote the union of all analytic irreducible components of $X$ of dimension
$k.$
%*************************************************************************************************
\section{Approximation}\label{proofs}
%*****************************************************
%****************************************************************
\subsection{The main theorem}\label{mainIprof}
\begin{theorem}\label{main}Let $U$ be an open subset of
$\mathbf{C}^n$ and let $F:U\rightarrow\mathbf{C}^{\hat{m}}$ be a
holomorphic map that satisfies a system of equations $Q(F(x))=0$
for $x\in U,$ where
$Q:\mathbf{C}^{\hat{m}}\rightarrow\mathbf{C}^q$ is a polynomial
map. Then for every $x_0\in U $ there are an open neighborhood
$U_0\subset U$ and a sequence
$\{F^{\nu}:U_0\rightarrow\mathbf{C}^{\hat{m}}\}$ of Nash maps
converging uniformly to $F|_{U_0}$ such that $Q(F^{\nu}(x))=0$ for
every $x\in U_0$ and $\nu\in\mathbf{N}.$
\end{theorem}
As said in Section \ref{introduction}, this theorem is known. The proof given in the present
paper is simpler than the previous ones and it allows us to design an
algorithm for computing the approximating maps.

The proof will be divided into two
parts. In part 1, we shall show how to reduce the problem to the case (c1)
specified as follows:\vspace*{2mm}\\
(c1) $q=1,$ and
$Q:\mathbf{C}^{\hat{m}}\approx\mathbf{C}^m_u\times\mathbf{C}_z\rightarrow\mathbf{C}$ is a
monic polynomial in $z$ such\linebreak \hspace*{7mm}that $R_Q(f_1,\ldots,f_m)\neq 0,$ where
$R_Q\in\mathbf{C}[u]$ is the discriminant of $Q$ {as a\linebreak \hspace*{7mm}polynomial in
$z,$} whereas
$f_1,\ldots,f_m$ are the first $m$ components of $F.$\vspace*{2mm}\\
In part 2, we shall show that given $Q, F$ (as in Theorem \ref{main}) such that (c1) holds,
we can produce a holomorphic map $g$ depending on $n-1$ variables and a polynomial
map $T$ with $T\circ g=0$ such that if $g$ can be
locally approximated by Nash maps $g^{\nu}$ with $T\circ g^{\nu}=0,$ then $F$ can be locally
approximated by Nash maps $F^{\nu}$ with $Q\circ F^{\nu}=0.$ For $n=1,$ $g$ will be a constant vector and we shall
take $g^{\nu}=g.$ Once parts 1, 2 are completed, the proof of Theorem \ref{main} will be completed as well (by induction on $n$).

The following lemma will be useful in part 1 of the proof. Let $U$ be a domain in $\mathbf{C}^n.$
(For the notion of an optimal polynomial see Section \ref{setswithprop}.)

\begin{lemma}\label{reducone}Assume we are given:\vspace*{2mm}\\
(1) a holomorphic map $(f_1,\ldots,f_m,f_{m+1},\ldots,f_{m+s}):U\rightarrow V,$
where $V\subset\mathbf{C}^m\times\mathbf{C}^s$ is an algebraic variety of pure dimension $m$ with proper
projection onto $\mathbf{C}^m,$\vspace*{2mm}\\
(2) a $\mathbf{C}$--linear map $L:\mathbf{C}^s\rightarrow\mathbf{C},$ $L\neq 0,$
such that the generic fiber in $V$ over $\mathbf{C}^m$ and the generic fiber in $\Phi_L(V)\subset\mathbf{C}^m\times\mathbf{C}$ over $\mathbf{C}^m$ have the same cardinality, where $\Phi_L:\mathbf{C}^m_u\times\mathbf{C}^s_v\rightarrow\mathbf{C}^m_u\times\mathbf{C}_z$ is defined by $\Phi_L(u,v)=(u,L(v)).$\vspace*{2mm}\\
Assume that $R_L(f_1,\ldots,f_m)\neq 0,$ where $R_L\in\mathbf{C}[u]$ is the discriminant of the optimal polynomial
$P_L\in(\mathbf{C}[u])[z]$ with $P_L^{-1}(0)=\Phi_L(V).$
Then for all sequences
$\{f_1^{\nu}\},\ldots,\{f_m^{\nu}\},\{\tilde{f}^{\nu}\}$
of functions, holomorphic on $U,$ converging locally uniformly to
$f_1,\ldots,f_m,L(f_{m+1},\ldots,f_{m+s}),$ respectively, such that\vspace*{2mm}\\
$\mathrm{(3.1)}$\hspace*{20mm}$P_L(f_1^{\nu},\ldots,f_m^{\nu},\tilde{f}^{\nu})=0, \mbox{ for almost all }\nu\in\mathbf{N},$\vspace*{2mm}\\
the following holds. There
exist sequences $\{f_{m+1}^{\nu}\},\ldots,\{f_{m+s}^{\nu}\}$ of functions, holomorphic on $U,$
converging locally uniformly to $f_{m+1},\ldots,f_{m+s},$ respectively, such that the
image of the map
$(f^{\nu}_1,\ldots,f^{\nu}_m,f^{\nu}_{m+1},\ldots,f^{\nu}_{m+s})$
is contained in $V$ for almost all $\nu\in\mathbf{N}.$
\end{lemma}
\textit{Proof of Lemma \ref{reducone}.} Recall that $\Phi_L(V)$ is algebraic (cf. Section \ref{setswithprop}). Since $V$ has proper projection onto $\mathbf{C}^m,$
$\Phi_L(V)$ also has proper projection onto $\mathbf{C}^m.$

For any holomorphic
mapping $H:E\rightarrow\mathbf{C}^m,$ where $E$ is an open subset
of $\mathbf{C}^n,$ and any algebraic subvariety $X$ of
$\mathbf{C}^m\times\mathbf{C}^s,$ consider the analytic set
$$\mathcal{V}(X,H)=\{(x,v)\in E\times\mathbf{C}^s:(H(x),v)\in X\}.$$
If $X$ has pure dimension $m$ and proper projection onto $\mathbf{C}^m,$
then $\mathcal{V}(X,H)$ has pure dimension $n$ and proper projection onto $E.$
Indeed, fix any compact $K\subset E.$ Then $H(K)\subset\mathbf{C}^m$ is also compact. Since $X$ has proper projection
onto $\mathbf{C}^m,$ the fibers in $X$ over
$H(K)$ are uniformly bounded. Consequently, by definition of $\mathcal{V}(X,H),$ the fibers
in $\mathcal{V}(X,H)$ over $K$ are also uniformly bounded. Thus $\mathcal{V}(X,H)$
has proper projection onto $E.$

Now $\mathcal{V}(X,H)$ has dimension $n$ because its projection onto { the open set $E\subset
\mathbf{C}^n$ is proper and surjective. The latter is true because the projection of $X$ onto
$\mathbf{C}^m$ is surjective, which in turn follows from the fact that $X$ has dimension $m$
and proper projection onto $\mathbf{C}^m$}. It remains to check that $\mathcal{V}(X,H)$ is of
pure dimension (i.e. none of its analytic irreducible components has dimension strictly
smaller than $n$). Suppose that this is not true. Then there is an open polydisc $B\times C$
{relatively compact in} $E\times\mathbf{C}^s$ such that $\mathcal{V}(X,H)\cap(B\times
C)\neq\emptyset,$ $\mathcal{V}(X,H)\cap(\overline{B}\times
\partial C)=\emptyset$ and $\mathrm{dim}(\mathcal{V}(X,H)\cap(B\times C))<n.$ Then, by
definition of $\mathcal{V}(X,H),$ we have $X\cap(H({B})\times C)\neq\emptyset$ and
$X\cap(H(\overline{B})\times\partial C)=\emptyset.$ {This implies that for every compact
$K\subset H(B)$ we have $(K\times\mathbf{C}^s)\cap X\cap(H(B)\times
C)=(K\times\overline{C})\cap X,$ hence $(K\times\mathbf{C}^s)\cap X\cap(H(B)\times C)$ is a
compact set, so the projection} {of $X\cap(H({B})\times C)$ onto $H(B)$ is proper.} Since $X$
has pure dimension $m$ and $H(B)\subset\mathbf{C}^m,$ the latter projection is surjective.
Consequently, by definition of $\mathcal{V}(X,H),$ the projection of
$\mathcal{V}(X,H)\cap(B\times C)$ onto $B$ is surjective which contradicts the fact that
$\mathrm{dim}(\mathcal{V}(X,H)\cap(B\times C))<n.$

Next put $\Psi_L(x,v)=(x,L(v))$ for any $x\in\mathbf{C}^n, v\in\mathbf{C}^s.$ Assume the
notation of Section \ref{setswithprop}. Then we have the following
{\begin{lemma}\label{funnyobvious} Let $Z\subset E \times\mathbf{C}^s$ be an analytic subset
of pure dimension $n$ with proper projection onto a domain $E\subset\mathbf{C}^n$ such that
the generic fiber in $Z$ and the generic fiber in ${\Psi_{L}(Z)}$ over $E$ have the same
cardinality. Then, for every analytic irreducible component $\Sigma$ of $\Psi_{L}(Z)$ there
exists an analytic irreducible component $\Gamma$ of $Z$ with $\Psi_{L}(\Gamma)=\Sigma$ such
that the generic fiber in ${\Gamma}$ and the generic fiber in ${\Sigma}$ over $E$ have the
same cardinality.
\end{lemma}}
\noindent\textit{Proof of {Lemma} \ref{funnyobvious}. }Let $\bigcup_{j=1}^lZ_j$ be the decomposition of $Z$ into
pairwise distinct analytic irreducible components. Then every $\Psi_L(Z_j)$ is irreducible. Let,
$\mu,\nu,\mu_j,\nu_j$ denote the cardinalities of the generic fibers in $Z,\Psi_L(Z),$ $Z_j,$ $\Psi_L(Z_j)$ over
$E,$ respectively. Clearly, $\mu=\mu_1+\cdots+\mu_l$ and $\mu_j\geq\nu_j$ for every $j.$ Since
$\Psi_L(Z)=\bigcup_{j=1}^l\Psi_L(Z_j),$ we have $\nu_1+\cdots+\nu_l\geq\nu.$ Now, by the hypothesis, $\nu=\mu.$
Hence, $\nu_j=\mu_j$ for every $j,$ as required.\qed\vspace*{2mm}

The lemma allows us to complete the proof of Lemma \ref{reducone}. Put
$\tilde{F}=(f_1,\ldots,f_m),$ $\tilde{F}^{\nu}=(f_1^{\nu},\ldots,f_m^{\nu}),$
$G=(f_{m+1},\ldots,f_{m+s}).$ {Note that, by the second and the third paragraph of the proof,
$\mathcal{V}(V,\tilde{F})$ and $\mathcal{V}(V,\tilde{F}^{\nu})$ have pure dimension $n$ and
proper projection onto $U.$}

{Next observe that the cardinalities of the generic fibers in
$\Psi_{L}(\mathcal{V}(V,\tilde{F})),$ $\mathcal{V}(V,\tilde{F}),$
$\Psi_{L}(\mathcal{V}(V,\tilde{F}^{\nu}))$ and in $\mathcal{V}(V,\tilde{F}^{\nu})$ over $U$
are equal for large $\nu.$} {Indeed, since $R_L\circ\tilde{F}\neq 0,$ we also have
$R_L\circ\tilde{F}^{\nu}\neq 0,$ for large $\nu,$ and therefore the cardinalities of the
generic fibers in $\mathcal{V}(V,\tilde{F})$ and in $\mathcal{V}(V,\tilde{F}^{\nu})$ over $U$
are equal. Now, by assumption (2) of Lemma \ref{reducone} and by  the fact that
$R_L\circ\tilde{F}\neq 0,$ the cardinalities of the generic fibers in
$\Psi_{L}(\mathcal{V}(V,\tilde{F}))$ and in $\mathcal{V}(V,\tilde{F})$ over $U$ are equal.
Clearly, the latter is also true for $\Psi_{L}(\mathcal{V}(V,\tilde{F}^{\nu}))$ and
$\mathcal{V}(V,\tilde{F}^{\nu})$ for $\nu$ large enough.}

{In view of $\mathrm{(3.1)}$ we have
$\mathrm{graph}(\tilde{f}^{\nu})\subset\Psi_{L}(\mathcal{V}(V,\tilde{F}^{\nu})).$ Hence, by
Lemma \ref{funnyobvious} with $Z=\mathcal{V}(V,\tilde{F}^{\nu})$, there is an analytic
irreducible component $\Gamma^{\nu}$ of $\mathcal{V}(V,\tilde{F}^{\nu})$ with
$\Psi_L(\Gamma^{\nu})=\mathrm{graph}(\tilde{f}^{\nu})$ such that the generic fiber in
$\Gamma^{\nu}$ over $U$ is one-element.} Therefore, by (1), (2), (3) of Section
\ref{setswithprop}, $\mathcal{S}(\Gamma^{\nu})=\emptyset$ and
$\Gamma^{\nu}=\mathrm{graph}(G^{\nu})$ for some holomorphic map
$G^{\nu}:U\rightarrow\mathbf{C}^s.$

Let us show that $\{G^{\nu}\}$ converges to $G$ locally uniformly. {This will be done in
three steps.} First, we will show that $\{G^{\nu}\}$ is locally uniformly bounded. Next we
will prove that $\{G^{\nu}\}$ converges to $G$ pointwise. {Finally, we will recall the proof
of the fact that the local uniform convergence follows from the first two steps (cf. also
\cite{WV}).}

Observe that $\{G^{\nu}\}$ is locally uniformly bounded.
Indeed, by $\Gamma^{\nu}\subset\mathcal{V}(V,\tilde{F}^{\nu}),$ the image of $(\tilde{F}^{\nu},G^{\nu})$
is contained in $V.$ Now, on one hand, $\{\tilde{F}^{\nu}\}$ is locally uniformly bounded (because it
is a convergent sequence of holomorphic functions). On the other hand,
$G^{\nu}(x)$ is contained in the fiber of $V$ over $\tilde{F}^{\nu}(x),$ for every $x\in U.$
Since $V$ has proper projection onto $\mathbf{C}^m,$ the sequence $\{G^{\nu}\}$ is also locally
uniformly bounded.

{Now suppose that there is $x_0\in U$ such that $\{G^{\nu}(x_0)\}$ does not converge to $G(x_0).$
Since $\{G^{\nu}(x_0)\}$ is bounded, there exists a subsequence $\{G^{\nu_i}(x_0)\}$ converging
to some $b\neq G(x_0)$.}
Let $K\subset U$ be a compact ball {with $\mathrm{int}(K)\neq\emptyset,$ $x_0\in K$.}
By Montel's Theorem (cf. \cite{KaK}, p. 17) and {by the choice of $x_0$}
{we may assume that $\{G^{\nu_i}\}$ converges uniformly on $K$ to some $\bar{G}\neq G|_K$}
with $\mathrm{graph}(\bar{G})\subset\mathcal{V}(V,\tilde{F}).$
{Now,
$\Psi_L(\mathcal{V}(V,\tilde{F}))\cap(\{x\}\times\mathbf{C})$ and
$\mathcal{V}(V,\tilde{F})\cap(\{x\}\times\mathbf{C}^s)$ have equal
cardinality for generic $x\in U,$ so we can pick $x=x_1$
for which these sets do have equal cardinality and moreover, $\bar{G}(x_1)\neq G(x_1)$.
On one hand, $$\Psi_L|_{\mathcal{V}(V,\tilde{F})\cap(\{x_1\}\times\mathbf{C}^s)}:\mathcal{V}(V,\tilde{F})\cap(\{x_1\}\times\mathbf{C}^s)\rightarrow \Psi_L(\mathcal{V}(V,\tilde{F}))\cap(\{x_1\}\times\mathbf{C})$$ is injective,
because it is surjective and its domain and range are finite and have equal cardinality.
On the other hand,
$L\circ\bar{G}=L\circ G$ because $L\circ G^{\nu_i}=\tilde{f}^{\nu_i},$ so
$L(\bar{G}(x_1))=L(G(x_1))$ which contradicts the injectivity of
$\Psi_L|_{\mathcal{V}(V,\tilde{F})\cap(\{x_1\}\times\mathbf{C}^s)}$.
Thus $\{G^{\nu}\}$ indeed converges to $G$ pointwise.}

{If $\{G^{\nu}\}$ does not converge to $G$ locally uniformly, then there are $\varepsilon>0,$
a subsequence $\{G^{\nu_j}\}$ and a sequence $\{x_{\nu_j}\}$ contained in some compact subset
$K\subset U$ such that $||G^{\nu_j}(x_{\nu_j})-G(x_{\nu_j})||>\varepsilon$ for every $j.$
Then from $\{G^{\nu_j}\}$ one cannot choose a subsequence converging uniformly on $K$ which,
in view of the facts that $\{G^{\nu_j}\}$ is locally uniformly bounded and pointwise
convergent, contradicts Montel's Theorem.}

Now we can define $f_{m+i}^{\nu}$ to be the $i$'th component
of the map $G^{\nu},$ for $i=1,\dots,s.$\qed\vspace*{2mm}\\
\textit{Proof of Theorem \ref{main}. } \textit{Part 1.} First, following \cite{B6}, we carry
out some preparations. Since the problem is local, it is sufficient to consider the case
where $U$ is connected. Then $F(U)$ is contained in an irreducible component of the algebraic
variety $V=\{y\in\mathbf{C}^{\hat{m}}: Q(y)=0\}$, so we may assume that $V$ is of pure
dimension, say $m.$ We may also assume that
$V\subset\mathbf{C}^{\hat{m}}\thickapprox\mathbf{C}^m\times\mathbf{C}^{s}$ is with proper
projection onto $\mathbf{C}^m.$ Indeed, for a generic $\mathbf{C}$-linear isomorphism
$J:\mathbf{C}^{m+s}\rightarrow\mathbf{C}^{m+s}$ the image $J(V)$ is with proper projection
onto $\mathbf{C}^m.$ Thus if there exists a sequence $H^{\nu}:U_0\rightarrow J(V)$ of Nash
mappings converging to $J\circ F|_{U_0}$ then the sequence $\{J^{-1}\circ H^{\nu}\}$
satisfies the assertion of the {theorem}.

Now the problem will be reduced to the case where $V$ is a hypersurface (see \cite{B6},
compare also \cite{vD}, p. 394). Let $L:\mathbf{C}^s\rightarrow\mathbf{C}$ be any
$\mathbf{C}$--linear form such that the generic fibers in $\Phi_L(V)$ and in $V$ over
$\mathbf{C}^m$ have the same cardinality, where
$\Phi_{L}:\mathbf{C}^m_u\times\mathbf{C}^s_v\rightarrow\mathbf{C}^m_u\times\mathbf{C}_z$ is
defined by $\Phi_L(u,v)=(u,L(v)).$ (The generic $L$ has this property. It is clear that
$\Phi_L(V)$ is an algebraic hypersurface of $\mathbf{C}^m_u\times\mathbf{C}_z$ with proper
projection onto $\mathbf{C}^m_u.$) Let $P_L\in\mathbf(\mathbf{C}[u])[z]$ be the monic
polynomial in $z$ with nonzero discriminant $R_L\in\mathbf{C}[u]$ and with
$P_L^{-1}(0)=\Phi_L(V).$ (Such $P_L$ is called the optimal polynomial for $\Phi_L(V)$, {cf.
Section \ref{setswithprop}}.)

{Let} $f_1,\ldots,f_m,f_{m+1},\ldots,f_{m+s}$ denote the coordinates of $F.$ {Observe} that
we may assume $R_L(f_1,\ldots,f_m)\neq 0.$ Indeed, otherwise we return to the very beginning
of part 1 and repeat the whole construction with $V$ replaced by $V\cap \{R_L=0\}.$ Since the
latter variety is of pure dimension $m-1,$ we eventually reach the required condition. (The
image of the projection of $V\cap \{R_L=0\}$ onto $\mathbf{C}^m$ is $\{R_L=0\}.$ For more
information on discriminant hypersurfaces cf.~\cite{LMC}.)

If, in a neighborhood of some point, there exist Nash approximations\linebreak
$f_1^{\nu},\ldots,f_m^{\nu},\tilde{f}^{\nu}$ of $f_1,\ldots,f_m,$ $L(f_{m+1},\ldots,f_{m+s}),$ respectively, such
that $$P_L(f_1^{\nu},\ldots,f_m^{\nu},\tilde{f}^{\nu})=0,$$ then there exists a Nash approximation $F^{\nu}$ of
$F$ with image contained in $V$. This is because, by Lemma \ref{reducone}, there exist holomorphic functions
$f^{\nu}_{m+1},\ldots,f^{\nu}_{m+s}$ approximating $f_{m+1},\ldots,f_{m+s},$ respectively, such that the image of
the map $(f^{\nu}_1,\ldots,f^{\nu}_m,f^{\nu}_{m+1},\ldots,f^{\nu}_{m+s})$ is contained in $V.$ Hence, it is
sufficient to observe that $f^{\nu}_{m+1},\ldots,f^{\nu}_{m+s}$ must be Nash functions and next to define
$F^{\nu}=(f^{\nu}_1,\ldots,f^{\nu}_m,f^{\nu}_{m+1},\ldots,f^{\nu}_{m+s}).$ But these functions are Nash because
$$T_i(f^{\nu}_1,\ldots,f^{\nu}_m,f^{\nu}_{m+i})=0,$$ where $T_i\in(\mathbf{C}[u])[z]$ is the optimal polynomial
for $\pi_i(V),$ where
$\pi_i:\mathbf{C}^m\times\mathbf{C}_1\times\cdots\times\mathbf{C}_s\rightarrow\mathbf{C}^m\times\mathbf{C}_i$ is
the natural projection, so $\mathrm{graph}({f}^{\nu}_{m+i})$ is contained in a Nash hypersurface
$\{(x,z):T_i(f^{\nu}_1(x),\ldots,f^{\nu}_m(x),z)=0\}$, for $i=1,\ldots,s$ ({cf. Corollary \ref{functHypNash}}).

Since $R_L(f_1,\ldots,f_m)\neq 0,$ we have completed part 1, i.e. we have reduced the problem
to the case where it is sufficient to consider  $Q=P_L$ and approximate the map
$(f_1,\ldots,f_m,$ $L(f_{m+1},\ldots,f_{m+s})).$\vspace*{2mm}\\
%*************************************************************************************
\textit{Part 2.} Assume that $Q, F$ satisfy (c1) and fix $x_0\in U.$ Without loss of generality
we assume that $x_0$ is the origin in $\mathbf{C}^n.$ We shall produce a holomorphic map $g$ depending on $n-1$ variables and a polynomial map $T$ with $T\circ g=0$ such that if $g$ can be
locally approximated by Nash maps $g^{\nu}$ with $T\circ g^{\nu}=0,$ then, in some neighborhood of
the origin, $F$ can be approximated by Nash maps $F^{\nu}$ with $Q\circ F^{\nu}=0.$ (For $n=1,$ $g$ will be
a constant vector, and then we take $g^{\nu}=g.$) Once this is completed, the proof
of Theorem \ref{main} will be completed as well (by induction on $n$).

The reduction from $n$ to $n-1$ will be carried out in a similar way to \cite{Ar1}
(by applying the Weierstrass Preparation Theorem). The Tougeron
Implicit Functions Theorem, which often appears in the context, is
replaced here, as in \cite{vD}, by the following lemma.

\begin{lemma}\label{vddlem}\em{(\cite{vD}, p. 393). }\em
{For all positive integers $d, M$, }there is $\varepsilon>0$ such that {for every compact polydisc
$K\subset\mathbf{C}^n$ centered at the origin the following\linebreak hold. For every
$A=a_0z^d+a_1z^{d-1}+\cdots+a_d\in\mathcal{O}(K)[z]$ whose coefficients $a_i, i=0,\ldots,d,$ satisfy $\sup_{x\in K
}|a_i(x)|<M,$ and for all $\alpha, c\in\mathcal{O}(K)$ with $\sup_{x\in K}|\alpha(x)|<M,$ $\sup_{x\in
K}|c(x)|<\varepsilon$} and $A(\alpha)=c\cdot(\frac{\partial A}{\partial z}(\alpha))^2$ there is
{$b\in\mathcal{O}(K)$} such that $A(b)=0,$ and $\mathrm{ord}_0(b-\alpha)\geq\mathrm{ord}_0c\cdot\frac{\partial
A}{\partial z}(\alpha),$ and {$$\sup_{x\in K}|b(x)-\alpha(x)|\leq 2\sup_{x\in K}|c(x)(\frac{\partial A}{\partial
z}(\alpha))(x)|.$$}
\end{lemma}
\begin{remark}\label{vdrema}\em The constant $\varepsilon$ in Lemma \ref{vddlem}
depends only on $M,d$ and can be effectively calculated (see
\cite{vD}, the proofs of Lemma 1.5 and Lemma 1.6, pp 392-393).
\end{remark}

It is more convenient for us to keep the notation of part 1. Hence, instead of $Q$ we write $P_L.$
The discriminant of $P_L$ will be denoted by $R_L,$ and the components of $F$ will be denoted by
$f_1,\ldots,f_m,\tilde{f}.$ Set $\tilde{F}=(f_1,\ldots,f_m).$

Let us turn to constructing $g, T.$ Since $R_L\circ\tilde{F}\neq 0$ and
$P_L(\tilde{F},\tilde{f})=0,$ we have $\frac{\partial
P_L}{\partial z}(\tilde{F},\tilde{f})\neq 0.$ Let $x', x=(x',x_n)$ be the tuples of the coordinates
in $\mathbf{C}^{n-1}, \mathbf{C}^n,$ respectively.
We may assume that $\frac{\partial
P_L}{\partial z}(\tilde{F}(o,\cdot),\tilde{f}(o,\cdot))$ has a zero of finite order, say $d,$
at $x_n=0,$ where $o$ denotes the origin in $\mathbf{C}^{n-1}.$ (Otherwise we apply a linear change of the coordinates in $\mathbf{C}^n.$)

By the Weierstrass
Preparation Theorem, $\frac{\partial P_L}{\partial
z}(\tilde{F}(x),\tilde{f}(x))=\hat{H}(x)W(x)$ in some neighborhood of the origin,
where $\hat{H}$ is a non-vanishing holomorphic function, and
$$W(x)=x_n^{d}+x_n^{d-1}a_{1}(x')+\cdots+a_{d}(x')$$ is the Weierstrass polynomial. (For every $l=1,\ldots,d,$ the function $a_l$
is holomorphic in some neighborhood of $o,$ and $a_l(o)=0.$)

Dividing $f_j,\tilde{f}$ by $W^2,$ one obtains, for $j=1,\ldots,m$:\\
$$f_j(x)=H_j(x)W(x)^2+r_j(x),$$ $$\tilde{f}(x)=\tilde{H}(x)W(x)^2+\tilde{r}(x)$$
in some neighborhood of the origin, where $H_j,\tilde{H}$ are holomorphic functions,
and\vspace*{2mm}\\
$r_j(x)=x_n^{2d-1}b_{j,0}(x')+x_n^{2d-2}b_{j,1}(x')+\cdots+b_{j,2d-1}(x'),$\\
$\tilde{r}(x)=x_n^{2d-1}c_{0}(x')+x_n^{2d-2}c_{1}(x')+\cdots+c_{2d-1}(x')$\vspace*{2mm}\\
are polynomials with coefficients holomorphic in some neighborhood of $o.$

Replacing the coefficients
$$a_{1},\ldots,a_{d},b_{j,0},\ldots,b_{j,2d-1},
c_{0},\ldots,c_{2d-1},$$ for all $j,$ in $W,r_j,\tilde{r}$ by new
variables denoted by the same letters we obtain polynomials
$\omega,\rho_j,\tilde{\rho}$. Define:
$$\phi_j=\chi_j\omega^2+\rho_j,\mbox{ }\tilde\phi=\tilde\chi\omega^2+\tilde\rho,$$ for
$j=1,\ldots,m,$ where $\chi_j,\tilde\chi$ are new variables. Now divide
$P_L(\phi_1,\ldots,\phi_m,\tilde\phi)$ by $\omega^2$ (treated as a
polynomial in $x_n$ with polynomial coefficients) and divide
$\frac{\partial P_L}{\partial z }(\phi_1,\ldots,\phi_m,\tilde\phi)$
by $\omega$ to
obtain\vspace*{2mm}\\
$\mathrm{(3.2)}$\hspace*{9mm}$P_L(\phi_1,\ldots,\phi_m,\tilde\phi)=\tilde{W}\omega^2+x_n^{2d-1}T_1+x_n^{2d-2}T_2+\cdots+T_{2d},$\vspace*{2mm}\\
$\mathrm{(3.3)}$\hspace*{6mm}$\frac{\partial P_L}{\partial z
}(\phi_1,\ldots,\phi_m,\tilde\phi)=\tilde{S}\omega+x_n^{d-1}T_{2d+1}+x_n^{d-2}T_{2d+2}+\cdots+T_{3d},$\vspace*{2mm}\\
where $\tilde{W},\tilde{S},T_1,\ldots,T_{3d}$ are polynomials such
that $T_1,\ldots,T_{3d}$ depend only on the variables
$$a_{1},\ldots,a_{d},b_{j,0},\ldots,b_{j,2d-1},
c_{0},\ldots,c_{2d-1}.$$

Let
$$a_{1},\ldots,a_{d},b_{j,0},\ldots,b_{j,2d-1},
c_{0},\ldots,c_{2d-1},$$ for $j=1,\ldots,m,$ now denote the
holomorphic coefficients of $W, r_1,\ldots,r_m, \tilde{r}.$ The tuple consisting of all these coefficients
will be denoted by $g$, and this is the map announced in the first paragraph of part 2. Set $T=(T_1,\ldots,T_{3d}).$

By uniqueness of the Weierstrass Division Theorem we have $T\circ g=0.$ It remains to check
that if $g$ can be locally approximated by a Nash map $g^{\nu}$ with $T\circ g^{\nu}=0,$ then,
in some neighborhood of the origin, $F$ can be approximated by a Nash map $F^{\nu}$ with $P_L\circ F^{\nu}=0.$

Let $\{g^{\nu}\}$ be a sequence of Nash maps approximating
$g$ uniformly in some neighborhood of
$o$ with $T\circ g^{\nu}=0.$
The components of $g^{\nu}$ will be denoted by
$$a_{1}^{\nu},\ldots,a_{d}^{\nu},b_{j,0}^{\nu},\ldots,b_{j,2d-1}^{\nu},
c_{0}^{\nu},\ldots,c_{2d-1}^{\nu},$$ where $j=1,\ldots,m.$

Define:\vspace*{2mm}\\
$W_{\nu}(x)=x_n^{d}+x_n^{d-1}a_{1}^{\nu}(x')+\cdots+a_{d}^{\nu}(x'),$\\
$r_{j,\nu}(x)=x_n^{2d-1}b_{j,0}^{\nu}(x')+x_n^{2d-2}b_{j,1}^{\nu}(x')+\cdots+b_{j,2d-1}^{\nu}(x'),$\\
$\tilde{r}_{\nu}(x)=x_n^{2d-1}c_{0}^{\nu}(x')+x_n^{2d-2}c_{1}^{\nu}(x')+\cdots+c_{2d-1}^{\nu}(x'),$\vspace*{2mm}\\
and, for $j=1,\ldots,m,$ define
$$f_j^{\nu}(x)=H_j^{\nu}(x)(W_{\nu}(x))^2+r_{j,\nu}(x),$$ $$\bar{f}^{\nu}(x)=\tilde{H}^{\nu}(x)(W_{\nu}(x))^2+\tilde{r}_{\nu}(x).$$
Here $\{H^{\nu}_{j}\},\{\tilde{H}^{\nu}\}$ are any sequences of
polynomials converging uniformly to $H_{j},\tilde{H},$ respectively, in some
neighborhood of the origin.

Now it is easy to see that by $\mathrm{(3.2)},$ $\mathrm{(3.3)}$ and the way
$f^{\nu}_1,\ldots,f^{\nu}_m,\bar{f}^{\nu}$ are defined, there is a
neighborhood of $0\in\mathbf{C}^n$ in which, for sufficiently large $\nu,$ the
following holds:
$$P_L(f_1^{\nu},\ldots,f_m^{\nu},\bar{f}^{\nu})=C^{\nu}\cdot(\frac{\partial{P_L}}{\partial{z}}(f_1^{\nu},\ldots,f_m^{\nu},\bar{f}^{\nu}))^2,$$
where $\{C^{\nu}\}$ is a sequence of holomorphic functions
converging to zero {uniformly. Uniform convergence of $\{C^{\nu}\}$ requires a brief explanation.
First, there is a closed polydisc $\overline{E'}\times \overline{E''}\subset\mathbf{C}^{n-1}\times\mathbf{C}$ centered at $0$ such that $\{f_1^{\nu}\},\ldots,\{f_m^{\nu}\},\{\bar{f}^{\nu}\}$ converge uniformly on $\overline{E'}\times \overline{E''}$ and, for $\nu$ large enough,
$$(\overline{E'}\times\partial E'')\cap \{x:\frac{\partial{P_L}}{\partial{z}}(f_1^{\nu}(x),\ldots,f_m^{\nu}(x),\bar{f}^{\nu}(x))=0\}=\emptyset.$$ Then
$\{C^{\nu}\}$ converges to $0$ uniformly on $\overline{E'}\times\partial E''$ so, by the Maximum Principle,
it converges to $0$ uniformly on $E'\times E''.$}

{In view of the previous paragraph,} it suffices to apply Lemma
\ref{vddlem} with $A=P_L(f_1^{\nu},\ldots,f_m^{\nu},z),$
$\alpha=\bar{f}^{\nu}$ and $c=C^{\nu}$ (for sufficiently large
$\nu$) to obtain
$$P_L(f_1^{\nu},\ldots,f_m^{\nu},\tilde{f}^{\nu})=0,$$ in some neighborhood $U_0$ of the origin, where $\{\tilde{f}^{\nu}\}$ is a
sequence of holomorphic functions converging to $\tilde{f}$ in
$U_0.$ Since $f_1^{\nu},\ldots,f_m^{\nu}$ are Nash functions, $\tilde{f}^{\nu}$ is a Nash
function as well (because its graph is contained in a Nash hypersurface $\{(x,z):P_L(f_1^{\nu}(x),\ldots,f_m^{\nu}(x),z)=0\}$
{; cf. Corollary \ref{functHypNash}}).

Defining $F^{\nu}=(f_1^{\nu},\ldots,f_m^{\nu},\tilde{f}^{\nu})$ we complete part 2.\qed
%*****************************************************
%
\subsection{Algorithms}\label{examp}
\subsubsection{Computing with holomorphic functions}\label{model}

In Section \ref{examp}, we work under the assumption that given a complex number $z$ we can tell whether $z=0$ or not. The assumption is justified by the fact that in practice $\mathbf{C}$ is replaced by a computable subfield
$K$ such that for $z\in K$ one can decide whether $z=0$ or not.

In this subsection the model of computation for our algorithms is described. It is clear that not every
holomorphic function can constitute (a part of) input data. Only objects which can be encoded as finite sequences
of symbols can be considered. Therefore, we assume (cf. \cite{vH2}, \cite{vH4}) that every function $f$ depending
on $x_1,\ldots,x_n$ is given by a finite procedure $Expand_f$ which for every tuple
$(k_1,\ldots,k_n)\in\mathbf{N}^n$ returns the coefficient of the monomial $x_1^{k_1}\cdots x_n^{k_n}$ of the
Taylor expansion of $f$ around zero.  The input data corresponding to the function $f$ consist of the procedure
$Expand_f$, the size of a polydisc $U_f$ centered at zero and a constant $M_f$ such that: $\sup_{x\in
U_f}|f(x)|<M_f$ and the Taylor series of $f$ is convergent on (some domain containing) $U_f.$ Then we can control
accuracy of polynomial approximation of $f.$ More precisely, the Cauchy integral formula yields estimates
(majorants) for the coefficients of the Taylor expansion of $f.$ Using these estimates we can compute, for a given
polydisc $\tilde{U}$ {relatively compact in} $U_f$ centered at zero and for $\varepsilon>0,$ an integer $N$ such
that $\sup_{x\in\tilde{U}}|\tilde{f}(x)-f(x)|<\varepsilon,$ where $\tilde{f}$ is the Taylor polynomial of $f$ at
zero of order $N.$ (For more information on Taylor models, majorant series and accuracy control in various
settings see \cite{Hef}, \cite{MaBe}, \cite{MS}, \cite{Neh}, \cite{vH4}, \cite{vH5} and references therein.)

Observe that for functions represented in this way, it is possible to perform ring operations
and compute differentials (cf. \cite{vH1}, \cite{vH2} and references therein). Namely, having
input data for two functions $f,g$ one can recover $Expand_{f+g},$ $Expand_{f\cdot g},$
$Expand_{\frac{\partial f}{\partial x_i}},$ for every $i\in\{1,\ldots,n\},$ and
$Expand_{\frac{1}{f}},$ if $f(0)\neq 0.$ One can also set $U_{f+g}=U_f\cap U_g,$
$M_{f+g}=M_f+M_g$ and $U_{f\cdot g}=U_f\cap U_g,$ $M_{f\cdot g}=M_f\cdot M_g.$ Moreover, the
Cauchy integral formula implies $\sup_{x\in\theta U_f}|\frac{\partial f}{\partial
x_i}(x)|<\frac{M_f}{r_i(1-\theta)},$ for every $\theta\in(0,1),$ where $\theta U_f$ is the
polydisc centered at the origin of polyradius $(\theta r_1,\ldots,\theta r_n),$ where
$(r_1,\ldots,r_n)$ is the polyradius of $U_f.$ So one can define $U_{\frac{\partial
f}{\partial x_i}}=\theta U_f$ and $M_{\frac{\partial f}{\partial
x_i}}=\frac{M_f}{r_i(1-\theta)}$ for some $\theta.$ Finally, once we know the bounds for the
{partial} derivatives of $f,$ we can compute a polydisc $\tilde{U}\subset U_f$ such that
$|f(x)|>\frac{1}{2}|f(0)|,$ for every $x$ in the closure of $\tilde{U},$ provided that
$f(0)\neq 0.$ {Indeed, we can clearly compute a constant $M$ such that, after shrinking $U_f$
slightly, the derivatives $D_xf$ satisfy $\mathrm{sup}_{x\in U_f}||D_xf||<M.$ Then
$|f(x)-f(0)|\leq M||x||$ for $x\in U_f$, so $|f(x)|\geq |f(0)|-M||x||.$ Then it is sufficient
to choose the radius of} {$\tilde{U}$ so small that $|f(0)|-M||x||>\frac{|f(0)|}{2}$ for
every $x$ in the closure of $\tilde{U}.$} Therefore we can set $U_{\frac{1}{f}}=\tilde{U},$
$M_{\frac{1}{f}}=\frac{2}{|f(0)|}.$

Let us recall that in this model we have effective versions of the Weierstrass Preparation and Division Theorems
(which will be used in the algorithms). Since the former is a consequence of the latter (see \cite{GP}, p. 319),
it is sufficient to discuss the Division Theorem. Let $f, g$ be functions (depending on $x_1,\ldots,x_n$)
holomorphic in some polydiscs $U_f, U_g$ centered at the origin, such that $f(o,x_n)=x_n^du(x_n),$ for some
$d\in\mathbf{N},$ where $u$ is holomorphic, $u(0)\neq 0,$ and $o$ is the origin in $\mathbf{C}^{n-1}.$ Then,
according to the Division Theorem, there are a holomorphic function $q$ and a polynomial $r$ in $x_n$ with
$\mathrm{deg}(r)<d$ and with holomorphic coefficients such that, in some neighborhood of the origin, $g=qf+r.$ We
need to show that given data representing $f,g,$ we can recover the data representing $q,r.$ The existence of the
procedures $Expand_q,$ $Expand_r$ is a consequence of the proof of the Division Theorem (see \cite{GP}, pp.
318-319; such procedures are in fact written, cf. \cite{GP}, p. 544).

It remains to compute  $U_q, U_r, M_q, M_r.$ Since $f(o,x_n)=x^d_nu(x_n),$ $u(0)\neq 0,$
there are a polydisc $A\times
B\subset\mathbf{C}^{n-1}_{x_1,\ldots,x_{n-1}}\times\mathbf{C}_{x_n}$ centered at the origin
and a constant $c$ with $\inf_{\overline{A}\times\partial B}|f|>c.$ Such $A, B, c$ can be
computed by means of $Expand_f, U_f, M_f.$ {Indeed, as above, we can compute a disc
$B_t\subset\mathbf{C}$ (centered at zero of radius $t$) such that
$|u(x_n)|\geq\frac{|u(0)|}{2}$ on $B_t$ and, moreover, a polydisc
$A_s\subset\mathbf{C}^{n-1}$ (centered at zero of radius $s$) such that $A_s\times B_t$ is
relatively compact in $U_f$. Then $|f(o,x_n)|\geq t^d\frac{|u(0)|}{2}$ on $\partial B_t.$
Since we know the bounds for the partial derivatives of $f$ on $A_s\times B_t,$ we can
compute $M$ such that for any $\tilde{s}\leq s$ and any $(x,y)\in A_{\tilde{s}}\times\partial
B_t$ we have $|f(x,y)-f(o,y)|\leq M\tilde{s}.$ Then $|f(x,y)|\geq
t^d\frac{|u(0)|}{2}-M\tilde{s},$ so we can take $B=B_t,$ $A=A_{\tilde{s}},$
$c=t^d\frac{|u(0)|}{2}-M\tilde{s}$ for $\tilde{s}$ small enough.}

Observe that $\overline{A\times B}$ may be additionally assumed to be contained in $U_{f}\cap
U_g.$ Let $C\subset B\subset \mathbf{C}$ be another disc centered at zero of radius, say,
$\frac{1}{2}$ of the radius of $B,$ and let $U_q=U_r=A\times C.$ Then the Taylor series of
$q,r$ at zero are convergent on $U_q.$ Now, by the integral representations of $q,r$ (see
\cite{Lo}, p. 112), for $(x',x_n)\in A\times C,$ where $x'=(x_1,\ldots,x_{n-1}),$ we have
$$q(x',x_n)=\frac{1}{2\pi i}\int_{\partial B}\frac{g(x',s)}{f(x',s)}\frac{1}{s-x_n}ds,$$
$$r(x',x_n)=\frac{1}{2\pi i}\int_{\partial B}\frac{g(x',s)}{f(x',s)}\frac{f(x',s)-f(x',x_n)}{s-x_n}ds.$$
Since $|g|, |f|$ are bounded from above on $\overline{A\times B}$ by $M_g, M_f,$ and $|f|$ is bounded from below
on $A\times\partial B$ by $c$, and $|s-x_n|$ is bounded from below for $x_n\in C, s\in\partial B$ by $\frac{1}{2}$
of the radius of $B,$ we easily obtain $M_q=\frac{2M_g}{c}, M_r=\frac{4M_gM_f}{c}.$

For effective versions of the Division Theorem for special classes of functions we refer to \cite{vH6} and the
references therein.

{In such a general model as the one considered in the present paper we cannot effectively
decide whether a given holomorphic function equals zero identically. That is, the zero-test
is not computable. It is only semi-computable. Namely, we can compute the coefficients of the
Taylor expansion of higher and higher} {orders to check whether some of them are non-zero. If
we come across a non-zero one then the function is non-zero. But if the function equals zero,
then the procedure does not stop. (Here let us mention that the zero-test is computable for
certain subclasses of our class of functions. For details see \cite{vH3} and references
therein.)} {In particular, given a holomorphic map $f$ and an algebraic set $V,$ we cannot
decide whether the image of $f$ is contained in $V.$ On the other hand, we have the
following}
\begin{remark}\label{interalgebr}{Let $V_1,V_2\subset\mathbf{C}^m$ be algebraic subsets and
let $U$ be an open connected neighborhood of $0\in\mathbf{C}^n.$ Let
$f=(f_1,\ldots,f_m):U\rightarrow\mathbf{C}^m$ be a holomorphic map, where $f_1,\ldots,f_m$
are represented in the data structure introduced above. Assume $f(U)\nsubseteq V_1\cap V_2.$
Then we can effectively point at (at least one) $i\in\{1,2\}$ such that $f(U)\nsubseteq
V_i.$}
\end{remark}
{We shall use this fact in the sequel. To see that it is true, note that at least one of the
finite number of polynomials defining $V_1$ or of polynomials defining $V_2$ is a non-zero
function after being composed with $f.$ If we know that in a finite family of functions (at
least) one of the functions does not equal zero identically, then, clearly, detecting such a
function is effectively computable.}

In the next subsection, we will present a partial algorithm relying on the {zero-test} which, as mentioned above,
is not computable. This problem can be partially handled as follows. Given {a holomorphic function $\hat{R}$ and}
$\theta\in (0,1),$ for every $\varepsilon>0$ we can compute an integer $N(\theta,\varepsilon,\hat{R})$ such that
the following holds. If the coefficients of the monomials of order smaller than $N(\theta,\varepsilon,\hat{R})$
(of the Taylor expansion of $\hat{R}$) all equal zero, then $\sup_{\theta U_{\hat{R}}}|\hat{R}|<\varepsilon.$
(Recall that we have assumed that for $z\in\mathbf{C}$ we can decide whether $z=0$ or not.) The number
$N(\theta,\varepsilon,\hat{R})$ can be computed in the following way. Let
$\hat{R}=\sum_{\alpha\in\mathbf{N}^n}a_{\alpha}z^{\alpha}$ be the Taylor expansion of $\hat{R}.$ Suppose that for
some $N(\theta,\varepsilon,\hat{R})$ the coefficients $a_{\alpha}=0$ for every $|\alpha|\leq
N(\theta,\varepsilon,\hat{R}).$ Then by the Cauchy inequalities we have $$\sup_{\theta
U_{\hat{R}}}|\hat{R}|<M_{\hat{R}}\sum_{|\alpha|>N(\theta,\varepsilon,\hat{R})}\theta^{|\alpha|}.$$ Now for every
$\varepsilon>0$ we can compute how large $N(\theta,\varepsilon,\hat{R})$ {has to} be to ensure that the righthand
side of the last inequality is smaller than $\varepsilon.$

 This implies that after computing a finite number of the coefficients
of the Taylor expansion of $\hat{R}$ we can decide whether:\vspace*{2mm}\\
(i) $\hat{R}\neq 0$ (if at least one of these coefficients is non-zero)\\
(ii) $\sup_{\theta U_{\hat{R}}}|\hat{R}|<\varepsilon$ (else).\vspace*{2mm}\\
Therefore we assume $\hat{R}=0$ if the coefficients of the monomials of order smaller than
$N(\theta,\varepsilon,\hat{R})$ all equal zero (i.e. we assume $\hat{R}=0$ if (ii) holds),
where $\varepsilon$ is some fixed very small positive real number. It may happen that (ii)
holds and $\hat{R}\neq 0,$ and then the assumption $\hat{R}=0$ may give rise to a non-correct
output of the algorithm. A complete solution, obtained by refining the ideas presented below,
will be given in Section \ref{completfor1}, where the zero-test will be avoided.

{\subsubsection{Partial algorithm}\label{partialalg} This subsection} is devoted to a
recursive algorithm of Nash appro\-ximation of a holomorphic mapping $F:U\rightarrow
V\subset\mathbf{C}^{\hat{m}},$ where $U$ is a neighborhood of zero in $\mathbf{C}^n,$ for any
$n\in\mathbf{N},$ and $V$ is an algebraic variety. We will first state the algorithm and then
comment on the main steps. The algorithm relies on zero-test, which is semi-computable in our
model of computation. For this reason, for some input it may not return correct output data
{(and that is why we call it partial)}. This occurs when the image of $F$ is very close to
the singular locus $\mathrm{Sing}(V)$ of the target variety but not contained in
$\mathrm{Sing}(V).$ If every zero-test returns the negative answer, which is usually the
case, then the output we obtain is correct. In Section \ref{completfor1} we will show how to
avoid the zero-test, which will be done by refining the algorithm discussed below.

Components of $F$ will be represented as described
in Section \ref{model}.\vspace*{2mm}\\
\textbf{Input:} a positive integer $\nu,$ an algebraic variety $V,$ and a holomorphic mapping
$F=(f_1,\ldots,f_{\hat{m}}):U\rightarrow
V\subset\mathbf{C}^{\hat{m}},$ where $U$ is an open connected
neighborhood of $0\in\mathbf{C}_x^n$.\vspace*{2mm}\\
\textbf{Output:} $P_i^{\nu}(x,z_i)\in(\mathbf{C}[x])[z_i],$
$P_i^{\nu}\neq 0$ for $i=1,\ldots,\hat{m},$
with the following properties:\vspace*{1mm}\\
(a) $P_i^{\nu}(x,f_i^{\nu}(x))=0$ for every $x\in U_0,$ where
$F^{\nu}=(f^{\nu}_1,\ldots,f^{\nu}_{\hat{m}}):U_0\rightarrow V$ is
a holomorphic mapping such that $||F-F^{\nu}||_{U_0}<\frac{1}{\nu}$ and
$U_0$ is an open neighborhood of
$0\in\mathbf{C}^n$ independent of $\nu,$\\
(b) $P_i^{\nu}$ is a monic polynomial in $z_i$ of degree in $z_i$ bounded by a constant
independent of $\nu$.\vspace*{2mm}

First let us note that if $n=0$ or $\mathrm{dim}(V)=0,$ then the map $F$ is constant and $F^{\nu}=F.$ The
algorithm will reduce the problem to the case $n\cdot\mathrm{dim}(V)=0$  by recursive calls.

Before discussing the problem in general, let us look at the special case (in which no recursive calls are
necessary) where $V$ is of pure dimension $m$ for some $m>0$ and $F(0)\in\mathrm{Reg}(V).$ Then, after a generic
linear change of the coordinates, $V\subset\mathbf{C}^{\hat{m}}\approx\mathbf{C}^m\times\mathbf{C}^s$ has proper
projection onto $\mathbf{C}^m.$ Moreover, generic fibers in $V$ and in ${\Phi(V)}$ over $\mathbf{C}^m$ have the
same cardinalities, where $\Phi(z_1,\ldots,z_m,z_{m+1},\ldots,z_{m+s})$ $=(z_1,\ldots,z_m,z_{m+1}).$ Furthermore,
the fiber in $\Phi(V)$ over $(f_1,\ldots,f_m)(0)\in\mathbf{C}^m$ has the same number of elements as the generic
fiber in $\Phi(V)$ over $\mathbf{C}^m.$

Now calculate the optimal polynomial
$P\in(\mathbf{C}[z_1,\ldots,z_m])[z_{m+1}]$ describing
$\Phi(V)\subset\mathbf{C}^m\times\mathbf{C}$ and the discriminant $R$ of $P.$ By the previous paragraph, $R(f_1,\ldots,f_m)$ and
$\frac{\partial P}{\partial z_{m+1}}(f_1,\ldots,f_{m+1})$  do not vanish in some
neighborhood of $0\in\mathbf{C}^n.$ {Moreover, $P(f_1,\ldots,f_{m+1})=0,$
so for Taylor polynomials $f_1^{\nu},\ldots,f_{m}^{\nu}, \bar{f}_{m+1}^{\nu}$
close to $f_1,\ldots,f_m, f_{m+1},$ respectively, we have a holomorphic function
$$C^{\nu}(x)=\frac{P(f_1^{\nu}(x),\ldots,f_m^{\nu}(x),\bar{f}_{m+1}^{\nu}(x))}{(\frac{\partial P}{\partial z_{m+1}}(f^{\nu}_1(x),\ldots,f_m^{\nu}(x),\bar{f}_{m+1}^{\nu}(x)))^2}$$
such that $|C^{\nu}(x)|$ is small.} Consequently, by {Lemma \ref{vddlem}}, there is a Nash function
$f_{m+1}^{\nu}$ close to $f_{m+1}$ such that $P(f_1^{\nu},\ldots,f_{m+1}^{\nu})=0.$ Using {Remark \ref{vdrema}} we
can estimate how close $f_1^{\nu},\ldots, f_m^{\nu}$ to $f_1,\ldots,f_m$ {have to} be to ensure that
$f_{m+1}^{\nu}$ approximates $f_{m+1}$ with the required precision on some neighborhood of $0\in\mathbf{C}^n.$

For $f_1^{\nu},\ldots,f_m^{\nu}$ close enough to $f_1,\ldots,f_m$ we also have $R(f_1^{\nu},\ldots,f_m^{\nu})\neq
0.$ {Then, there are holomorphic functions $f_{m+2}^{\nu},\ldots,f_{m+s}^{\nu}$ such that the image of
$(f_1^{\nu},\ldots,f_m^{\nu}, f_{m+1}^{\nu},\ldots,f_{m+s}^{\nu})$ is contained in $V.$ Indeed, set
$H=(f_1^{\nu},\ldots,f_m^{\nu}),$ $z=(z_{m+1},\ldots,z_{m+s}),$ $\mathcal{V}(V,H)=\{(x,z):(H(x),z)\in V\},$
$L(z)=z_{m+1},$ and $\Psi_L(x,z)=(x,L(z))$ and confirm (previous paragraph) that
$$\Psi_L(\mathcal{V}(V,H))=\{(x,z_{m+1}):P(f_1^{\nu}(x),\ldots,f_m^{\nu}(x),z_{m+1})=0\}.$$
(For properties of $\mathcal{V}(V,H)$ cf. the discussion preceding Lemma \ref{funnyobvious}.)
Since $\Psi_L(\mathcal{V}(V,H))$ contains the graph of $f_{m+1}^{\nu},$ then, by Lemma
\ref{funnyobvious} with $Z=\mathcal{V}(V,H),$ we have that $\mathcal{V}(V,H)$ contains the
graph of $(f_{m+1}^{\nu},f_{m+2}^{\nu},\ldots,f_{m+s}^{\nu})$ for some holomorphic functions
$f_{m+2}^{\nu},\ldots,f_{m+s}^{\nu}.$}

Now there are $Q_{2},\ldots, Q_{s}\in(\mathbf{C}[z_1,\ldots,z_m])[z_{m+1}]$ such that
$$V\setminus R^{-1}(0)=\{R\neq 0, P=0, z_{m+j}R=Q_j,  j=2,\ldots,s\}.$$
(The polynomials $Q_2,\ldots, Q_s$ can be effectively
computed, cf. \cite{JS} p. 233, \cite{GH}, \cite{KP}.)
Consequently, the functions $f_{1}^{\nu},\ldots,f_{m+s}^{\nu}$ satisfy the equations
$$f_{m+j}^{\nu}R(f_1^{\nu},\ldots,f_{m}^{\nu})=Q_j(f_1^{\nu},\ldots,f_{m+1}^{\nu}),\mbox{ for } j=2,\ldots,s$$
(which in particular implies that $f_{m+2}^{\nu},\ldots,f_{m+s}^{\nu}$ are not only holomorphic but Nash as well).

Since the functions $f_1,\ldots,f_{m+s}$ also satisfy these equations, we can estimate (given $Q_j, R$) how close
$f_1^{\nu},\ldots,f_{m+1}^{\nu}$ to $f_1,\ldots,f_{m+1}$ {have to} be to ensure that
$f_{m+2}^{\nu},\ldots,f_{m+s}^{\nu}$ approximate $f_{m+2},\ldots,$ $f_{m+s}$ with the required precision. When we
know that the approximation exists, the polynomials $P_{i}^{\nu}$ can be computed as follows.

For $i=1,\ldots,m,$ set $P_i^{\nu}(x,z_i)=z_i-f_{i}^{\nu}(x).$ Next define
$$V^{\nu}=\{(x,z)\in\mathbf{C}^n_x\times\mathbf{C}^{m+s}_{z} : z\in
V, z_i=f_{i}^{\nu}(x), \mbox{ for } i=1,\ldots,m\},$$ where
$z=(z_1,\ldots,z_m,z_{m+1},\ldots,z_{m+s}).$ Finally, for $i=1,\ldots,s,$
take
$P^{\nu}_{m+i}\in$\linebreak $(\mathbf{C}[x])[z_{m+i}]$ to be the optimal
polynomial describing the image of the projection of $V^{\nu}$
to $\mathbf{C}^n_x\times\mathbf{C}_{z_{m+i}}.$ (This projection is proper so the image
is an algebraic variety.)

In the general case we have the following (partial) algorithm. {Let us recall that $V_{(k)}$
denotes the union of all $k$-dimensional irreducible components of $V.$}
\vspace*{2mm}\\
\textbf{Algorithm 1:}\vspace*{2mm}\\
\textbf{1.} Replace $V$ by $V_{(m)}$ such that
$F(U)\subset V_{(m)}.$\\
\textbf{2.} After a generic linear change of the coordinates in $\mathbf{C}^{\hat{m}}$
we have:\\
(x) $V\subset\mathbf{C}^{\hat{m}}\approx\mathbf{C}^m\times\mathbf{C}^s$ has proper
projection onto $\mathbf{C}^m,$\\
(y) generic fibers in $V$ and in ${\Phi(V)}$ over $\mathbf{C}^m$ have the same
cardinalities,\\
\hspace*{5.5mm}where $\Phi(z_1,\ldots,z_m,z_{m+1},\ldots,z_{m+s})=(z_1,\ldots,z_m,z_{m+1}).$\\
Apply such a change of the coordinates.\\
\textbf{3.} Calculate the optimal polynomial
$P\in(\mathbf{C}[z_1,\ldots,z_m])[z_{m+1}]$ describing
$\Phi(V)\subset\mathbf{C}^m\times\mathbf{C}.$ Calculate the discriminant $R\in\mathbf{C}[z_1,\ldots,z_m]$ of $P.$\\
\textbf{4.} If $R(f_1,\ldots,f_m)=0,$ then  go to Step 2 with $m,$ $s,$ $V$ replaced by
$m-1,$ $s+1,$ $V\cap\{R=0\},$ respectively. Otherwise $\frac{\partial P}{\partial z_{m+1}
}(f_{1},\ldots,f_{m},{f}_{m+1})\neq 0.$\\
\textbf{5.} If $\frac{\partial P}{\partial z_{m+1}}(f_1,\ldots,f_{m+1})(0)\neq 0,$ then
$F(0)\in\mathrm{Reg}(V)$ (this case has been discussed above). If $\frac{\partial P}{\partial z_{m+1}}(f_1,\ldots,f_{m+1})(0)= 0,$ then
apply a linear change of the coordinates in
$\mathbf{C}^n$ after which the following holds: $\frac{\partial
P}{\partial z_{m+1} }(f_1,\ldots,f_m, f_{m+1})(x)$$=\hat{H}(x)W(x)$
in some neighborhood of $0\in\mathbf{C}^n,$ where $\hat{H}$ is a
holomorphic function, $\hat{H}(0)\neq 0$ and $W$ is a monic
polynomial in $x_n$ with holomorphic coefficients
vanishing at $0\in\mathbf{C}^{n-1},$ depending on
$x'=(x_1,\ldots,x_{n-1}).$ Put $d=\mathrm{deg}(W).$\\
\textbf{6.} Divide $f_i$ by $(W)^2$ to obtain
$f_i=(W)^2H_i+r_i$ in some
neighborhood of $0\in\mathbf{C}^n,$ $i=1,\ldots,m+1.$ Here
$H_i$ are holomorphic functions and $r_i$ are
polynomials in
$x_n$ with holomorphic coefficients depending on $x'$ such that $\mathrm{deg}(r_i)<2d.$\\
\textbf{7.} Treating $H_i,$ $i=1,\ldots,m+1,$ and all the
coefficients of $W,r_1,\ldots,r_{m+1}$ as new variables
(except for the leading coefficient $1$ of
$W$) apply the division procedure for
polynomials to obtain:\vspace*{2mm}\\
$P(W^2H_1+r_1,\ldots,W^2H_{m+1}+r_{m+1})$
$=\tilde{W}W^2+x_n^{2d-1}T_1+x_n^{2d-2}T_2+\cdots+T_{2d},$\\
$\frac{\partial P}{\partial
z_{m+1}}(W^2H_1+r_1,\ldots,W^2H_{m+1}+r_{m+1})$\\
\hspace*{51mm}$=\tilde{S}W+x_n^{d-1}T_{2d+1}+x_n^{d-2}T_{2d+2}+\cdots+T_{3d}.$\\
Here $T_1,\ldots,T_{3d}$ are polynomials depending only on the variables standing for the coefficients of $W,$
$r_1,\ldots,r_{m+1}.$ Moreover, $T_1\circ g=\cdots=T_{3d}\circ g=0,$ where $g$ is the mapping whose components are
all these coefficients
(cf. Section \ref{mainIprof}).\\
\textbf{8.} If $n>1,$ then apply Algorithm 1 recursively with $\nu, F, V$ replaced by $\mu(\nu),
g,$$\{T_1=\cdots=T_{3d}=0\},$ respectively. (How to choose $\mu(\nu)$ is described below.) As a result, for every
component $c(x')$ of $g(x')$ one obtains a monic polynomial $Q_c^{\mu}(x',t_c)\in(\mathbf{C}[x'])[t_c]$ which put
in place of $P_i^{\nu}(x,z_i)$ satisfies (a) and (b) above with $\nu,$ $x,$ $z_i,$ $f_i^{\nu}, F, F^{\nu}$
replaced by $\mu,$ $x',$ $t_c,$ $c^{\mu},$ $g,$ $g^{\mu},$ respectively. Here, for every $c,$ $c^{\mu}$ is a Nash
function approximating $c,$ in some neighborhood of the origin in $\mathbf{C}^{n-1}$ such that the mapping
$g^{\mu}$ obtained by replacing every $c$\\
\\
of $g$ by $c^{\mu}$ satisfies $T_1\circ g^{\mu}=\cdots=T_{3d}\circ g^{\mu}=0.$ If $n=1,$ then $g$ is
constant, $g=g^{\mu},$ and one gets the $Q^{\mu}_c$'s immediately.\\
\textbf{9.} Approximate $H_i$ for $i=1,\ldots,m,$ by
polynomials $H_i^{\mu}.$ (Accuracy of this approximation is discussed
below.) Let $W_{\mu},r_{1,{\mu}},$
$\ldots,$ $r_{m+1,{\mu}}$ be the polynomials in $x_n$ defined by
replacing the coefficients of $W,$ $r_1,\ldots,r_{m+1}$ by
their Nash approxima\-tions (i.e. the components of $g^{\mu}$)
determined in Step 8. Using $Q_c^{\mu}$ (for all $c$) and
$H^{\mu}_i$ one can calculate $P_i^{\nu}\in(\mathbf{C}[x])[z_i],$
for $i=1,\ldots,m,$ satisfying (b) and (a) with
$f_i^{\nu}=H_i^{\mu}(W_{\mu})^2+r_{i,{\mu}}$ being the $i$'th
component of the mapping $F^{\nu}$ (whose last $\hat{m}-m$
components are determined by
$P_{m+1}^{\nu},\ldots,P_{\hat{m}}^{\nu}$ obtained in the next
step). {To calculate $P^{\nu}_i,$ for $i=1,\ldots,m,$ perform the ring operations
in the construction of $f_{i}^{\nu}$ in their implicit representation (see Lemma \ref{algNashfunct})}.\\
\textbf{10.} Put $V^{\nu}=\{(x,z)\in\mathbf{C}^n_x\times\mathbf{C}^{m+s}_{z} : z\in V, P^{\nu}_i(x,z_i)=0 \mbox{
for } i=1,\ldots,m\},$ where $z=(z_1,\ldots,z_m,z_{m+1},\ldots,z_{m+s}).$ For $i=1,\ldots,s,$ take
$P^{\nu}_{m+i}\in$\linebreak $(\mathbf{C}[x])[z_{m+i}]$ to be the optimal polynomial describing the image of the
projection of $V^{\nu}$ to $\mathbf{C}^n_x\times\mathbf{C}_{z_{m+i}}.$
\begin{remark}\label{alter}Observe that Step 4 can be replaced by:\vspace*{2mm}\\
{\em\textbf{4'.} If $\frac{\partial P}{\partial z_{m+1} }(f_{1},\ldots,f_{m+1})=0,$ then go
to Step 2 with $m,$ $s,$ $V$ replaced by $m-1,$ $s+1,$ $V\cap\{R=0\},$
respectively.}\vspace*{2mm}\\
Note that $\frac{\partial P}{\partial z_{m+1} }(f_{1},\ldots,f_{m+1})=0$ implies
$R(f_1,\ldots,f_m)=0$ but in general the converse is  false.
\end{remark}
%*************************************************************

In the remaining part of this section we comment on the main ideas of the algorithm.
Then we prove that one can effectively control the error of approximation of $F$
by $F^{\nu}.$ Finally, we show that the degree of $P_i^{\nu}$ in $z_i$ is bounded by a constant independent
of $\nu$ (i.e. independent of accuracy of approximation). The latter fact is important if one considers
efficiency of performing ring operations for functions $f_i^{\nu}$ in their implicit representation.

The aim of Step 1 is to reduce the problem to the case where $V$ is of pure dimension. Since
a generic algebraic variety is irreducible, $V$ is almost always purely dimensional, and Step
1 need not be performed. If $V$ is not purely dimensional, then we proceed as follows. Since
$U$ is connected, there is an integer $m$ such that $F(U)\subset V_{(m)}.$ {(Recall that
$V_{(m)}$ denotes the union of all $m$-dimensional irreducible components of $V.$)} To
compute $V_{(m)}$ we first decompose $V$ into equidimensional parts. (For algorithmic
equidimensional decomposition see \cite{DL}, p. 104, \cite{GP}, p. 258, \cite{JKSS},
\cite{JS}, \cite{Lec} and references therein, see also \cite{ScP}, \cite{SVW}.) Then it
remains to check for a fixed $m,$ whether $F(U)$ is contained in $V_{(m)}.$ This is
equivalent to testing whether $u_{m,l}\circ F$ is identically zero, for $l=1,\ldots, j_m,$
where $u_{m,1},\ldots,u_{m,j_m}$ are polynomials describing $V_{(m)}.$ (Here is the first
appearance of the zero-test discussed in Section \ref{model}.) If we have the extra knowledge
that $F(U)$ is not contained in $\mathrm{Sing}(V),$ then the problem is decidable because
there is precisely one $m$ such that $F(U)\subset V_{(m)}.$ (Then $m$ can be found in a
finite number of steps by excluding the other equidimensional components, {cf. Remark
\ref{interalgebr}}.)

As for Step 2, the set of linear isomorphisms
$\hat{J}:\mathbf{C}^m\times\mathbf{C}^s\rightarrow\mathbf{C}^m\times\mathbf{C}^s$ such that $\hat{J}(V)$
has proper projection onto $\mathbf{C}^m$ is dense and open in the set of all linear maps. {Moreover, given $\hat{J},$ one can test algorithmically whether
$\hat{J}(V)$ does have
proper projection onto $\mathbf{C}^m.$ Thus, we
can simply choose $\hat{J}$ at random with very high probability and then check whether it is suitable}. (This is discussed in detail in \cite{JS}, p. 235.)
Now if $V$ has proper projection onto $\mathbf{C}^m,$ then one can choose a linear form $L:\mathbf{C}^s\rightarrow\mathbf{C}$ such that
generic fibers in $\Phi_L(V)$ and in $V$ over $\mathbf{C}^m$
have the same finite number of elements, where $\Phi_L(z_1,\ldots,z_m,z_{m+1},\ldots,z_{m+s})=(z_1,\ldots,z_m,L(z_{m+1},\ldots,z_{m+s})).$
For the problem of choosing $L$ effectively, see \cite{JS}, p. 233 or \cite{GH}, Section 3.4.7, or \cite{KP}, Proposition 27. For simplicity of notation we assume  $L(z_{m+1},\ldots,z_{m+s})=z_{m+1},$ which can be achieved
by a linear change of coordinates in $\mathbf{C}^m\times\mathbf{C}^s$.

Once we have
$\hat{J}:\mathbf{C}^m\times\mathbf{C}^s\rightarrow\mathbf{C}^m\times\mathbf{C}^s$ such that
(x), (y) {in Step 2} are satisfied with $V$ replaced by $\hat{J}(V),$ we also replace $F$ by
$\hat{J}\circ F.$ One thing, which requires explanation, concerns the fact that the output of
the algorithm consists of hypersurfaces containing the graphs of approximating maps. Thus the
question is how to recover the output for our original $F:U\rightarrow V$ if we proceed with
$\hat{J}\circ F:U\rightarrow\hat{J}(V)$ instead of $F$ and obtain
$\hat{F}_{\nu}:U_0\rightarrow\hat{J}(V).$ Observe that the components of the Nash map
$F_{\nu}=\hat{J}^{-1}\circ\hat{F}_{\nu}:U_0\rightarrow V$ approximating $F|_{U_0}$ are linear
combinations of the components of the map $\hat{F}_{\nu}:U_0\rightarrow\hat{J}(V)$
approximating $\hat{J}\circ F.$ {This implies, by Lemma \ref{algNashfunct}, that the output
for $F$ can be recovered from the output obtained for $\hat{J}\circ F$.}

Set
$F^{*}=(f_1,\ldots,f_m,f_{m+1}):U\rightarrow\Phi(V).$ Now the idea is to approximate the map $F^{*}:U\rightarrow\Phi(V)$
and then recover the output data for $F$ having such data for
$F^{*}.$ To do this, we first need to calculate (Step 3) the
optimal polynomial $P$ for $\Phi(V)$
and the discriminant $R$ of $P.$
Computing the optimal polynomial $P$ for $\Phi(V)$ (given $V$) is discussed in \cite{JS}, p. 240-241 (where $P$
is called the minimal polynomial). Let us note that
the problem of effective elimination of variables (an instance of which is computing $\Phi(V)$) has been discussed in several works (cf. \cite{BeWa}, \cite{CLOS}, Chapter 3, or \cite{GP}, pp. 69-73).

We need $R(f_1,\ldots,f_m)\neq 0,$ which can be achieved by repeating Steps 2-4 in a loop.
More precisely, if in Step 4, {we have} $R(f_1,\ldots,f_m)=0,$ then we can repeat Steps 2, 3
with $m,$ $s,$ $V$ replaced by $m-1,$ $s+1,$ $V_1=V\cap\{R=0\}$ (the new variety is purely
$(m-1)$--dimensional) and check the condition again. In this process, the dimension of the
target variety drops, so finally the required condition is satisfied. (In Step 4 we have the
zero-test discussed in Section \ref{model}.) Note that if {$R$ and the discriminant $R^{*}$
of $R$  vanish identically on the images of the corresponding maps,} {then after Steps~2, 3,
4 performed with $m-1,$ $s+1,$ $V_1,$ the algorithm will return to Step~2 again. For such
$R$, instead of performing Steps 2, 3 with $m-1, s+1, V_1$ (to obtain some new variety $V_2$
in Step 4) we can alternatively define $V_2=V_1\cap\{R^{*}=0\}$} (then
$\mathrm{dim}(V_2)=\mathrm{dim}(V_1)-1$). If {in addition} $R^{*}$ is reducible, then we can
{reduce} the complexity of $V_2$ by taking (instead of $R^{*}$) a factor of $R^{*}$ whose
zero-set contains the image of the map. Here it is important that there are direct methods of
computing factors of iterated discriminants (such as $R^{*}$) which in many cases are much
more efficient than computing and factorizing iterated discriminants (for details see
\cite{LMC}). Let us make the {alternative} construction of $V_2$ more precise.

First observe that $V_R=\{(z_1,\ldots,z_m)\in\mathbf{C}^m: R(z_1,\ldots,z_m)=0\}$ may be
assumed to have proper projection onto $\mathbf{C}^{m-1}_{z_1\ldots z_{m-1}}.$ Then
$V_1=V\cap\{R=0\}$ has proper projection onto $\mathbf{C}^{m-1}_{z_1\ldots z_{m-1}}.$ Choose
$R^*\in\mathbf{C}[z_1,\ldots,z_{m-1}]$ such that $\{R^*=0\}\subset\mathbf{C}^{m-1}$ is the
set of points over which the fiber in $V_R$ has not the maximal cardinality. If the
discriminant of $R$ is non-zero, then define $R^*$ to be the discriminant of $R.$ It is clear
that if $R^*(f_1,\ldots,f_{m-1})=0,$ then the fibers in $V_1$ over the image of
$(f_1,\ldots,f_{m-1})$ have smaller cardinality than generic fibers in $V_1$ over
$\mathbf{C}^{m-1}.$ Hence, the discriminant $R_1$ that would be obtained if we performed Step
3 with $V_1$ ($R_1$ is related to $V_1$ in the same way as $R$ is related to $V$) would
vanish identically on the image of $(f_1,\ldots,f_{m-1}).$ In this case, replace $V_1$ by
$V_1\cap\{S=0\}$ (of strictly smaller dimension), where $S$ is a reduced factor of $R^{*}$
vanishing on the image of $(f_1,\ldots,f_{m-1}).$

Finally note that if we know that the image of
$F$ is not contained in $\mathrm{Sing}(V),$ then we need not replace $V$ by a smaller variety
in order to achieve\linebreak $R(f_1,\ldots,f_m)\neq 0,$ because the generic linear change of variables
in Step 2 gives this condition.

Now to approximate the map $F^{*}|_{U_0}$ by some Nash map $F^{*,\nu}:U_0\rightarrow\Phi(V),$ where $U_0$ is some
neighborhood of zero, we use Lemma \ref{vddlem} which implies that it is sufficient to find Nash functions
$f_1^{\nu},\ldots, f_m^{\nu} ,\bar{f}^{\nu} $ approximating $f_{1},\ldots,f_m,{f}_{m+1},$
respectively, on $U_0$ such that\vspace*{2mm}\\
$\mathrm{(3.4)}$\hspace*{17mm}$P(f_1^{\nu},\ldots,f_m^{\nu},\bar{f}^{\nu})=
C^{\nu}(\frac{\partial{P}}{\partial{z_{m+1}}}(f_1^{\nu},\ldots,f_m^{\nu},\bar{f}^{\nu}))^2,$
\vspace*{2mm}\\for some holomorphic
function $C^{\nu}$ with small norm.

To find such functions we need to handle the zero-set of $\frac{\partial P}{\partial
z_{m+1}}(f_1,\ldots,f_{m+1}).$ (Recall that $R(f_1,\ldots,f_m)\neq 0,$ so $\frac{\partial
P}{\partial z_{m+1}}(f_1,\ldots,f_m,{f}_{m+1})\neq 0.$) For that reason we change the
coordinates in a neighborhood of zero in $\mathbf{C}^n$ (Step~5) after which the Weierstrass
Preparation Theorem can be applied to $\frac{\partial P}{\partial
z_{m+1}}(f_1,\ldots,f_{m+1}).$ More precisely, $(f_{1},\ldots,{f}_{m+1})$ may have to be
replaced  by $(f_{1}\circ J ,\ldots,f_{m+1}\circ J),$ where
$J:\mathbf{C}^n\rightarrow\mathbf{C}^n$ is a generic linear isomorphism. Since the
coordinates are changed in the domain of the approximated map (and not in its range), the
difficulties which we discussed in Step 2 do not appear here. In view of the Weierstrass
Preparation Theorem (the effective version of which is discussed in Section \ref{model}), we
assume that in some neighborhood of $0\in\mathbf{C}^n,$ the zero-set of the function
$\frac{\partial P}{\partial z_{m+1}}(f_1,\ldots,{f}_{m+1})$ is given by the zero-set of a
monic polynomial $W$ in $x_n$ with holomorphic coefficients depending on
$x'=(x_1,\ldots,x_{n-1}),$ where $x_1,\ldots,x_n$ are the coordinates in $\mathbf{C}^n.$

The problem of finding $f_1^{\nu},\ldots,f_m^{\nu},\bar{f}^{\nu}$
satisfying the equation $\mathrm{(3.4)}$ can be solved using recursion. The aim
of Steps 6 and 7 is to prepare the setup for this recursion. More
precisely, we define a new map $g$ depending only on the variables
$x_1,\ldots,x_{n-1}$ whose image is contained in some algebraic
variety and whose approximation enables us to obtain the solution
to $\mathrm{(3.4)}$. To define $g$, we need the effective version of the Division
Theorem (see Section \ref{model}). Step 8 is devoted to the recursive application of the
algorithm.

The aim of Step 9 is to recover the solution to $\mathrm{(3.4)}$ from the data
obtained by the recursive application of the algorithm. More
precisely, we recover the functions $f_1^{\nu},\ldots,f_m^{\nu}$, which will
constitute the first $m$ components of the map $F^{\nu}$ approximating the
original map $F.$
This means that, for $i=1,\ldots,m,$ we construct a polynomial
$P_i^{\nu}(x,z_i)$ whose zero-set contains the graph of
$f_i^{\nu}.$  As for the function $\bar{f}^{\nu}$, it is sufficient
to know that such a function satisfying $\mathrm{(3.4)}$ together with
$f_1^{\nu},\ldots,f_m^{\nu}$ exists. Here one can simply take
$\bar{f}^{\nu}=H^{\mu}_{m+1}W_{\mu}^2+r_{m+1,\mu},$ where $H_{m+1}^{\mu}$ is a polynomial
approximation of $H_{m+1}.$ Then, by Lemma \ref{vddlem},
there exists a Nash approximation ${f}_{m+1}^{\nu}$ of ${f}_{m+1}$
such that $P(f_1^{\nu},\ldots,f_m^{\nu},{f}_{m+1}^{\nu})=0$ which allows us
to recover the output data for the remaining functions
$f^{\nu}_{m+1},\ldots,f^{\nu}_{m+s}$ satisfying the required
properties. How to construct these last $s$ components of
$F^{\nu}$ is discussed in Step 10.

Consider the algebraic variety $V^{\nu}$ defined in Step 10. Since $V\subset\mathbf{C}^m\times\mathbf{C}^s$ has
proper projection onto $\mathbf{C}^m$ and, for $i=1,\ldots,m,$ $P_i^{\nu}(x,z_i)$ are monic in $z_i,$ $V^{\nu}$
has proper projection onto $\mathbf{C}^n_x.$ Hence, $V^{\nu}$ also has proper projection to
$\mathbf{C}^n_x\times\mathbf{C}_{z_{m+i}},$ for $i=1,\ldots,s.$ Therefore the image $V^{\nu,i}$ of the projection
of $V^{\nu}$ to $\mathbf{C}^n_x\times\mathbf{C}_{z_{m+i}}$ is an algebraic hypersurface with proper projection
onto $\mathbf{C}^n_x.$ Now observe that there are holomorphic functions $f_{m+2}^{\nu},\ldots,f_{m+s}^{\nu}$ such
that the image of $(f_1^{\nu},\ldots,f_{m+1}^{\nu},\ldots,f_{m+s}^{\nu})$ is contained in $V.$ This follows by the
fact that $P(f_1^{\nu},\ldots,f_m^{\nu},{f}_{m+1}^{\nu})=0$ and by Lemma~\ref{funnyobvious} with
$L(z_{m+1},\ldots,z_{m+s})=z_{m+1}$ and $Z=\mathcal{V}(V,H)$ where $H=(f_1^{\nu},\ldots,f_m^{\nu})$ (cf. the
discussion preceding Algorithm~1). Consequently, $V^{\nu,i}$ contains the graph of $f_{m+i}^{\nu},$ for
$i=1,\ldots,s$; in particular, every $f_{m+i}^{\nu}$ is Nash. Finally, define $P_{m+i}^{\nu}$ to be the optimal
polynomial for $V^{\nu,i},$ for $i=1,\ldots,s.$

In the following subsection we check that $f_{m+i}^{\nu}$ do
approximate $f_{m+i}$ and show how to control the error of this approximation.\vspace*{3mm}\\
{\textbf{3.2.2.1 Controlling accuracy of approximation}}\vspace*{2mm}\\
Assume in Step 5 that $\frac{\partial P}{\partial z_{m+1}}(f_1,\ldots,f_{m+1})(0)=0,$ since otherwise we have
$F(0)\in\mathrm{Reg}(V),$ which has been discussed already. It is sufficient to find upper bounds for the distance
between $g$ and $g^{\mu}$ (that is, to compute suitable $\mu(\nu)$ in Step~8) and for the distance between
$H_i^{\mu}$ and $H_i$ (Step 9) which guarantee that $||F-F^{\nu}||_{U_0}<\frac{1}{\nu}.$

First, using the definition of $f_i^{\nu}$ (Step 9), for $i=1,\ldots,m,$ one can determine how close $H_i^{\mu},
g^{\mu}$ {have to} be to $H_i,g$ to ensure that $||f_i-f_i^{\nu}||_{U_0}<\frac{1}{\nu}.$

As for $f_{m+1},$ observe that, by the formulas in Step 7 and by Steps 8, 9 we have:
$$P(f_1^{\nu},\ldots,f_{m}^{\nu}, \bar{f}^{\nu})=
C^{\nu}(\frac{\partial P}{\partial z_{m+1}}(f_1^{\nu},\ldots,f_m^{\nu},\bar{f}^{\nu}))^{2},$$ where
$\bar{f}^{\nu}=H_{m+1}^{\mu}W_{\mu}^{2}+r_{m+1,\mu},$ and $H_{m+1}^{\mu}$ is a polynomial approximation of
$H_{m+1}$, and $C^{\nu}$ is a holomorphic function. More precisely, $C^{\nu}$ is the composition of
${\tilde{W}}/{\tilde{S}^2}$ with $g^{\mu},H_1^{\mu},\ldots,H_{m+1}^{\mu}.$ Observe that the composition of
$\tilde{S}^2$ with $g, H_1,\ldots, H_{m+1}$ is non-vanishing in some neighborhood {of the origin, whereas the
composition of $\tilde{W}$ with the same functions is identically zero. Moreover, in some neighborhood of every
point at which $\tilde{S}$ is non vanishing we have bounds} for the partial derivatives of
${\tilde{W}}/{\tilde{S}^2}.$ Therefore we can compute {an open bounded polydisc $U_1$ in a neighborhood of
$0\in\mathbf{C}^n$ in which the algorithm is performed,} such that for every $\varepsilon>0$ we can compute
$\delta>0$ such that if the distance between $g^{\mu}, H_i^{\mu}$ and $g, H_{i}$ is smaller than $\delta$, then
{$||C^{\nu}||_{U_1}<\varepsilon.$} For $\varepsilon$ small enough, by Lemma \ref{vddlem} and Remark \ref{vdrema},
we have our Nash function $f^{\nu}_{m+1}$ with $P(f_1^{\nu},\ldots,f_m^{\nu},f_{m+1}^{\nu})=0$ and
{$||f_{m+1}^{\nu}-f_{m+1}||_{U_1}<\frac{1}{\nu}.$}

As for $f_{m+j},$ for $j=2,\ldots,s,$ since $V\subset\mathbf{C}^m\times\mathbf{C}^s$ satisfies (x) and (y)
(Step~2), there are {effectively computable} $Q_{2},\ldots, Q_{s}\in(\mathbf{C}[z_1,\ldots,z_m])[z_{m+1}]$ such
that
$$V\setminus R^{-1}(0)=\{R\neq 0, P=0, z_{m+j}R=Q_j,  j=2,\ldots,s\},$$
where $R\in\mathbf{C}[z_1,\ldots,z_m]$ is the discriminant of the optimal polynomial
$P\in(\mathbf{C}[z_1,\ldots,z_m])[z_{m+1}]$ for $\Phi(V)$ {(cf. \cite{JS} p. 233, \cite{GH},
\cite{KP})}.

Denote $\hat{R}(x)=R(f_1(x),\ldots,f_m(x)).$ After Step 4, $\hat{R}\neq 0,$ so we may assume that
$\hat{R}(o,\cdot)$ has a zero of finite order at $x_n=0,$ where $o\in\mathbf{C}^{n-1}$ is the origin. Now we can
compute polydiscs $A\subset\mathbf{C}^{n-1}, B\subset\mathbf{C},$ centered at zero, and a real number $c>0$ such
that {$\overline{A\times B}\subset U_1$}, and $\inf_{\overline{A}\times\partial B}|\hat{R}|>c$ {(cf. the
discussion on the Division Theorem in Section \ref{model}).} We know how close $g^{\mu}, H_i^{\mu}$ to $g, H_i$
{have to} be to ensure that $\mathrm{inf}_{\overline{A}\times\partial B}|R(f_1^{\nu},\ldots,f_{m}^{\nu})|> 0.$
Then, once we have a Nash function $f_{m+1}^{\nu}$ with $P(f_1^{\nu},\ldots,f_m^{\nu},f_{m+1}^{\nu})=0,$ {by
Lemma~\ref{funnyobvious} with $L(z_{m+1},\ldots,z_{m+s})=z_{m+1}$ and $Z=\mathcal{V}(V,H)$ where
$H=(f_1^{\nu},\ldots,f_m^{\nu})$, we obtain (cf. the discussion preceding Algorithm 1)} holomorphic functions
$f_{m+2}^{\nu},\ldots,f_{m+s}^{\nu}$ such that the image of $(f_1^{\nu},\ldots,f_{m+s}^{\nu})$ is contained in
$V.$

On (some neighborhood of) $\overline{A}\times\partial B$ we have, for $j=2,\ldots,s,$
$$f_{m+j}=\frac{Q_j(f_1,\ldots,f_m,{f}_{m+1})}{R(f_1,\ldots,f_m)},\mbox{ }
f_{m+j}^{\nu}=\frac{Q_j(f_1^{\nu},\ldots,f_m^{\nu},{f}_{m+1}^{\nu})}{R(f_1^{\nu},\ldots,f_m^{\nu})}.$$
Consequently, ($f_{m+2}^{\nu},\ldots,f_{m+s}^{\nu}$ are Nash and) we can compute {$\delta>0$} such that if
$||{f}_{i}-{f}_i^{\nu}||_{{A}\times B}<{\delta},$ for $i=1,\ldots,m+1,$ then
$||f_{m+j}-f_{m+j}^{\nu}||_{\overline{A}\times\partial B}<\frac{1}{\nu}$ (hence, by the Maximum Principle,
$||f_{m+j}-f_{m+j}^{\nu}||_{{A}\times B}<\frac{1}{\nu}$) for $j=2,\ldots,s.$ But, as shown above, we know how
close $H_i^{\mu}, g^{\mu}$ to $H_i, g$ {have to} be to ensure that $||f_i-f_i^{\nu}||_{A\times B}<{\delta},$ for
$i=1,\ldots,m+1.$ This means that we can control the error of approximation of $F$ by $F^{\nu}$ on $A\times
B.$\vspace*{2mm}\\
{\textbf{3.2.2.2 Bound for the degree of output polynomials}}\vspace*{2mm}\\
We prove the claim: there is a bound for the degree of $P_i^{\nu}$ in $z_i,$ $i=1,\ldots,\hat{m},$ independent of
$\nu$. If $h(x)$ is a Nash function defined in some open neighborhood of $0\in\mathbf{C}^n,$ then by an implicit
form of $h$ we mean a polynomial $P_h(x,z_i)$ monic in $z_i$ with  $P_h(x,h(x))=0$ in some neighborhood of
$0\in\mathbf{C}^n.$ By the degree of the implicit form $P_h(x,z_i)$ we mean the degree of $P_h$ in $z_i.$

If $n=0$ or $m=\mathrm{dim}(V)=0,$ then the claim is clearly true because $F$ is a constant map and $F^{\nu}=F$.

{Consider the linear ordering of $\mathbf{N}^2$ defined by $(l,k)<(l',k')$ iff $l<l'$ or ($l=l'$ and $k<k'$).
Assume that $(0,0)<(n,m).$ Let $V$ be an algebraic subvariety of $\mathbf{C}^{\hat{m}}$ with $\mathrm{dim}(V)=m$
and let $F:U\rightarrow V$ be any holomorphic map, where $U$ is an open neighborhood of $0$ in $\mathbf{C}^n.$
Assume that the claim is true for {every} holomorphic map $\tilde{F}:\tilde{U}\rightarrow\tilde{V},$ where
$\tilde{U}$ is an open neighborhood of $0$ in $\mathbf{C}^l$ and $\mathrm{dim}(\tilde{V})=k$ and $(l,k)<(n,m).$ We
will show that the claim is also true for $F$, which, by induction, will complete the proof.}

First recall that in Step 2 we apply a linear change of the coordinates in $\mathbf{C}^{\hat{m}}$ resulting in
replacing the components of the input map by linear combinations of these components. Here we need to know that if
the degrees of implicit forms of Nash approximations of the components do not depend on $\nu,$ then the degrees of
implicit forms of linear combinations of these Nash approximations do not depend on $\nu$ either. This is an
immediate consequence of Lemma \ref{algNashfunct}.

We may assume that in Step 4 we have $F(U)\nsubseteq\{R=0\}$ because otherwise $F(U)$ is contained in some variety
of dimension smaller than $m$ and the induction hypothesis completes the proof. This implies that $\frac{\partial
P}{\partial z_{m+1}}(f_1,\ldots,f_m,f_{m+1})\neq 0$ and either $\frac{\partial P}{\partial
z_{m+1}}(f_1,\ldots,f_m,f_{m+1})(0)\neq 0$ or $\frac{\partial P}{\partial z_{m+1}}(f_1,\ldots,f_m,f_{m+1})(0)=0.$

If $\frac{\partial P}{\partial z_{m+1}}(f_1,\ldots,f_m,f_{m+1})(0)\neq 0,$ then $F(0)\in \mathrm{Reg}(V_{(m)})$
and the polynomial $P_i^{\nu}$ is linear in $z_i$ for $i=1,\ldots,m.$ (Recall that in the case $F(0)\in
\mathrm{Reg}(V_{(m)})$ the algorithm has been discussed separately before.) For $i=m+1,\ldots,m+s,$ the degree of
$P_i^{\nu}$ in $z_i$ is bounded by the cardinality of the generic fiber in $V$ over $\mathbf{C}^m.$ Indeed, by
construction of $P_i^{\nu}$ for $i=m+1,\ldots,m+s,$ the degree of $P_i^{\nu}$ in $z_i$ does not exceed the number
of points in the generic fiber in ${V}^{\nu}$ over $x$ which (by linearity of $P_i^{\nu}$ in $z_i$ for
$i=1,\ldots,m$) is bounded by the cardinality of the generic fiber in $V$ over $\mathbf{C}^m.$

Now suppose that $\frac{\partial P}{\partial z_{m+1}}(f_1,\ldots,f_m,f_{m+1})(0)=0.$ The linear change of the
coordinates in Step~5 regards the domain of the approximated map and therefore does not affect the bounds we look
for. First, let us discuss $P_i^{\nu}$ for $i=1,\ldots,m.$ Since $H_i^{\mu}$ obtained in Step 9 are polynomials,
they have implicit forms of degree $1.$ The same is true for $x_n^j,$ for $j\in\mathbf{N}.$ By induction
hypothesis, the coefficients of $W_{\mu}, r_{1,\mu},\ldots, r_{m,\mu}$ have implicit forms of degrees bounded by a
constant independent of $\mu.$ {Now recall that to calculate $P^{\nu}_i,$ for $i=1,\ldots,m,$ we perform the ring
operations in the construction of $f_{i}^{\nu}=H_i^{\mu}W_{\mu}^2+r_{i,\mu}$ in their implicit representation.}
{By Lemma~\ref{algNashfunct}, for Nash} functions $f, g$ with implicit forms $P_f, P_g,$ we obtain implicit forms
$P_{f+g}, P_{f\cdot g}$ whose degrees are bounded by a constant depending only on the degrees of implicit forms
$P_f, P_g.$ This shows the claim for $i=1,\ldots,m.$

By Step 10, for $i=1,\ldots,s,$ $P_{m+i}^{\nu}$ is the optimal polynomial describing the image $V^{\nu,i}$ of the
projection of $V^{\nu}$ to $\mathbf{C}_x^n\times\mathbf{C}_{z_{m+i}}.$ So the degree of $P_{m+i}^{\nu}$ in
$z_{m+i}$ equals the cardinality of generic fibers in $V^{\nu,i}$ over $\mathbf{C}^n_x.$ By definition of
$V^{\nu},$ the latter number is bounded by a constant depending only on the cardinality of generic fibers in $V$
over $\mathbf{C}^m$ and on the degree of $P_i^{\nu}$ in $z_i,$ for $i=1,\ldots,m.$ In view of the previous
paragraph, we obtain a bound independent of $\nu.$
\begin{remark}\label{comparis} \emph{The main difference between Algorithm 1 and the me\-thod of approximation presented in \cite{B6} is
that Algorithm 1 does not rely on factorization of polynomials with holomorphic coefficients which appears in Step
7 of the method of \cite{B6}. The factorization is not computable in the model considered in the present paper.}
\emph{To see this, take any holomorphic function $a(x)$ and $W(x,y)=y(y-a(x)).$ Then either $W=W_1^2,$ where
$W_1(x,y)=y$ (which occurs iff $a=0$) or $W=W_1\cdot W_2,$ where $W_1(x,y)=y$ and $W_2(x,y)=y-a(x)$ (which occurs
iff $a\neq 0$), hence, computability of the factorization would imply the decidability of the zero-test for $a.$
In other words, given two very close (possibly equal) factors of $W$ it is not possible to distinguish whether
they are equal or not.}

\emph{Step 7 (of the method of \cite{B6}) could be effectively performed if we knew that $W$ (obtained in Step 5) does not have multiple
factors. But $W$ can have multiple factors (even if $R$ obtained in Step 4 is replaced
by a reduced polynomial).} \emph{Indeed, let us consider any algebraic hypersurface
$V\subset\mathbf{C}^{\hat{m}}$
with $\mathrm{Sing}(V)=\{0\}$ and any non-constant holomorphic map
$f=(f_1,\ldots,f_{\hat{m}}):\mathbf{C}^2\supset U\rightarrow V,$ where $U$ is an open connected neighborhood of $(0,0)$ and, for $j=1,\ldots,\hat{m},$  $$f_j(x_1,x_2)=g_j(x_1,x_2)u(x_1,x_2),$$
for some holomorphic $g_j, u$ with $u(0,0)=0.$}
\emph{With these data let us try the method of \cite{B6} to compute approximation for $f.$}

\emph{First we apply a generic linear change of the coordinates after which
$V\subset\mathbf{C}^{m}_{y}\times\mathbf{C}$ has proper projection onto $\mathbf{C}^{m}_y.$
After this change we also have that the discriminant $R(y)$ of the optimal polynomial for $V$
satisfies $R(f_1,\ldots,f_{{m}})\neq 0$ (because $f(U)\nsubseteq\mathrm{Sing}(V)$).}
\emph{Clearly, $$R(f_1(x_1,x_2),\ldots,f_{{m}}(x_1,x_2))=G(x_1,x_2)u(x_1,x_2),$$ where $G\neq
0$ is holomorphic. (The latter formula remains true, maybe with some other $G,$ after
replacing $R$ by the reduced polynomial having the same zero-set as $R$). Therefore, $W$
obtained in Step 5 has multiple factors if $u$ has multiple factors.}
\end{remark}
\subsubsection{Complete algorithm}\label{completfor1}
Here we improve Algorithm 1 to obtain a complete algorithm of approximation
(without semi-computable steps). Every component of the approximated map will be represented as
described in Section \ref{model}.\vspace*{2mm}\\
\textbf{Input:} a positive integer $\nu,$ an algebraic
variety $V\subset\mathbf{C}^{\hat{m}},$ and a holomorphic mapping
$F=(f_1,\ldots,f_{\hat{m}}):U\rightarrow\mathbf{C}^{\hat{m}},$ where $U$ is an open connected
neighborhood of $0\in\mathbf{C}_x^n$ .\vspace*{2mm}\\
The algorithm either detects that the image of $F$ is not contained in $V$ and then returns "$F(U)\nsubseteq V$"
or computes a Nash map $F^{\nu}:U_0\rightarrow V$ such that $||F^{\nu}-F||_{U_0}<\frac{1}{\nu}.$ Note that the
problem of testing whether $F(U)\nsubseteq V$ is not decidable (it is semi-decidable), but one can perform at
least one of two tasks stated above. In other words, {the algorithm may compute} an approximation of $F$ into $V$
even if $F(U)\nsubseteq V$ (but $F(U)$ is very close to $V$). The idea to give the algorithm the choice "either
detect or approximate" comes from the fact that the existence of a non-exact solution $(f_1,\ldots,f_{m+1})$ of
the equation $P(z_1,\ldots,z_{m+1})=0$ satisfying (3.4) (with small $C^{\nu}$) implies the existence of an\\
\\ exact solution (cf. Lemma \ref{vddlem}). If we have $F(U)\subset V,$ then the algorithm necessarily computes
an approximation.\vspace*{2mm}\\
\textbf{Output:} Either "$F(U)\nsubseteq V$" or
$P_i^{\nu}(x,z_i)\in(\mathbf{C}[x])[z_i],$
$P_i^{\nu}\neq 0$ for $i=1,\ldots,\hat{m},$
with the following properties:\vspace*{1mm}\\
(a) $P_i^{\nu}(x,f_i^{\nu}(x))=0$ for every $x\in U_0,$ where
$F^{\nu}=(f^{\nu}_1,\ldots,f^{\nu}_{\hat{m}}):U_0\rightarrow V$ is
a holomorphic mapping such that $||F-F^{\nu}||_{U_0}<\frac{1}{\nu}$ and
$U_0$ is an open neighborhood of
$0\in\mathbf{C}^n$ independent of $\nu,$\\
(b) $P_i^{\nu}$ is a monic polynomial in $z_i$ of degree in $z_i$ bounded by a constant
independent of $\nu$.\vspace*{2mm}

Let us first discuss two simple cases to which the algorithm reduces the problem by recursive calls.
Recall that we work under the assumption that the zero test can be performed in $\mathbf{C}$
(cf. Section \ref{model}).

If $n=0$ (i.e. $F$ is a constant map), then we obtain output data immediately. Now observe that if $V$ is zero-dimensional, then, for any $n,$
we obtain output data almost immediately. Indeed, the algorithm works as follows:
check whether $F(0)\in V.$
If $F(0)\notin V$, then return "$F(U)\nsubseteq V$".
Otherwise compute $U_0$ and $\kappa\in\mathbf{N}$
such that if all the coefficients (except for the terms of order zero) of the Taylor expansion
of $F$ at $0$ up to order $\kappa$ vanish, then
$||F-F(0)||_{U_0}<\frac{1}{\nu}.$ Check whether these coefficients vanish. If this is true, then
the constant map $x\mapsto F(0)$ is the required approximation. Otherwise, $F$ is
a non-constant map defined on an open connected set so its image is not contained in any zero-dimensional variety (return "$F(U)\nsubseteq V$").

Denote $m=\mathrm{dim}(V).$ As mentioned above, the problem will be reduced to the case
where $n\cdot m=0$. More precisely, running with $V, F,$ the algorithm will be called
recursively either with the same map but with a target variety of strictly smaller dimension than
$m,$
or it will be called with some map depending on strictly fewer variables (and then the dimension
of the target variety may even increase).

Before discussing the problem in general let us look at one more special case (where no
recursive calls are necessary). Namely, assume that $V$ is of pure dimension $m>0$ and
$F(0)\in\mathrm{Reg}(V).$ Then, after a generic linear change of the coordinates,
$V\subset\mathbf{C}^{\hat{m}}\approx\mathbf{C}^m\times\mathbf{C}^s$ has proper projection
onto $\mathbf{C}^m.$ Moreover, generic fibers in $V$ and in ${\Phi(V)}$ over $\mathbf{C}^m$
have the same cardinalities, where $\Phi(z_1,\ldots,z_m,z_{m+1},\ldots,z_{m+s})$
$=(z_1,\ldots,z_m,z_{m+1}).$ Furthermore, the fiber in $\Phi(V)$ over
$(f_1,\ldots,f_m)(0)\in\mathbf{C}^m$ has the same number of elements as generic fibers in
$\Phi(V)$ over $\mathbf{C}^m.$

Now calculate the optimal polynomial $P\in(\mathbf{C}[z_1,\ldots,z_m])[z_{m+1}]$ describing
$\Phi(V)\subset\mathbf{C}^m\times\mathbf{C}$ and the discriminant $R$ of $P.$ By the previous paragraph,
{$R(f_1(0),\ldots,f_m(0))\neq 0\neq \frac{\partial P}{\partial z_{m+1}}(f_1(0),\ldots,f_{m+1}(0)),$ so we can
compute a polydisc $E$ centered at $0\in\mathbf{C}^n$ and $\tilde{c}>0$ such that
$\inf_{E}|R(f_1,\ldots,f_m)|>\tilde{c}$} and {$\inf_E|\frac{\partial P}{\partial
z_{m+1}}(f_1,\ldots,f_{m+1})|>\tilde{c},$} so
$$C(x)=\frac{P(f_1(x),\ldots,f_{m+1}(x))}{(\frac{\partial P}{\partial
z_{m+1}}(f_1(x),\ldots,f_{m+1}(x)))^2}$$ is a holomorphic function in $E.$

Let us show that for every $\delta>0,$ we can compute $\kappa$ such that if only the coefficients of the Taylor
expansion of $P(f_1,\ldots,f_{m+1})$ up to order $\kappa$ vanish, then there are polynomials $f_1^{\nu},\ldots,
f_m^{\nu}$ and a Nash function $f_{m+1}^{\nu}$ whose distance from $f_1,\ldots,f_{m+1},$ respectively, is bounded
by $\delta$ and $P(f_1^{\nu},\ldots,f_{m+1}^{\nu})=0.$

For every $\varepsilon>0$ we can compute $\kappa$ such that if the coefficients of the Taylor expansion of
$P(f_1,\ldots,f_{m+1})$ at zero up to order $\kappa$ vanish, then $||C||_E<\frac{\varepsilon}{2}.$ Indeed, {it is
sufficient to compute $\kappa$ such that if the coefficients of the Taylor expansion of $P(f_1,\ldots,f_{m+1})$ at
zero up to order $\kappa$ vanish, then $||P(f_1,\ldots,f_{m+1})||_E<\frac{\varepsilon\tilde{c}^2}{2}.$} This can
be done in the same way as computing $N(\theta,\varepsilon,\hat{R})$ in the last but one paragraph of Section
\ref{model}.

{If $||C||_E<\frac{\varepsilon}{2},$ then there is $\gamma>0$ such that for Taylor polynomials
$f_1^{\nu},\ldots,f_m^{\nu}$ of $f_1,\ldots,f_m,$ respectively, with $||f_i-f_{i}^{\nu}||_E<\gamma$,
$i=1,\ldots,m,$} we have
$$\Big{|}\Big{|}\frac{P(f_1^{\nu},\ldots,f_m^{\nu},{f}_{m+1})}{(\frac{\partial P}{\partial
z_{m+1}}(f_1^{\nu},\ldots,f_m^{\nu},{f}_{m+1}))^2}\Big{|}\Big{|}_{E}<\varepsilon.$$ {The number $\gamma$ is
effectively computable (by means of bounds for the partial derivatives of ${P}$ and $\frac{\partial P}{\partial
z_{m+1}}$ on the polydisc in $\mathbf{C}^{m+1}$ centered at the origin, of radius $\tilde{K}+1,$ where $\tilde{K}$
is a common bound for $f_1,\ldots,f_{m+1}$ on $E$).}

{By the previous paragraph, Lemma \ref{vddlem} and Remark~\ref{vdrema}, for every $\delta>0,$\linebreak we can
compute $\varepsilon>0$ such that if $||C||_E<\frac{\varepsilon}{2}$,} {then for polynomials
$f_1^{\nu},\ldots,f_m^{\nu}$ approximating $f_{1},\ldots,f_m$ close enough, there is a Nash function
$f_{m+1}^{\nu}$ such that $||f_{m+1}-f_{m+1}^{\nu}||_E<\delta,$} {and $P(f_1^{\nu},\ldots,f_{m+1}^{\nu})=0.$
Hence, for every $\delta>0,$\linebreak we can compute $\kappa$ such that if the coefficients of the Taylor
expansion of $P(f_1,\ldots,f_{m+1})$ up to order $\kappa$ vanish, then there are polynomials $f_1^{\nu},\ldots,
f_m^{\nu}$ and a Nash function $f_{m+1}^{\nu}$ whose distance from $f_1,\ldots,f_{m+1},$ respectively, is bounded
by $\delta$ and $P(f_1^{\nu},\ldots,f_{m+1}^{\nu})=0$}.

{Recall there are  effectively computable} $Q_{2},\ldots, Q_{s}\in(\mathbf{C}[z_1,\ldots,z_m])[z_{m+1}]$ such that
$$V\setminus R^{-1}(0)=\{R\neq 0, P=0, z_{m+j}R=Q_j,  j=2,\ldots,s\}$$
{(cf. \cite{JS} p. 233, \cite{GH}, \cite{KP}).}

{Now compute nonnegative $\delta<\frac{1}{\nu}$ so small that if $||f_i-f_i^{\nu}||_E<\delta,$ for
$i=1,\ldots,m+1,$ then $R(f_1^{\nu},\ldots,f_m^{\nu})$ does not vanish in $E$ and }
{$$\Big{|}\Big{|}\frac{Q_j(f_1,\ldots,f_{m+1})}{R(f_1,\ldots,f_m)}-
\frac{Q_j(f_1^{\nu},\ldots,f_{m+1}^{\nu})}{R(f_1^{\nu},\dots, f_m^{\nu})}\Big{|}\Big{|}_E<\frac{1}{2\nu}.$$}

{Next compute $\kappa$ such that if only the coefficients of the Taylor expansion of $P(f_1,\ldots,f_{m+1})$ up to
order $\kappa$ vanish, then there are polynomials $f_1^{\nu},\ldots, f_m^{\nu}$ and a Nash function
$f_{m+1}^{\nu}$ whose distance from $f_1,\ldots,f_{m+1},$ respectively, on $E$ is bounded by $\delta$ and
$P(f_1^{\nu},\ldots,f_{m+1}^{\nu})=0.$ Check whether all coefficients of the Taylor expansion of
$P(f_1,\ldots,f_{m+1})$ up to order $\kappa$ vanish; if not, then return "$F(U)\nsubseteq V$".}

{Let  $f_1^{\nu},\ldots, f_m^{\nu}$ be polynomials and $f_{m+1}^{\nu}$ a Nash function
 whose distance from $f_1,\ldots,f_{m+1},$ respectively, on $E$ is bounded by
$\delta$ and $P(f_1^{\nu},\ldots,f_{m+1}^{\nu})=0$.} {Then, by Lemma~\ref{funnyobvious} with
$L(z_{m+1},\ldots,z_{m+s})=z_{m+1}$ and $Z=\mathcal{V}(V,H)$ where $H=(f_1^{\nu},\ldots,f_m^{\nu})$, we conclude
(cf. the discussion preceding Algorithm~1) that there are holomorphic functions
$f_{m+2}^{\nu},\ldots,f_{m+s}^{\nu}$ such that the image of $(f_1^{\nu},\ldots,f_{m+s}^{\nu})$ is contained in
$V.$} Consequently, the functions $f_{1}^{\nu},\ldots,f_{m+s}^{\nu}$ satisfy the equations
$$f_{m+j}^{\nu}R(f_1^{\nu},\ldots,f_{m}^{\nu})=Q_j(f_1^{\nu},\ldots,f_{m+1}^{\nu}),\mbox{ for } j=2,\ldots,s$$
(which in particular implies that $f_{m+2}^{\nu},\ldots,f_{m+s}^{\nu}$ are not only holomorphic, but Nash as
well).

The functions $f_1,\ldots,f_{m+s}$ may not satisfy these equations (if $F(U)\nsubseteq V$). But if the
coefficients of the Taylor expansion of
$$f_{m+j}R(f_1,\ldots,f_{m})-Q_j(f_1,\ldots,f_{m+1})$$ at zero up to sufficiently high order $\beta$ vanish,
for $j=2,\ldots,s,$ then $f_{m+2}^{\nu},\ldots,$ $f_{m+s}^{\nu}$ are close to $f_{m+2},\ldots,f_{m+s}$ if
$f_1^{\nu},\ldots,f_{m+1}^{\nu}$ are close to $f_1,\ldots,f_{m+1},$ respectively. {More precisely, we compute
$\beta$ such that if the vanishing condition is satisfied up to order $\beta,$ then
$$||f_{m+j}R(f_1,\ldots,f_{m})-Q_j(f_1,\ldots,f_{m+1})||_E<\frac{\tilde{c}}{2\nu}.$$}
{This can be done in the same way as computing $N(\theta,\varepsilon,\hat{R})$ in the last but one paragraph of
Section \ref{model}.} {If the vanishing condition is not satisfied, then we have $F(U)\nsubseteq V$. Otherwise,
$||f_{m+j}-f^{\nu}_{m+j}||_E\leq$
$$\Big{|}\Big{|}f_{m+j}-\frac{Q_j(f_1,\ldots,f_{m+1})}{R(f_1,\ldots,f_m)}\Big{|}\Big{|}_E+\Big{|}\Big{|}\frac{Q_j(f_1,\ldots,f_{m+1})}{R(f_1,\ldots,f_m)}-
\frac{Q_j(f_1^{\nu},\ldots,f_{m+1}^{\nu})}{R(f_1^{\nu},\dots,
f_m^{\nu})}\Big{|}\Big{|}_E<\frac{1}{\nu},$$ for $j=2,\ldots,s.$}

When we know that the approximation exists ({on $U_0=E$}), the polynomials $P_{i}^{\nu},$ can be computed as
follows. For $i=1,\ldots,m,$ set $P_i^{\nu}(x,z_i)=z_i-f_{i}^{\nu}(x).$ Next define
$$V^{\nu}=\{(x,z)\in\mathbf{C}^n_x\times\mathbf{C}^{m+s}_{z} : z\in
V, z_i=f_{i}^{\nu}(x), \mbox{ for } i=1,\ldots,m\},$$ where
$z=(z_1,\ldots,z_m,z_{m+1},\ldots,z_{m+s}).$ Finally, for $i=1,\ldots,s,$
take
$P^{\nu}_{m+i}\in$\linebreak $(\mathbf{C}[x])[z_{m+i}]$ to be the optimal
polynomial describing the image of the projection of $V^{\nu}$
to $\mathbf{C}^n_x\times\mathbf{C}_{z_{m+i}}.$ (This projection is proper so the image
is an algebraic variety.)

Let us present the main algorithm solving the problem in the general case. {First, the
algorithm reduces the problem, by recursive calls, to the case where $V$ is purely
dimensional. Then it either detects that $F(U)\nsubseteq V$ or reduces the problem to the
case where $V$ is a hypersurface. It is possible that after the reduction the output will
turn out to be "$F(U)\nsubseteq V$". By reducing to the case of a hypersurface we mean that
the algorithm produces a single polynomial $P\in(\mathbf{C}[z_1,\ldots,z_m])[z_{m+1}]$ monic
in $z_{m+1}$ and a new holomorphic map $F^{\star}:U\rightarrow\mathbf{C}^{m+1}$ such that the
following hold. The composition $\frac{\partial P}{\partial z_{m+1}}\circ F^{\star}$ is a
non-zero function and if $P\circ F^{\star}\neq 0$ then $F(U)\nsubseteq V.$ Moreover, from any
Nash approximation $F^{\star,\nu}$ of $F^{\star}$ satisfying $P\circ F^{\star,\nu}=0$ and
$\frac{\partial P}{\partial z_{m+1}}\circ F^{\star, \nu}\neq 0$  we can recover a Nash map
$F^{\nu}$ into $V$ such that $F^{\nu}$ approximates $F$ if $F(U)$ is close to $V$ (even if
$F(U)\nsubseteq V$).}

{If the algorithm has produced $P, F^{\star}$ as described above, then, using $P,
\frac{\partial P }{\partial z_{m+1}},$ $F^{\star}$ it either recognizes that $F(0)$ is a
regular point of $V$ (this case has been discussed separately above) or performs an operation
the outcome of which is one of the two following things. The first thing is the detection of
$F(U)\nsubseteq V$, the second one is confirmation of the fact that $F(U)$ is close to $V$
together with a Nash map $F^{\star,\nu}$ approximating $F^{\star}$ and satisfying $P\circ
F^{\star,\nu}=0$ and $\frac{\partial P}{\partial z_{m+1}}\circ F^{\star, \nu}\neq 0.$}

{The operation consists of two stages. The first stage is checking the fulfillment of a
condition which implies that $F(U)$ is close to $V$, and whose negation gives $F(U)\nsubseteq
V.$ The second stage is producing another holomorphic map $g$ depending on smaller number of
variables than $F$ and a new algebraic set $X$ with the following properties. If the image of
$g$ is not contained in $X$ then $P\circ F^{\star}\neq 0$ (and then $F(U)\nsubseteq V$). If
there is a Nash approximation of $g$ into $X,$ then there is a Nash approximation $\hat{F}$
of $F^{\star}$ into $\mathbf{C}^{m+1}$ and a Nash function $C$ close to zero such that
$P\circ \hat{F}=C\cdot(\frac{\partial P}{\partial z_{m+1}}\circ\hat{F})^2$ (so, by Lemma
\ref{vddlem}, we have a Nash approximation $F^{\star,\nu}$ of $F^{\star}$ with $P\circ
F^{\star,\nu}=0$ and $\frac{\partial P}{\partial z_{m+1}}\circ F^{\star, \nu}\neq 0,$ and
then, as mentioned above, we obtain a Nash map $F^{\nu}$ into $V$ approximating $F$). }

{ To obtain a Nash approximation of $g$ into $X$ or to find out that the image of $g$ is not contained in $X$ the
algorithm calls itself recursively with $g$ and $X.$ The details are as follows.}
\vspace*{2mm}\\
\textbf{Algorithm 2}\vspace*{2mm}\\
\textbf{1.} If $n=0$ or $m=0,$ then proceed as discussed above.\\
\textbf{2.} Compute the equidimensional decomposition $V_{(0)}\cup\ldots\cup V_{(m)}$
of $V.$\\
\textbf{3.} For $j=0,\ldots,m-1$ call the algorithm recursively with $\nu, V_{(j)}, F$
and if it returns an approximating map for at least one $j,$ then stop. If for every $j=0,\ldots,m-1$ the algorithm
returns "$F(U)\nsubseteq V_{(j)}$", then go to Step 4.\\
\textbf{4.} After a generic linear change of the coordinates in $\mathbf{C}^{\hat{m}}$
we have:\\
(x) $V_{(m)}\subset\mathbf{C}^{\hat{m}}\approx\mathbf{C}^m\times\mathbf{C}^s$ has proper
projection onto $\mathbf{C}^m,$\\
(y) generic fibers in $V_{(m)}$ and in ${\Phi(V_{(m)})}$ over $\mathbf{C}^m$ have the same
cardinalities,\\
\hspace*{5.5mm}where $\Phi(z_1,\ldots,z_m,z_{m+1},\ldots,z_{m+s})=(z_1,\ldots,z_m,z_{m+1}).$\\
Apply such a change of the coordinates.\\
\textbf{5.} Calculate the optimal polynomial
$P\in(\mathbf{C}[z_1,\ldots,z_m])[z_{m+1}]$ describing
$\Phi(V_{(m)})\subset\mathbf{C}^m\times\mathbf{C}.$ Calculate the discriminant $R\in\mathbf{C}[z_1,\ldots,z_m]$ of $P.$\\
\textbf{6.} Call the algorithm recursively with $\nu, V_{(m)}\cap\{R=0\}, F.$ If it returns an approximation, then
stop. Otherwise one has $F(U)\nsubseteq V_{(m)}\cap\{R=0\},$ and we {decide whether} $F(U)\nsubseteq V_{(m)}$
(return "$F(U)\nsubseteq V$") or $R(f_1,\ldots,f_m)\neq 0$ {(cf. Remark \ref{interalgebr}).} If
$R(f_1,\ldots,f_m)\neq 0$ then $P(f_1,\ldots,f_m,f_{m+1})\neq 0$ (return "$F(U)\nsubseteq V$") or $\frac{\partial
P}{\partial z_{m+1}
}(f_{1},\ldots,f_{m},{f}_{m+1})\neq 0$ (go to Step 7).\\
\textbf{7.} If $\frac{\partial P}{\partial z_{m+1}}(f_1,\ldots,f_{m+1})(0)\neq 0,$
then either $F(0)\notin V_{(m)}$ (return "$F(U)\nsubseteq V$") or
$F(0)\in\mathrm{Reg}(V_{(m)})$ (case discussed above).
If $\frac{\partial P}{\partial z_{m+1}}(f_1,\ldots,f_{m+1})(0)=0,$ then
apply a linear change of the coordinates in
$\mathbf{C}^n$ after which the following holds: $\frac{\partial
P}{\partial z_{m+1} }(f_1,\ldots,f_{m+1})(x)$$=\hat{H}(x)W(x)$
in some neighborhood of $0\in\mathbf{C}^n,$ where $\hat{H}$ is a
holomorphic function, $\hat{H}(0)\neq 0$ and $W$ is a monic
polynomial in $x_n$ with holomorphic coefficients
vanishing at $0\in\mathbf{C}^{n-1},$
 depending on
$x'=(x_1,\ldots,x_{n-1}).$ Put $d=\mathrm{deg}(W).$\\
\textbf{8.} Choose polydiscs $E_1\subset\mathbf{C}^{n-1}, E_2\subset\mathbf{C},$ centered at zero,
and a real number $\tilde{c}>0$ such that
$\overline{E_1\times E_2}$ is contained in the open neighborhood of the origin in which Step 7 has been performed, $\mathrm{inf}_{\overline{E_1\times E_2}}|\hat{H}|>0,$
$\inf_{\overline{E_1}\times\partial E_2}|R(f_1,\ldots,f_{m})|>\tilde{c},$ and
$\inf_{\overline{E_1}\times\partial E_2}|\frac{\partial P}{\partial z_{m+1}}(f_1,\ldots,f_{m+1})|>\tilde{c}.$\\
\textbf{9.} Compute
$Q_{2},\ldots, Q_{s}\in(\mathbf{C}[z_1,\ldots,z_m])[z_{m+1}]$ such that
$$V_{(m)}\setminus R^{-1}(0)=\{R\neq 0, P=0, z_{m+j}R=Q_j,  j=2,\ldots,s\}.$$
Compute $\kappa\in\mathbf{N}$ such that:\vspace*{0mm}\\
(u) if all the coefficients of the Taylor expansion of $P(f_1,\ldots,f_{m+1})$ at zero\linebreak
\hspace*{6mm}up to order $\kappa$
vanish, then {$||P(f_1,\ldots,f_{m+1})||_{\overline{E_1\times E_2}}<\theta(\nu),$}\\
(v) if all the coefficients of the Taylor expansion of\\
\hspace*{6mm}$f_{m+j}R(f_1,\ldots,f_m)-Q_j(f_1,\ldots,f_{m+1})$ at zero
up to order $\kappa$ vanish, then\linebreak
\hspace*{6mm}$||f_{m+j}R(f_1,\ldots,f_m)-Q_j(f_1,\ldots,f_{m+1})||_{\overline{E_1\times E_2}}<\eta(\nu),$ for $j=2,\ldots,s.$\\
(How to choose {$\theta(\nu), \eta(\nu)$} is described below, where accuracy of approximation is discussed;
{$\kappa$ can be computed similarly to computing $N(\theta,\varepsilon,\hat{R})$ in the last but one paragraph of
Section \ref{model}.}) Check whether all coefficients of the expansions
vanish. If not, then return "$F(U)\nsubseteq V$", otherwise go to Step 10.\\
\textbf{10.} Divide $f_i$ by $(W)^2$ to obtain
$f_i=(W)^2H_i+r_i$ in some
neighborhood of $\overline{E_1\times E_2}\subset\mathbf{C}^n,$ $i=1,\ldots,m+1.$ Here
$H_i$ are holomorphic functions and $r_i$ are
polynomials in
$x_n$ with holomorphic coefficients depending on $x'$ such that $\mathrm{deg}(r_i)<2d.$\\
\textbf{11.} Treating $H_i,$ $i=1,\ldots,m+1,$ and all the
coefficients of $W,r_1,\ldots,r_{m+1}$ as new variables
(except for the leading coefficient $1$ of
$W$) apply the division procedure for
polynomials to obtain:\vspace*{2mm}\\
$P(W^2H_1+r_1,\ldots,W^2H_{m+1}+r_{m+1})$
$=\tilde{W}W^2+x_n^{2d-1}T_1+x_n^{2d-2}T_2+\cdots+T_{2d},$\\
$\frac{\partial P}{\partial
z_{m+1}}(W^2H_1+r_1,\ldots,W^2H_{m+1}+r_{m+1})$\\
\hspace*{51mm}$=\tilde{S}W+x_n^{d-1}T_{2d+1}+x_n^{d-2}T_{2d+2}+\cdots+T_{3d}.$\\
Here $T_1,\ldots,T_{3d}$ are polynomials depending only on the
variables standing for the coefficients of $W,$
$r_1,\ldots,r_{m+1}.$ Let $g$ be the mapping (in $n-1$ variables) whose
components are all these coefficients.\\
\textbf{12.} Call the algorithm recursively with $\nu, V, F$ replaced by $\mu(\nu), $$\{T_1=\cdots=T_{3d}=0\}, g,$
respectively. (How to choose $\mu(\nu)$ is described below.) If the algorithm detects that the image of $g$ is not
contained in $\{T_1=\cdots=T_{3d}=0\},$ then $P(f_1,\ldots,f_{m+1})\neq 0$ (return "$F(U)\nsubseteq V$").
Otherwise, for every component $c(x')$ of $g(x')$ one obtains a monic polynomial
$Q_c^{\mu}(x',t_c)\in(\mathbf{C}[x'])[t_c]$ which put in place of $P_i^{\nu}(x,z_i)$ satisfies (a) and (b) above
with $\nu,$ $x,$ $z_i,$ $f_i^{\nu}, F, F^{\nu}$ replaced by $\mu,$ $x',$ $t_c,$ $c^{\mu},$ $g,$ $g^{\mu},$
respectively. Here, for every $c,$ $c^{\mu}$ is a Nash function approximating $c$ in $\overline{E_1'}\subset E_1$
(where $E_1'$ is an open polydisc centered at the origin and independent of $\mu$) such that the mapping $g^{\mu}$
obtained by replacing every $c$ of $g$ by $c^{\mu}$
satisfies $T_1\circ g^{\mu}=\cdots=T_{3d}\circ g^{\mu}=0.$\\
\textbf{13.} Approximate $H_i$ for $i=1,\ldots,m,$ by
polynomials $H_i^{\mu}.$ (Accuracy of this approximation is discussed
below.) Let $W_{\mu},r_{1,{\mu}},$
$\ldots,$ $r_{m+1,{\mu}}$ be the polynomials in $x_n$ defined by
replacing the coefficients of $W,$ $r_1,\ldots,r_{m+1}$ by
their Nash approxima\-tions (i.e. the components of $g^{\mu}$)
determined in Step 12. Using $Q_c^{\mu}$ (for all $c$) and
$H^{\mu}_i$ one can calculate $P_i^{\nu}\in(\mathbf{C}[x])[z_i],$
for $i=1,\ldots,m,$ satisfying (b) and (a) with
$f_i^{\nu}=H_i^{\mu}(W_{\mu})^2+r_{i,{\mu}}$ being the $i$'th
component of the mapping $F^{\nu}$ (whose last $\hat{m}-m$
components are determined by
$P_{m+1}^{\nu},\ldots,P_{\hat{m}}^{\nu}$ obtained in the next
step). {To calculate $P^{\nu}_i,$ for $i=1,\ldots,m,$ perform the ring operations
in the construction of $f_{i}^{\nu}$ in their implicit representation
(see Lemma \ref{algNashfunct}).}\\
\textbf{14.} Put $V^{\nu}=\{(x,z)\in\mathbf{C}^n_x\times\mathbf{C}^{m+s}_{z} : {z\in V_{(m)}}, P^{\nu}_i(x,z_i)=0
\mbox{ for } i=1,\ldots,m\},$ where $z=(z_1,\ldots,z_m,z_{m+1},\ldots,z_{m+s}).$ For $i=1,\ldots,s,$ take
$P^{\nu}_{m+i}\in$\linebreak $(\mathbf{C}[x])[z_{m+i}]$ to be the optimal polynomial describing the image of the
projection of $V^{\nu}$ to $\mathbf{C}^n_x\times\mathbf{C}_{z_{m+i}}.$ \vspace*{2mm}

Let us list the main differences between Algorithm 1 and Algorithm 2.
First we need Steps 1, 2, 3 in Algorithm 2 because without zero-test we cannot
check which equidimensional components of $V$ (if any) contain the image of $F.$
In Algorithm 1 we pick such a component in Step 1.

Next, Algorithm 2 calls itself recursively in Step 6 because without zero-test
one cannot check whether $R(f_1,\ldots,f_m)\neq 0.$ Observe that if $F(U)\subset V_{(m)},$
then detecting that $R(f_1,\ldots,f_m)\neq 0$ (see Step 4 of Algorithm 1) is equivalent to
detecting that $F(U)\nsubseteq V_{(m)}\cap\{R=0\}$ (see Step 6 of Algorithm 2).

Moreover, in Step 9 of Algorithm 2 we compute {$\kappa, \theta(\nu), \eta(\nu)$ as there is still the possibility
that $f_{m+j}R(f_1,\ldots,f_m)-Q_j(f_1,\ldots,f_{m+1})\neq 0,$ for some $j,$ or $P(f_1,\ldots,f_m,f_{m+1})$ $\neq
0$}
 which implies
$F(U)\nsubseteq V.$ This does not appear in Algorithm 1, where we know that $F(U)\subset
V.$\vspace*{2mm}\\
{\textbf{3.2.3.1 Controlling accuracy of approximation}}\vspace*{2mm}\\
Let us explain how to control accuracy of approximation of $F$ by $F^{\nu}$ on
$U_0=E_1'\times E_2.$ In Step 7, we assume that $\frac{\partial P}{\partial
z_{m+1}}(f_1,\ldots,f_{m+1})(0)=0,$ since in the other case we have either $F(U)\nsubseteq V$
or $F(0)\in\mathrm{Reg}(V_{(m)})$ which has already been discussed.
%%%%%%%%%%%%%%%%%%%%%%%%%%%%%%%%%%%%%%%%%%%%%%%%%%%%%%%%%%%%%%%%%%%%%%%%%%%%%%%%%%%%%%%
{We will show how to choose $\theta(\nu), \eta(\nu)$ in Step 9  and how to chose $\mu(\nu)$ and the bounds for the
distance between $H_i$ and $H_i^{\mu}$ in Steps 12, 13, respectively, to ensure the required precision of
approximation.}

{First compute $\delta<\frac{1}{\nu}$ and $\eta(\nu)$ such that if the inequality in (v), Step 9, holds and if
$||f_i-f^{\nu}_i||_{\overline{E_1'\times E_2}}<\delta,$ $i=1,\ldots,m+1$ , then for $f_{m+j}^{\nu},$
$j=2,\ldots,s,$ satisfying
$$f_{m+j}^{\nu}{R(f_1^{\nu},\ldots,f_m^{\nu})}={Q_j(f_1^{\nu},\ldots,f_m^{\nu},{f}_{m+1}^{\nu})}$$
on $\overline{E_1'\times E_2}$ we have $||f_{m+j}-f_{m+j}^{\nu}||_{\overline{E_1'\times E_2}}<\frac{1}{\nu}.$ To
do this, first note that, since  $\inf_{\overline{E_1}\times\partial E_2}|R(f_1,\ldots,f_m)|>\tilde{c}$ (cf. Step
8), we can effectively bound $\delta$ from above to ensure that $\inf_{\overline{E_1'}\times\partial
E_2}|R(f_1^{\nu},\ldots,f_m^{\nu})|>\frac{\tilde{c}}{2}.$ Next observe that if the inequality in (v) holds, then
$||f_{m+j}-\frac{{Q_j(f_1,\ldots,f_m,{f}_{m+1})}}{R(f_1,\ldots,f_m)}||_{\overline{E_1}\times\partial E_2
}<\frac{\eta(\nu)}{\tilde{c}},$} {so for $f_{m+j}^{\nu}$ as above we have
$||f_{m+j}-f_{m+j}^{\nu}{||}_{\overline{E_1'\times E_2}}\leq$\\$\frac{\eta(\nu)}{\tilde{c}}+
\big{|}\big{|}\frac{{Q_j(f_1,\ldots,f_m,{f}_{m+1})}}{R(f_1,\ldots,f_m)}-\frac{{Q_j(f_1^{\nu},\ldots,f_m^{\nu},{f}_{m+1}^{\nu})}}
{R(f_1^{\nu},\ldots,f_m^{\nu})}\big{|}\big{|}_{\overline{E_1'}\times\partial E_2}.$ Now it is sufficient to take
$\eta(\nu)<\frac{\tilde{c}}{2\nu}$ and $\delta$ so small that the second term of the righthand side of the last
inequality is bounded by $\frac{1}{2\nu}.$}

{Next, compute $\varepsilon$ in Lemma \ref{vddlem} (cf. Remark \ref{vdrema}) with $d$ equal to the degree of $P$
in $z_{m+1}$ and $M$ such that $\frac{M}{2}$ is a common bound for $f_{m+1}$ and the coefficients of
$P(f_1,\ldots,f_m,z_{m+1})\in\mathcal{O}(\overline{E_1\times E_2})[z_{m+1}]$. Define
$\theta(\nu):=\min\{\frac{\varepsilon\tilde{c}^2}{8},\frac{\delta\tilde{c}^2}{32N}\},$ where $\frac{N}{2}$ is a
bound for $\frac{\partial P}{\partial z_{m+1}}(f_1,\ldots,f_m, f_{m+1}).$ We will show that with $\eta(\nu),
\theta(\nu)$ picked above and $\frac{1}{\mu(\nu)},$ $||H_{i}^{\mu}-H_i||_{\overline{E_1'\times E_2}}$ small enough
we have $||f_i-f^{\nu}_i||_{\overline{E_1'\times E_2}}<\delta,$ for $i=1,\ldots,m+1.$}

{From $f_i=H_i(W)^2+r_{i}$ and $f_i^{\nu}=H_i^{\mu}(W_{\mu})^2+r_{i,\mu}$ (cf. Steps 10, 13), for $i=1,\ldots,m,$
we can clearly determine how close $H_i^{\mu}, g^{\mu}$ to $H_i,g$ have to be to ensure that
$||f_i-f_i^{\nu}||_{\overline{E_1'\times E_2}}<\delta,$ for $i=1,\ldots,m.$ That is, we have upper bounds for
$\frac{1}{\mu(\nu)}$ and $||H_{i}^{\mu}-H_i||_{\overline{E_1'\times E_2}}$ which imply
$||f_i-f_i^{\nu}||_{\overline{E_1'\times E_2}}<\delta,$ for $i=1,\ldots,m.$ Let us show that (for possibly smaller
bounds) we also have $||f_{m+1}-f_{m+1}^{\nu}||_{\overline{E_1'\times E_2}}<\delta.$}

{To do this, denote} $\bar{f}^{\nu}=H_{m+1}^{\mu}W_{\mu}^{2}+r_{m+1,\mu},$ where $H_{m+1}^{\mu}$ is a polynomial
approximation of $H_{m+1}.$ {We can compute an upper bound for $\frac{1}{\mu(\nu)}$ and
$||H_i^{\mu}-H_i||_{\overline{E_1'\times E_2}}$ which guarantee that} {$$P(f_1^{\nu},\ldots,f_{m}^{\nu},
\bar{f}^{\nu})= C^{\nu}(\frac{\partial P}{\partial z_{m+1}}(f_1^{\nu},\ldots,f_m^{\nu},\bar{f}^{\nu}))^{2},$$ on
$\overline{E'_1\times E_2},$ where $C^{\nu}$ is a holomorphic function. Indeed,} {observe that $W(0,x_n)=x_n^d$
and $W_{\mu}(0,x_n)$ has the form $W_{\mu}(0,x_n)=x_n^d+w_{\mu}(x_n).$ Compute how close $g^{\mu}, H_i^{\mu}$ to
$g, H_i$ {have to be} to ensure that $W_{\mu}(0,\cdot)$ has $d=\mathrm{deg}(W)$ roots in $E_2$ (counted with
multiplicities}{; by the Rouch\'e Theorem it is sufficient to ensure that $|x_n|^d>|w_{\mu}(x_n)|$ on $\partial
E_2$})
and\vspace*{2mm}\\
(3.5)\hspace*{9mm}$||\frac{\partial P}{\partial z_{m+1}}(f^{\nu}_1,\ldots,f^{\nu}_m,\bar{f}^{\nu})-\frac{\partial P}{\partial z_{m+1}}(f_1,\ldots,f_m,f_{m+1})||_{\overline{E_1'\times E_2}}<\frac{\tilde{c}}{2}.$\vspace*{2mm}\\

If the last two conditions hold, then, by the fact that $$\mathrm{inf}_{\overline{E_1'}\times\partial
E_2}|\frac{\partial P}{\partial z_{m+1}}(f_1,\ldots,f_{m+1})|>\tilde{c}$$ and by the Rouch\'e Theorem,
$\frac{\partial P}{\partial z_{m+1}}(f_1^{\nu},\ldots,f_m^{\nu},\bar{f}^{\nu})(x',\cdot)$ has $d$ roots (counted
with multiplicities) in $E_2$ for every $x'\in \overline{E_1'}$ (and no root in $\partial E_2$). {But recall that
$\frac{\partial P}{\partial z_{m+1}}(f_1^{\nu},\ldots,f_m^{\nu},\bar{f}^{\nu})$ is divisible by $W_{\mu}$ (because
$T_j\circ g^{\mu}=0,$ for $j=2d+1,\ldots, 3d$, cf. Steps 11, 12, 13), and $W_{\mu}(0,\cdot)$ has $d$ roots. Hence,
for every $x'\in\overline{E_1'}$ we have that $\frac{\partial P}{\partial
z_{m+1}}(f_1^{\nu},\ldots,f_m^{\nu},\bar{f}^{\nu})(x',\cdot)$ and $W_{\mu}(x',\cdot)$ have the same roots with the
same multiplicities. On the other hand, $P(f_1^{\nu},\ldots,f_{m}^{\nu}, \bar{f}^{\nu})$ is divisible by
$(W_{\mu})^2$ (because $T_j\circ g^{\mu}=0,$ for $j=1,\ldots, 2d$, cf. Steps 11, 12, 13) so
$$P(f_1^{\nu},\ldots,f_{m}^{\nu}, \bar{f}^{\nu})=
C^{\nu}(\frac{\partial P}{\partial z_{m+1}}(f_1^{\nu},\ldots,f_m^{\nu},\bar{f}^{\nu}))^{2},$$ on
$\overline{E'_1\times E_2},$ where $C^{\nu}$ is a holomorphic function,} {as claimed above.}

{Now, improving the obtained upper bounds for $\frac{1}{\mu(\nu)}$ and
$||H_i-H_i^{\mu}||_{\overline{E'_1\times E_2}},$ which can be done effectively, we have that\vspace*{2mm}\\
$||\frac{\partial P}{\partial z_{m+1}}(f_1^{\nu},\ldots,f_m^{\nu},\bar{f}^{\nu})||_{\overline{E'_1\times E_2}}<N,$
$||P(f_1^{\nu},\ldots,f_m^{\nu},\bar{f}^{\nu})||_{\overline{E_1'\times E_2}}<2\theta(\nu)$ (recall (u) in Step 9)
and $||f_{m+1}-\bar{f}^{\nu}||_{\overline{E_1'\times E_2}}<\frac{\delta}{2},$ and
$||\bar{f}^{\nu}||_{\overline{E_1'\times E_2}}<M,$ and the coefficients of
$P(f_1^{\nu},\ldots,f_m^{\nu},z_{m+1})\in\mathcal{O}(\overline{E_1'\times E_2})[z_{m+1}]$ are bounded by $M.$
(Recall that $\frac{N}{2}$ is a bound for $\frac{\partial P}{\partial z_{m+1}}(f_1,\ldots,f_m,{f}_{m+1})$ and
$\frac{M}{2}$ is a common bound for $f_{m+1}$ and the coefficients of
$P(f_1,\ldots,f_m,z_{m+1})\in\mathcal{O}(\overline{E_1\times E_2})[z_{m+1}].$)\vspace*{2mm}\\
When this is done then, by (3.5) and by $\mathrm{inf}_{\overline{E_1'}\times\partial E_2}|\frac{\partial
P}{\partial z_{m+1}}(f_1,\ldots,f_{m+1})|>\tilde{c}$, we have
$$||C^{\nu}||_{\overline{E_1'\times E_2}}=||C^{\nu}||_{\overline{E_1'}\times
\partial E_2}\leq\frac{||P(f_1^{\nu},\ldots,f_{m}^{\nu}, \bar{f}^{\nu})||_{\overline{E_1'\times
E_2}}}{\inf_{\overline{E_1'}\times
\partial E_2}|(\frac{\partial P}{\partial
z_{m+1}}(f_1^{\nu},\ldots,f_m^{\nu},\bar{f}^{\nu}))^{2}|}<\frac{8\theta(\nu)}{\tilde{c}^2}.$$ Consequently, (cf.
the definition of $\theta(\nu)$) we obtain $||C^{\nu}||_{\overline{E_1'\times E_2}}<\varepsilon$ and
$$2||C^{\nu}\cdot \frac{\partial P}{\partial z_{m+1}}(f_1^{\nu},\ldots,f_m^{\nu},\bar{f}^{\nu})||_{\overline{E_1'\times
E_2}}\leq\frac{\delta}{2}.$$ Hence, by Lemma \ref{vddlem}, there is
$f_{m+1}^{\nu}\in\mathcal{O}(\overline{E_1'\times E_2})$ with $P(f_1^{\nu},\ldots,f_m^{\nu},f_{m+1}^{\nu})=0$ (so
$f_{m+1}^{\nu}$ is Nash as its graph is contained in a Nash hypersurface) and
$||f_{m+1}^{\nu}-\bar{f}^{\nu}||_{\overline{E_1'\times E_2}}\leq\frac{\delta}{2}$ therefore
$||f_{m+1}-{f}_{m+1}^{\nu}||_{\overline{E_1'\times E_2}}<\delta.$ Note that
$\{(x,z_{m+1}):P(f_1^{\nu}(x),\ldots,f_m^{\nu}(x),z_{m+1})=0\}$ is contained in the image of the projection of
$V^{\nu}$ defined in Step 14 onto $\mathbf{C}^n_x\times\mathbf{C}_{z_{m+1}}.$ Thus we can effectively bound
$\frac{1}{\mu(\nu)}$ and $||H_i-H_i^{\mu}||_{\overline{E_1'\times E_2}}$ from above to ensure that
$||f_{m+1}-f_{m+1}^{\nu}||_{\overline{E_1'\times E_2}}<\delta.$}

With $f_1^{\nu},\ldots, f_{m+1}^{\nu}$ as above, we know that there are holomorphic functions
$f_{m+2}^{\nu},\ldots,f_{m+s}^{\nu}$ such that the image of $(f_1^{\nu},\ldots,f_m^{\nu},f_{m+1}^{\nu},\ldots,
f_{m+s}^{\nu})$ is contained in $V_{(m)}.$ (This follows, in view of $\inf_{\overline{E_1'}\times\partial
E_2}|R(f_1^{\nu},\ldots,f_m^{\nu})|>\frac{\tilde{c}}{2},$ by Lemma~\ref{funnyobvious} with
$L(z_{m+1},\ldots,z_{m+s})=z_{m+1}$ and $Z=\mathcal{V}(V_{(m)},H)$ where $H=(f_1^{\nu},\ldots,f_m^{\nu})$;  cf.
the discussion preceding Algorithm~1). {Therefore we have
$$f_{m+j}^{\nu}{R(f_1^{\nu},\ldots,f_m^{\nu})}={Q_j(f_1^{\nu},\ldots,f_m^{\nu},{f}_{m+1}^{\nu})}$$
on $\overline{E_1'\times E_2},$ for $j=2,\ldots,s$ (cf. Step 9; in particular the functions
$f_{m+2}^{\nu},\ldots,$ $f_{m+s}^{\nu}$ are not only holomorphic but Nash as well). As we have proved above, this
implies that $||f_{m+j}^{\nu}-f_{m+j}||_{\overline{E_1'\times E_2}}<\frac{1}{\nu},$ for $j=2,\ldots,s,$ which
completes the discussion on error control.}\pagebreak

%%%%%%%%%%%%%%%%%%%%%%%%%%%%%%%%%%%%%%%%%%%%%%%%%%%%%%%%%%%%%%%%%%%%%%%%%%%%%%%%%%%%%%%%%%%%%%%%%%%%%%%
\noindent{\textbf{3.2.3.2 Bound for the degree of output polynomials}}\vspace*{2mm}\\
We prove the claim: if the algorithm returns an approximating map $F^{\nu},$ then the degree of $P_i^{\nu}$ in
$z_i,$ $i=1,\ldots,\hat{m},$ is bounded by a constant independent of $\nu$.

If $n=0$ or $m=\mathrm{dim}(V)=0,$ then the claim holds true because if the algorithm returns an approximating map
$F^{\nu}$, then $F^{\nu}$ is constant.

{Consider the linear ordering of $\mathbf{N}^2$ defined by $(l,k)<(l',k')$ iff $l<l'$ or ($l=l'$ and $k<k'$).
Assume that $(0,0)<(n,m).$ Let $V$ be an algebraic subvariety of $\mathbf{C}^{\hat{m}}$ with $\mathrm{dim}(V)=m$
and let $F:U\rightarrow V$ be any holomorphic map, where $U$ is an open neighborhood of $0$ in $\mathbf{C}^n.$
Assume that the claim is true for {every} holomorphic map $\tilde{F}:\tilde{U}\rightarrow\tilde{V},$ where
$\tilde{U}$ is an open neighborhood of $0$ in $\mathbf{C}^l$ and $\mathrm{dim}(\tilde{V})=k$ and $(l,k)<(n,m).$ We
will show that the claim is also true for $F$, which, by induction, will complete the proof.}

First note that if the algorithm returns $F^{\nu}$ in Step 3, then, by the induction hypothesis, the degrees of
the implicit forms of the components of $F^{\nu}$ are bounded by a constant independent of $\nu$.

Next recall that in Step 4 we apply a linear change of the coordinates in $\mathbf{C}^{\hat{m}}$ resulting in
replacing the components of the input map by linear combinations of these components. Here we need to know that if
the degrees of implicit forms of Nash approximations of the components do not depend on $\nu,$ then the degrees of
implicit forms of linear combinations of these Nash approximations do not depend on $\nu$ either. This is an
immediate consequence of Lemma \ref{algNashfunct}.

Now, if we obtain $F^{\nu}$ by the recursive call of the algorithm in Step 6, then,  by the induction hypothesis,
the degrees of the implicit forms of the components of $F^{\nu}$ satisfy the requirements. Therefore in Step 6 we
may assume that $\frac{\partial P}{\partial z_{m+1}}(f_1,\ldots,f_m,f_{m+1})\neq 0$ and either $\frac{\partial
P}{\partial z_{m+1}}(f_1,\ldots,f_m,f_{m+1})(0)\neq 0$ or $\frac{\partial P}{\partial
z_{m+1}}(f_1,\ldots,f_m,f_{m+1})(0)=0.$

If $\frac{\partial P}{\partial z_{m+1}}(f_1,\ldots,f_m,f_{m+1})(0)\neq 0,$ then we may assume $F(0)\in
\mathrm{Reg}(V_{(m)})$ (see Step 7). In this case the polynomial $P_i^{\nu}$ is linear in $z_i$ for
$i=1,\ldots,m.$ (Recall that in the case $F(0)\in \mathrm{Reg}(V_{(m)})$ the algorithm has been discussed
separately before.) For $i=m+1,\ldots,m+s,$ the degree of $P_i^{\nu}$ in $z_i$ is bounded by the cardinality of
the generic fiber in $V$ over $\mathbf{C}^m.$ Indeed, by construction of $P_i^{\nu}$ for $i=m+1,\ldots,m+s,$ the
degree of $P_i^{\nu}$ in $z_i$ does not exceed the number of points in the generic fiber in ${V}^{\nu}$ over $x$
which (by linearity of $P_i^{\nu}$ in $z_i$ for $i=1,\ldots,m$) is bounded by the cardinality of the generic fiber
in $V$ over $\mathbf{C}^m.$

Now suppose that $\frac{\partial P}{\partial z_{m+1}}(f_1,\ldots,f_m,f_{m+1})(0)=0.$ The linear change of the
coordinates in Step~7 regards the domain of the approximated map and therefore does not affect the bounds we look
for. First, let us discuss $P_i^{\nu}$ for $i=1,\ldots,m.$ Since $H_i^{\mu}$ obtained in Step 13 are polynomials,
they have implicit forms of degree $1.$ The same is true for $x_n^j,$ for $j\in\mathbf{N}.$ By the induction
hypothesis, the coefficients of $W_{\mu}, r_{1,\mu},\ldots, r_{m,\mu}$ have implicit forms of degrees bounded by a
constant independent of $\mu.$ {Now recall that to calculate $P^{\nu}_i,$ for $i=1,\ldots,m,$ we perform the ring
operations in the construction of $f_{i}^{\nu}=H_i^{\mu}W_{\mu}^2+r_{i,\mu}$ in their implicit representation.}
{By Lemma~\ref{algNashfunct}, for Nash} functions $f, g$ with implicit forms $P_f, P_g,$ we obtain implicit forms
$P_{f+g}, P_{f\cdot g}$ whose degrees are bounded by a constant depending only on the degrees of implicit forms
$P_f, P_g.$ This shows the claim for $i=1,\ldots,m.$

By Step 14, for $i=1,\ldots,s,$ $P_{m+i}^{\nu}$ is the optimal polynomial describing the image $V^{\nu,i}$ of the
projection of $V^{\nu}$ to $\mathbf{C}_x^n\times\mathbf{C}_{z_{m+i}}.$ So the degree of $P_{m+i}^{\nu}$ in
$z_{m+i}$ equals the cardinality of generic fibers in $V^{\nu,i}$ over $\mathbf{C}^n_x.$ By definition of
$V^{\nu},$ the latter number is bounded by a constant depending only on the cardinality of generic fibers in
$V_{(m)}$ over $\mathbf{C}^m$ and on the degree of $P_i^{\nu}$ in $z_i,$ for $i=1,\ldots,m.$ In view of the
previous paragraph, we obtain a bound independent of $\nu.$
%TU BYLA DODATKOWA {
\begin{remark}\emph{Suppose that $R$ computed in Step 5 is reducible: $R=R_1^{\alpha_1}\cdots R_k^{\alpha_k},$ where $R_1,\ldots,R_k$ are some relatively prime non-constant polynomials.
Observe that instead of calling the algorithm with $V_{(m)}\cap\{R=0\}$ in Step 6, we may
call it simultaneously (i.e. in parallel) with every $V_{(m)}\cap\{R_j=0\}$ (each of which is
simpler than the original $V_{(m)}\cap\{R=0\}$) and recover the output for
$V_{(m)}\cap\{R=0\}$ from the output for $V_{(m)}\cap\{R_j=0\},$ for $j=1,\ldots,k.$ To do
this, we need to compute non-trivial factors of $R.$ Let us discuss this problem when
$V_{(m)}$ is of the form $V_{(m)}=\tilde{V}_{(m+1)}\cap\{\tilde{R}=0\}$ where
$\tilde{V}_{(m+1)}\subset\mathbf{C}^{m+1}\times\mathbf{C}^{s-1}$ has pure dimension $m+1$ and
has proper projection onto $\mathbf{C}^{m+1}$ and $\tilde{R}$ is the discriminant of the
optimal polynomial for $\tilde{\Phi}(\tilde{V}_{(m+1)}),$ and
$\tilde{\Phi}(z_1,\ldots,z_{m+1},z_{m+2},\ldots,z_{m+s})=(z_1,\ldots,z_{m+1},z_{m+2})$; we
deal with such $V_{(m)}$ after calling Algorithm 2 recursively with  $\tilde{V}_{(m+1)}.$}

\emph{We may assume (changing the coordinates in
$\mathbf{C}^{m+1}\approx\mathbf{C}^m\times\mathbf{C}$) that
$\{\tilde{R}=0\}\subset\mathbf{C}^m\times\mathbf{C}$ has proper projection onto
$\mathbf{C}^m.$ Set $V_{\tilde{R}}=\{(z_1,\ldots,z_m,z_{m+1})\in\mathbf{C}^m\times\mathbf{C}:
\tilde{R}(z_1,\ldots,z_m,z_{m+1})=0\}.$ Choose $R^*\in\mathbf{C}[z_1,\ldots,z_{m}]$ such that
$\{R^*=0\}\subset\mathbf{C}^{m}$ is the set of points over which the fiber in $V_{\tilde{R}}$
has not the maximal cardinality. If the discriminant of $\tilde{R}$ is non-zero, then define
$R^*$ to be the discriminant of $\tilde{R}.$}

\emph{It is clear that (after a generic change of the coordinates in
$\mathbf{C}^{s}_{z_{m+1},\ldots,z_{m+s}}$), $V_{(m)}=\tilde{V}_{(m+1)}\cap\{\tilde{R}=0\}$
satisfies (x), (y) of Step 4 and that $\{R^*=0\}\subset\{R=0\}.$ Therefore (reduced) factors
of  $R^*$  are also factors of $R.$ But $R^*$ is an iterated discriminant and there are
direct methods for computing factors of iterated discriminants (cf. \cite{LMC}) which in many
cases are much faster than computing and then factorizing $R^*.$}\end{remark}

\subsubsection{Example}\label{examples}
{The aim of this subsection is to present an example illustrating Algorithm 1. We have chosen
input data in such a way that computations} {are short enough to be presented here and still
not quite trivial.}

Let
$$P(z_1,z_2,z_3,z_4,z_5,z_6)=z_6^5+z_6z_5+z_5^6+z_1^2+z_2^2+z_3^2+z_4^2,$$
$$V=\{(z_1,z_2,z_3,z_4,z_5,z_6)\in\mathbf{C}^6:P(z_1,z_2,z_3,z_4,z_5,z_6)=0\}.$$ Let $$m_1(x_1,x_2)=\frac{1}{2}x_2-\frac{1}{4}x_1\mathrm{sin}(\frac{1}{4}x_1),\mbox{ }
m_2(x_1,x_2)=\frac{1}{2}x_2-\frac{1}{4}x_1\mathrm{cos}(\frac{1}{4}x_1),$$
$$p=1+(1-m_1^{19})^{\frac{1}{6}}.$$ Define $$F(x_1,x_2)=(f_1(x_1,x_2),f_2(x_1,x_2),f_3(x_1,x_2),f_4(x_1,x_2),f_5(x_1,x_2),f_6(x_1,x_2))$$
by setting $$f_6=m_1^5,\mbox{ } f_5=m_1(1-m_1^{19})^{\frac{1}{6}},\mbox{ } f_1=ip^{\frac{1}{2}}m_2^3,\mbox{ } f_2=ip^{\frac{1}{2}}(m_1^3+m_2^3),$$
$$f_3=-\frac{i}{2}(m_1^3+m_2^3)+im_2^3p,\mbox{ } f_4=\frac{1}{2}(m_1^3+m_2^3)+m_2^3p.$$

It is easy to verify that the image of $F$ is contained in $V$ (but not in $\mathrm{Sing}(V)$). We will illustrate
Algorithm 1 by computing a Nash approximation of $F.$ Given the formulas above we know much more about $F$ than
just the input data described in Section \ref{model}. This extra information could be used to simplify the
computations, but our aim is to show how to compute approximations when we have nothing but $P, Expand_{f_j},
M_{f_j}, U_{f_j}.$ Therefore, when performing the algorithm, {we act as if the formulas defining the components of
$F$ could be used only for computing the coefficients of the Taylor expansion of these components (i.e. the output
of the procedure $Expand$). The reader might find the formulas useful to confirm (with help of a pocket
calculator) the correctness of numerical constants such as upper and lower bounds of functions appearing in the
sequel. Then the steps of the algorithm can be performed by hand without using computer programs. Note that in
practice explicit formulas might not exist, and when exist, they might not be easy to guess.}

{We assume $U_{f_j}=D_{\frac{3}{2}}\times D_{\frac{3}{2}},$ for every $j,$ where
$D_{r}\subset\mathbf{C}$ is the disc of radius $r$ centered at $0.$ On $D_{\frac{3}{2}}$ we
have $|\sin(\frac{1}{4}x_1)|<0.4,$ so $|m_1(x_1,x_2)|< 0.9$ and $|p(x_1,x_2)|<2.1.$ Also
$|m_2(x_1,x_2)|<\frac{9}{8},$ hence $M_{f_6}=0.6,$ $M_{f_5}=1,$ $M_{f_1}=2.5,$ $M_{f_2}=3.5,$
$M_{f_3}=4.5,$ $M_{f_4}=4.5$ are bounds for the components of $F$.}

{Let $P, V$ be fixed as above (then, in particular, $\frac{\partial P}{\partial z_6}$ depends
only on $z_5, z_6$). Note that} for any holomorphic map
$F(x_1,x_2)=(f_1(x_1,x_2),\ldots,f_6(x_1,x_2))$ into $V$ with $\frac{\partial P}{\partial
z_6}(f_5(x_1,x_2),f_6(x_1,x_2))=$ $\hat{H}(x_1,x_2)(x_2-A(x_1))$ (for some $A, \hat{H}$
holomorphic in a neighborhood of zero, $\hat{H}(0,0)\neq 0$) and such that {at least} one of
$f_1,\ldots,f_4$ does not vanish identically on the graph of $A,$ the computations would be
similar to those presented below. Hence, in fact, we discuss here the whole class of
examples.

We will slightly relax our requirements regarding the output. Namely,
to simplify the computations we will let some of the output polynomials $P_i^{\nu}(x_1,x_2,z_i)$ (cf. (b) of the output of Algorithm 1)
not be monic in $z_i$ (but their leading coefficients will be non-vanishing at
the origin and their degrees in $z_i$ will be independent of the accuracy of approximation). Our aim is to approximate
$F$ on $D_{\frac{1}{3}}\times D_{\frac{1}{3}}$ with precision $10^{-5}$ (i.e. $\nu=10^5$).

{Step~1}: there is nothing to do because $V$ is a hypersurface (i.e. $m=5$).

{Step~2}: no change of the coordinates is necessary as $P$ is monic in $z_6,$ so $V$ has
proper projection onto $\mathbf{C}^{5}_{z_1\ldots z_5}.$ Here
$\Phi=\mathrm{id}_{\mathbf{C}^6}.$ (Projecting along $\mathbf{C}_{z_j},$ where $j$ is one of
$1, 2, 3, 4,$ instead of projecting along $\mathbf{C}_{z_6}$ leads to a large system of
equations in Step 7. This is because all partial derivatives of {$\frac{\partial P}{\partial
z_j}(F(x_1,x_2))=2f_j(x_1,x_2),$} for $j=1,\ldots,4,$ up to order at least $2$ vanish at $0,$
so the monic polynomial $W$ obtained in Step 5 would be of degree at least $3.$)

{Step~3}: here $P$ is already optimal and $R=5^5(z_5^6+z_1^2+z_2^2+z_3^2+z_4^2)^4+4^4z_5^5.$

{Step~4'}: confirm that the partial derivative of
$\frac{\partial P}{\partial z_6}(f_5(x_1,x_2),f_6(x_1,x_2))$ (of order $1$) with respect to $x_2$
at $0$ is non-zero.

Step 5: since the partial derivative of
$\frac{\partial P}{\partial z_6}(f_5(x_1,x_2),f_6(x_1,x_2))$ with respect to $x_2$ at $0$ is non-zero, we have
$d=1$ and (by the Weierstrass Preparation Theorem)
$$\frac{\partial P}{\partial z_6}(f_5(x_1,x_2),f_6(x_1,x_2))=\hat{H}(x_1,x_2)(x_2-A(x_1)),$$ for
$\hat{H}, A$ holomorphic in some neighborhood of $0,$ $\hat{H}(0,0)\neq 0.$ Let us estimate the size of
the domain of $\hat{H}, A.$ Set $\tilde{h}(x_1,x_2)=\frac{\partial P}{\partial z_6}(f_5(x_1,x_2),f_6(x_1,x_2))$ and $\tilde{g}(x_2)=\frac{\tilde{h}(0,x_2)}{x_2}.$ Having $Expand_{f_j},$ $U_{f_j},$ $M_{f_j},$ for $j=5, 6,$
we also have $Expand_{\tilde{h}},$ $U_{\tilde{h}}=D_{\frac{3}{2}}\times D_{\frac{3}{2}},$
$M_{\tilde{h}}=1.65$ and  $Expand_{\tilde{g}},$ $U_{\tilde{g}}=D_{\frac{3}{2}},$
$M_{\tilde{g}}={1.2}.$ With these data we confirm that
$\inf_{D_1\times \partial D_1}|\tilde{h}|>0.4>0$ and $\inf_{D_{1}}|\tilde{g}|>0.4>0.$ Therefore, the number of the solutions $x_2\in D_{1}$ of $\tilde{h}(x_1,x_2)=0$ is constant (counted with multiplicities) when $x_1$ varies in $D_{1}.$ Due to the latter inequality,
$x_2=0$ is the only solution
of $\tilde{h}(0,x_2)=0$ in ${D_{1}}$ and it has multiplicity $1$ because $\frac{\partial\tilde{h}}{\partial x_2}(0,0)\neq 0.$ Hence, $\{\tilde{h}=0\}\cap(D_{1}\times D_{1})$ is the graph of a holomorphic function defined on $D_{1}.$ Thus, we can set $U_A=D_1,$ $U_{\hat{H}}=D_1\times D_1.$

Let us compute $M_A.$ Using $Expand_{\tilde{h}}, M_{\tilde{h}}$ we check that $\inf_{D_1\times \partial D_{\frac{1}{3}}}|\tilde{h}|>0.05>0.$ This implies that $\mathrm{graph}(A)\subset D_1\times D_{\frac{1}{3}},$
hence, we can define $M_A=\frac{1}{3}.$

Step 6: by the Weierstrass Division Theorem (cf. Section \ref{model})
we have
$f_i(x_1,x_2)=(x_2-A(x_1))^2H_i(x_1,x_2)+a_i(x_1)x_2+b_i(x_1),$ where $H_i,a_i,b_i$ are holomorphic functions on
$D_{1}\times D_{1},$ $D_{1},$ respectively,
for $i=1,\ldots,6.$

Let us compute $M_{a_i}, M_{b_i}.$ Set $r_i(x_1,x_2)=a_i(x_1)x_2+b_i(x_1).$ We have (cf. Section \ref{model})
$$r_i(x_1,x_2)=\frac{1}{2\pi i}\int_{\partial D_1}\frac{f_i(x_1,s)}{(s-A(x_1))^2}\frac{(s-A(x_1))^2-(x_2-A(x_1))^2}{s-x_2}ds=$$
$$=\frac{1}{2\pi i}\int_{\partial D_1}\frac{f_i(x_1,s)}{(s-A(x_1))^2}(s+x_2-2A)ds.$$
Thus, $$||r_i(x_1,0)||_{D_1}=||b_i||_{D_1}\leq \frac{15}{4}||f_i||_{D_1\times D_1}=:M_{b_i},$$
$$||r_i(x_1,1)||_{D_1}=||a_i+b_i||_{D_1}\leq 6||f_i||_{D_1\times D_1},$$
so $$||a_i||_{D_1}\leq \frac{39}{4}||f_i||_{D_1\times D_1}=:M_{a_i}.$$ Using $Expand_{f_i}, M_{f_i}$ (which bounds
$f_i$ on $D_{\frac{3}{2}}\times D_{\frac{3}{2}}$) we compute: $||f_6||_{D_1\times D_1}\leq 0.2,$ $||f_5||_{D_1\times D_1}\leq 0.7,$ $||f_1||_{D_1\times D_1}\leq 0.7,$ $||f_2||_{D_1\times D_1}\leq 1,$
$||f_3||_{D_1\times D_1}\leq 1.5,$ $||f_4||_{D_1\times D_1}\leq 1.5.$

Step 7: we obtain $$T_1=Q_1+2Aa_1^2+
2a_1b_1,$$
where\vspace*{2mm}\\ $Q_1=6a_5(a_5A+b_5)^5+
5a_6(a_6A+b_6)^4
+2a_5a_6A+b_5a_6+a_5b_6+2Aa_2^2+2Aa_3^2+2Aa_4^2+2a_2b_2+2a_3b_3+2a_4b_4.$ (Here $a_i, b_i$ denote
new variables corresponding to the functions $a_i(x_1), b_i(x_1),$ respectively.)

$$T_2=Q_2-a_1^2A^2+b_1^2,$$ where\vspace*{2mm}\\
$Q_2=
(Aa_6+b_6)^5+(Aa_6+b_6)(Aa_5+b_5)+(Aa_5+b_5)^6
-6Aa_5(Aa_5+b_5)^5
-5Aa_6(Aa_6+b_6)^4
-2A^2a_5a_6
-Ab_5a_6-Aa_5b_6
-A^2a_2^2+b_2^2
-A^2a_3^2+b_3^2
-A^2a_4^2+b_4^2.$

$$T_3=5(a_6A+b_6)^4+a_5A+b_5.$$
(In particular, $Q_1, Q_2, T_3$ do not depend on $a_1,b_1.$)

Step 8: we could call Algorithm 1 recursively to obtain the output satisfying
all the requirements. However,
$T_1, T_2$ are of degree $2$ in $a_1, b_1$ and $T_3$ is linear with respect to $b_5,$ therefore, it is
easier to compute Nash approximations
of $a_i, b_i, A$ by a different direct method, after slight relaxing the requirements regarding the output.
Namely, the output polynomials for the approximations of $a_1, b_1$ (cf. (b)) will not be monic. Consequently, $P_1^{\nu}, P_6^{\nu}$ corresponding to
the Nash approximations $f_1^{\nu}, f_6^{\nu}$ of $f_1, f_6,$ respectively, will not be
monic either but their leading coefficients will be non-vanishing at $0$ and their degrees in $z_1, z_6,$ respectively,
will be independent of $\nu$, as required in (b).
As before, we are allowed to use only
the data representing $A, a_i, b_i$ described in Section \ref{model} (obtained by applying the effective
versions of the Weierstrass Preparation and Division Theorems in Steps 5, 6)
and the equations defining the variety (Step~7).

First, using the equation $5(a_6A+b_6)^4+a_5A+b_5=0$ we can eliminate the variable $b_5$
from $Q_1, Q_2$ obtaining $\tilde{Q}_1,\tilde{Q}_2,$ respectively.
Now the system $$\tilde{Q}_1+2Aa_1^2+2a_1b_1=
\tilde{Q}_2-a_1^2A^2+b_1^2=0$$ has the following solutions: $a_1=\frac{-i\tilde{Q}_1}{2\sqrt{\tilde{Q}_1A+\tilde{Q}_2}}$ and
$b_1=\frac{-i(A\tilde{Q}_1+2\tilde{Q}_2)}{2\sqrt{\tilde{Q}_1A+\tilde{Q}_2}}.$

We will replace $a_5(x_1), A(x_1),$  $a_j(x_1), b_j(x_1),$ for $j=2, 3, 4, 6,$ in $\tilde{Q}_1, \tilde{Q}_2$ by polynomial approximations $a_{5,\mu}(x_1), A_{\mu}(x_1),$  $a_{j,\mu}(x_1), b_{j,\mu}(x_1),$
obtaining $\tilde{Q}_1^{\mu}, \tilde{Q}_2^{\mu}$ such that
$\frac{i\tilde{Q}_1^{\mu}}{2\sqrt{\tilde{Q}_1^{\mu}A_{\mu}+\tilde{Q}_2^{\mu}}}$ and
$\frac{i(A_{\mu}\tilde{Q}_1^{\mu}+2\tilde{Q}_2^{\mu})}{2\sqrt{\tilde{Q}_1^{\mu}A_{\mu}+\tilde{Q}_2^{\mu}}}$
are holomorphic in $D_{\frac{1}{3}}.$ To do this, we compute the order $\kappa$ of zero of $\delta(x_1):=\tilde{Q}_1(x_1)A(x_1)+\tilde{Q}_2(x_1)$ at $0$ to confirm that $\kappa=6.$ Since (after Steps
5, 6) we have
$Expand_{a_i}, Expand_{b_i}, Expand_{A},$ $M_{a_i}, M_{b_i}, M_A,$ we also have
$Expand_{{\delta(x_1)}/{x_1^6}}, M_{\delta(x_1)/x_1^6}$ (where $U_{\delta(x_1)/x_1^6}=D_1$).
Therefore, we can check that $\inf_{D_{\frac{1}{3}}}|\delta(x_1)/x_1^6|>0$ which implies that $0$ is
the only root of $\delta$ in $\overline{D_{\frac{1}{3}}}.$ Hence, it is sufficient to approximate $a_5, A,$
$a_j, b_j$ in such a way that $0$ is the root of $\tilde{Q}_1^{\mu}A_{\mu}+\tilde{Q}_2^{\mu}$ of order $6$ and the root of both numerators of order (at least) $3.$ Since
$a_j, b_j, A$ depend on one variable, this can be done as follows.

In general, one can write
$a_j(x_1)=x_1^{\kappa}\alpha_j(x_1)+p_j(x_1),$ $b_j(x_1)=x_1^{\kappa}\beta_j(x_1)+q_j(x_1),$
$A(x_1)=x_1^{\kappa}\gamma(x_1)+e(x_1),$ where $p_j, q_j, e$ are polynomials of degree smaller than
$\kappa,$ and then one can replace $\alpha_j, \beta_j, \gamma$ by polynomial approximations
close enough. In the case when $0$ is the root of order (at least) $3$ of
$Aa_j+b_j$ for $j=2,3,4,6,$ (which occurs in our example) one can proceed in a slightly different way:
we have\vspace*{2mm}\\
$\tilde{Q}_1=-6\cdot 5^5a_5(Aa_6+b_6)^{20}+a_5(Aa_6+b_6)+\sum_{j=2}^42a_j(Aa_j+b_j),$\\
$\tilde{Q}_2=-4(Aa_6+b_6)^5+5^6(Aa_6+b_6)^{24}+6\cdot 5^5Aa_5(Aa_6+b_6)^{20}+\\
\hspace*{50mm}-Aa_5(Aa_6+b_6)+\sum_{j=2}^4-(Aa_j+b_j)(Aa_j-b_j),$\\
$\tilde{Q}_1A+\tilde{Q}_2=-4(Aa_6+b_6)^5+5^6(Aa_6+b_6)^{24}+\sum_{j=2}^4(Aa_j+b_j)^2.$\vspace*{2mm}\\
Therefore, if $A(x_1)a_j(x_1)+b_j(x_1)=x_1^3\cdot\gamma_j(x_1),$ for some holomorphic $\gamma_j,$ then it is
sufficient to approximate $a_5, A, a_j, \gamma_j,$ by polynomials $a_{5,\mu}, A_{\mu}, a_{j,\mu}, \gamma_{j,\mu}$ close enough and
define $b_{j,\mu}=-A_{\mu}(x_1)a_{j,\mu}(x_1)+x_1^3\cdot\gamma_{j,\mu}(x_1),$
for $j=2,3,4,6.$

Finally, Nash approximations $a_{1,\mu}, b_{1,\mu}$ of $a_1, b_1$ and (a polynomial approximation) $b_{5,\mu}$ of $b_5$ are defined by\vspace*{2mm}\\ (3.6)\hspace*{15mm}$a_{1,\mu}=\frac{-i\tilde{Q}_1^{\mu}}{2\sqrt{\tilde{Q}_1^{\mu}A_{\mu}+\tilde{Q}_2^{\mu}}},\mbox{  }\mbox{  }b_{1,\mu}=\frac{-i(A_{\mu}\tilde{Q}_{1}^{\mu}+2\tilde{Q}_{2}^{\mu})}{2\sqrt{\tilde{Q}_{1}^{\mu}A_{\mu}+
\tilde{Q}_{2}^{\mu}}},$

$$b_{5,\mu}=-5(a_{6,\mu}A_{\mu}+b_{6,\mu})^4-a_{5,\mu}A_{\mu}.$$

How close $a_{j,\mu}, b_{j,\mu}, A_{\mu}$ {have to} be to $a_j, b_j, A$ is discussed below. Note that, for $P$
defined above, the method of approximation of $A, a_j, b_j$ presented in Step~8 of this example can be applied for
any initial map $F$ for which in Step 5 we have $d=1$ and for which in Steps 5, 6 we obtain $a_1, b_1, A$ such
that $Aa_1+b_1\neq 0.$ (If $d=1,$ then $T_1, T_2, T_3$ are of the same form as in Step 7 above. $Aa_1+b_1\neq 0$
implies that $\tilde{Q}_1A+\tilde{Q}_2\neq 0.$ Clearly, the requirement $Aa_1+b_1\neq 0$ can be replaced by
$Aa_j+b_j\neq 0$ for any $j\in\{1,\ldots,4\}$.)

Step 9: for $j=2,\ldots,5,$ set
$$f_{j}^{\nu}(x_1,x_2)=(x_2-A_{\mu}(x_1))^2H_j^{\mu}(x_1,x_2)+a_{j,\mu}(x_1)x_2+b_{j,\mu}(x_1),$$ where $H_j^{\mu}$ is a polynomial approximation of $H_j$ close enough. Since $A_{\mu}, a_{j,\mu}, b_{j,\mu}$ are polynomials, for $j=2,\ldots,5,$ we can define
$P_j^{\nu}(x_1,x_2,z_j)=z_j-f_j^{\nu}(x_1,x_2).$

Let us turn to $f_1.$ Here $H_1^{\mu}$ could also be
any polynomial (or Nash) approximation of $H_1$. However, using the results of Step 8 (following from $T_1, T_2, T_3$ being of low degrees in certain variables) one can choose $H_1^{\mu}$ more
carefully so that the degree of $P_{6}^{\nu}$ (obtained in Step 10) in $z_6$ is relatively small (i.e.
equal to $5$).
We have that ${\tilde{Q}_1^{\mu}A_{\mu}+\tilde{Q}_2^{\mu}}=x_1^{\kappa}\tau_{\mu}(x_1),$ where $\tau_{\mu}$ is a polynomial non-vanishing on $D_{\frac{1}{3}}$ (if the approximations performed in Step 8 are close enough). Approximate $H_1\cdot\sqrt{\tau_{\mu}}$ by a polynomial $\eta_{\mu}$ and
put $H_1^{\mu}=\frac{\eta_{\mu}}{\sqrt{\tau_{\mu}}},$ $$f_1^{\nu}(x_1,x_2)=
(x_2-A_{\mu}(x_1))^2H_1^{\mu}(x_1,x_2)+a_{1,\mu}(x_1)x_2+b_{1,\mu}(x_1).$$

The order of zero of $\tilde{Q}_1^{\mu}$ and of $(A_{\mu}\tilde{Q}_{1}^{\mu}+2\tilde{Q}_{2}^{\mu})$ at $0$
is at least $\frac{\kappa}{2}=3,$ so
$\omega_a(x_1)=\frac{-i\tilde{Q}_1^{\mu}}{x_1^3}\mbox{ and }\omega_b(x_1)=\frac{-i(A_{\mu}\tilde{Q}_{1}^{\mu}+2\tilde{Q}_{2}^{\mu})}{x_1^3}$ are polynomials.
Therefore, we can define $${P}_1^{\nu}(x_1,x_2,z_1)=4\tau_{\mu}z_1^2-
(2\eta_{\mu}\cdot (x_2-A_{\mu})^2+\omega_ax_2+\omega_b)^2.$$ Note that $P_1^{\nu}$ is not monic. More precisely,
dividing this polynomial by $4\tau_{\mu},$ we obtain a monic polynomial but with holomorphic (not polynomial) coefficients (i.e. then the coefficients are rational with
denominators not vanishing on $D_{\frac{1}{3}}$).

Step 10: recall that
$$P(f_1^{\nu},\ldots,f_{5}^{\nu}, \bar{f}^{\nu})=
C^{\nu}(\frac{\partial P}{\partial z_{6}}(f_1^{\nu},\ldots,f_5^{\nu},\bar{f}^{\nu}))^{2},$$
where $\bar{f}^{\nu}(x_1,x_2)=(x_2-A_{\mu}(x_1))^{2}H_{6}^{\mu}(x_1,x_2)+a_{6,\mu}(x_1)x_2+b_{6,\mu}(x_1),$ and $H_{6}^{\mu}$ is a polynomial
approximation of $H_{6}$, and $C^{\nu}$ is a holomorphic function. Moreover, $|C^{\nu}|$ is
small as $a_{j,\mu}, b_{j,\mu}, A_{\mu}$ are close to $a_{j}, b_{j}, A_{},$ respectively, hence,
by Lemma \ref{vddlem}, there is a Nash function $f^{\nu}_{6}$ close to $f_6$
with $P(f_1^{\nu},\ldots,f_5^{\nu},f_{6}^{\nu})=0.$ Recall also that
$$V^{\nu}=\{(x_1,x_2,z)\in\mathbf{C}^2_{x_1,x_2}\times\mathbf{C}^{6}:
P(z)=0, P^{\nu}_i(x_1,x_2,z_i)=0 \mbox{ for } i=1,\ldots,5\}.$$
Clearly, the graph of $f_6^{\nu}$ is contained in the image $V^{\nu,6}$ of the projection of $V^{\nu}$
to $\mathbf{C}^2\times\mathbf{C}_{z_6}.$ But now we do not know whether $V^{\nu,6}$ is algebraic
because $P_1^{\nu}$ is not monic in $z_1$ (i.e. the projection need not be
proper). Nevertheless, using $P_j^{\nu}$ for $j=1,\ldots,5,$ we can easily eliminate
$z_1,\ldots,z_5$ from $P(z_1,\ldots,z_6)$ to obtain
 $$P_6^{\nu}(x_1,x_2,z_6)=\tau_{\mu}\cdot z_6^5+\tau_{\mu}\cdot f_{5}^{\nu}\cdot z_6+\tau_{\mu}\cdot (f_{5}^{\nu})^6+
\lambda+\tau_{\mu}\cdot (f_2^{\nu})^2+\tau_{\mu}\cdot (f_3^{\nu})^2+\tau_{\mu}\cdot (f_4^{\nu})^2,$$ where
$\lambda=\frac{1}{4}(2\eta_{\mu}\cdot (x_2-A_{\mu})^2+\omega_ax_2+\omega_b)^2.$
Then $\{P_6^{\nu}(x_1,x_2,z_6)=0\}$ contains $V^{\nu,6},$ hence, also the graph of $f_6^{\nu}.$\vspace*{1mm}\\

Finally, let us estimate how close $a_{j,\mu}, b_{j,\mu}, A_{\mu},$ $H_j^{\mu}$ {have to} approximate $a_j, b_j,
A, H_j,$ for $j=1,\ldots,6,$ to ensure that $||F-F^{\nu}||_{D_{\frac{1}{3}}\times D_{\frac{1}{3}}}<10^{-5}.$
First, let us follow the proof of Lemma \ref{vddlem} to observe that our requirement is satisfied if (3.7)-(3.10)
below
are fulfilled:\vspace*{2mm}\\
(3.7)\hspace*{15mm} $||f_j^{\nu}-f_j||_{D_{\frac{1}{3}}\times D_{\frac{1}{3}}}<10^{-5},$ for $j=1,\ldots,5,$\\
(3.8)\hspace*{15mm} $||\bar{f}^{\nu}-f_6||_{D_{\frac{1}{3}}\times D_{\frac{1}{3}}}<5\cdot 10^{-6},$\\
(3.9)\hspace*{15mm}
$||C^{\nu}||_{D_{\frac{1}{3}}\times D_{\frac{1}{3}}}<10^{-2},$\\
(3.10)\hspace*{15mm} $||2C^{\nu}\cdot\frac{\partial P}{\partial z_6}(f_1^{\nu},\ldots,f_5^{\nu},\bar{f}^{\nu})||_{D_{\frac{1}{3}}\times D_{\frac{1}{3}}}<5\cdot 10^{-6}.$\vspace*{2mm}

Define $Q\in\mathcal{O}(D_{\frac{1}{3}}\times D_{\frac{1}{3}})[z_6]$ by $Q(z_6)=P(f_1^{\nu},\ldots,f_5^{\nu},z_6).$
By the derivatives of $Q$ we mean the derivatives of $P(f_1^{\nu},\ldots,f_5^{\nu},z_6)$ with respect
to $z_6$. Since
$Q(\bar{f}^{\nu})=C^{\nu}\cdot (Q'(\bar{f}^{\nu}))^2,$ we have (cf. Lemma 1.6 of \cite{vD})\vspace*{2mm}\\ $Q(\bar{f}^{\nu}+C^{\nu}\cdot Q'(\bar{f}^{\nu})\cdot Y)=\\
C^{\nu}\cdot (Q'(\bar{f}^{\nu}))^2\cdot(1+Y+\sum_{i=2}^5(Q^{(i)}(\bar{f}^{\nu})/i!)\cdot (Q'(\bar{f}^{\nu}))^{i-2}\cdot (C^{\nu})^{i-1}Y^i).$\vspace*{2mm}\\
Set $$\tilde{Q}(Y)=1+Y+\sum_{i=2}^5(Q^{(i)}(\bar{f}^{\nu})/i!)\cdot (Q'(\bar{f}^{\nu}))^{i-2}\cdot (C^{\nu})^{i-1}Y^i.$$

Using $Expand_{f_5}, Expand_{f_6}, M_{f_5}, M_{f_6}$ we confirm that $||f_5||_{D_{\frac{1}{3}}\times
D_{\frac{1}{3}}}< 0.2,$ $||f_6||_{D_{\frac{1}{3}}\times D_{\frac{1}{3}}}< 0.01.$ Then, by (3.7), (3.8), (3.9) we
have: $||\tilde{Q}(-1)||_{D_{\frac{1}{3}}\times D_{\frac{1}{3}}}\leq\frac{1}{2},$ and
$||\tilde{Q}'(t)-1||_{D_{\frac{1}{3}}\times D_{\frac{1}{3}}}\leq\frac{1}{4},$ and
$\sum_{i=2}^5||\tilde{Q}^{(i)}(t)/i!||_{D_{\frac{1}{3}}\times D_{\frac{1}{3}}}\leq\frac{1}{4},$ for every
$t\in\mathcal{O}(D_{\frac{1}{3}}\times D_{\frac{1}{3}})$ with $||t||_{D_{\frac{1}{3}}\times D_{\frac{1}{3}}}\leq
2.$ Therefore, there is $y\in\mathcal{O}({D_{\frac{1}{3}}\times D_{\frac{1}{3}}})$ with
$||y||_{D_{\frac{1}{3}}\times D_{\frac{1}{3}}}\leq 2$ such that $\tilde{Q}(y)=0$ (cf. the proof of Lemma 1.5 of
\cite{vD} with $f=\tilde{Q}$). Thus, $f_6^{\nu}=\bar{f}^{\nu}+C^{\nu}\cdot Q'(\bar{f}^{\nu})\cdot y$ satisfies
$P(f_1^{\nu},\ldots,f_5^{\nu},f_6^{\nu})=0,$ and, by (3.8), (3.10), we have
$||f_6-f_6^{\nu}||_{D_{\frac{1}{3}}\times D_{\frac{1}{3}}}<10^{-5}.$ So we have proved that it is sufficient to
fulfil (3.7)-(3.10).

Now using $Expand_{f_j}, M_{f_j},$ for $j=5,6,$ confirm that $$\sup_{D_{\frac{1}{3}}\times D_{\frac{1}{3}}}|\frac{\partial P}{\partial z_6}(f_5,f_6)|<0.2\mbox{  and }\inf_{D_{\frac{1}{3}}\times \partial D_{\frac{1}{3}}}|\frac{\partial P}{\partial z_6}(f_5,f_6)|>0.13.$$ Therefore, if (3.7), (3.8) hold, then
$$\sup_{D_{\frac{1}{3}}\times D_{\frac{1}{3}}}|\frac{\partial P}{\partial z_6}(f_5^{\nu},\bar{f}^{\nu})|<0.21\mbox{ and }\inf_{D_{\frac{1}{3}}\times \partial D_{\frac{1}{3}}}|\frac{\partial P}{\partial z_6}(f_5^{\nu},\bar{f}^{\nu})|>0.12,$$ and
both (3.9) and (3.10) are satisfied if $||C^{\nu}||_{D_{\frac{1}{3}}\times D_{\frac{1}{3}}}<1.25\cdot 10^{-5}.$
In view of the Maximum Principle, it is sufficient to obtain the last inequality
on $D_{\frac{1}{3}}\times \partial D_{\frac{1}{3}},$ on which we have $$C^{\nu}=\frac{P(f_1^{\nu},\ldots,f_5^{\nu},\bar{f}^{\nu})}{(\frac{\partial P}{\partial z_6}(f_5^{\nu},\bar{f}^{\nu}))^2}.$$ Since $\inf_{D_{\frac{1}{3}}\times \partial D_{\frac{1}{3}}}|\frac{\partial P}{\partial z_6}(f_5^{\nu},\bar{f}^{\nu})|>0.12,$ it is sufficient to require\vspace*{2mm}\\
(3.11)\hspace*{15mm}$\sup_{D_{\frac{1}{3}}\times D_{\frac{1}{3}}}|P(f_1^{\nu},\ldots,f_5^{\nu},\bar{f}^{\nu})|<1.8\cdot 10^{-7}.$\vspace*{2mm}\\

Using $Expand_{f_j},$ $M_{f_j}$ we obtain the following bounds: $||f_1||_{D_{\frac{1}{3}}\times D_{\frac{1}{3}}}< 0.03,$\linebreak $||f_2||_{D_{\frac{1}{3}}\times D_{\frac{1}{3}}}< 0.05,$ $||f_3||_{D_{\frac{1}{3}}\times D_{\frac{1}{3}}}<0.06,$ $||f_4||_{D_{\frac{1}{3}}\times D_{\frac{1}{3}}}< 0.06$; recall that $||f_5||_{D_{\frac{1}{3}}\times D_{\frac{1}{3}}}< 0.2,$ $||f_6||_{D_{\frac{1}{3}}\times D_{\frac{1}{3}}}< 0.01.$ These data and the definition of $P$ imply that
(3.11) holds
if\vspace*{2mm}\\
(3.12)\hspace*{14mm}$||\bar{f}^{\nu}-f_6||_{D_{\frac{1}{3}}\times D_{\frac{1}{3}}}<2.8\cdot 10^{-7}$\vspace*{2mm}\\
 and\vspace*{2mm}\\
(3.13)\hspace*{14mm}$||f_j^{\nu}-f_j||_{D_{\frac{1}{3}}\times D_{\frac{1}{3}}}<2.8\cdot 10^{-7},$
for $j=1,\ldots,5.$\vspace*{2mm}\\ Hence, if (3.12), (3.13) are satisfied, then
(3.7)-(3.10) are also satisfied.

Now recall that after Step 6,\vspace*{2mm}\\
(3.14)\hspace*{14mm}$f_j(x_1,x_2)=(x_2-A(x_1))^2H_j(x_1,x_2)+a_j(x_1)x_2+b_j(x_1),$\\
\hspace*{23mm}for $j=1,\ldots, 6.$\vspace*{2mm}\\
Moreover, (cf. Step 9) $f_j^{\nu}, \bar{f}^{\nu}$ are of the form\vspace*{2mm}\\
(3.15)\hspace*{14mm}$f_j^{\nu}(x_1,x_2)=(x_2-A_{\mu}(x_1))^2H_j^{\mu}(x_1,x_2)+a_{j,\mu}(x_1)x_2+b_{j,\mu}(x_1),$\\ \hspace*{23mm}for $j=1,\ldots,5,$\vspace*{2mm}\\
and\vspace*{2mm}\\
(3.16)\hspace*{14mm}$\bar{f}^{\nu}(x_1,x_2)=(x_2-A_{\mu}(x_1))^2H_6^{\mu}(x_1,x_2)+a_{6,\mu}(x_1)x_2+b_{6,\mu}(x_1).$\vspace*{2mm}\\
Confirm that $||A||_{D_{\frac{1}{3}}}<0.03$ and $||H_j||_{D_{\frac{1}{3}}\times D_{\frac{1}{3}}}<1.5$ for $j=1,\ldots,6.$ These bounds and (3.14)-(3.16) imply that
(3.12), (3.13) are satisfied if
$||A-A_{\mu}||_{D_{\frac{1}{3}}}<10^{-7},$
$||H_j-H_j^{\mu}||_{D_{\frac{1}{3}}\times D_{\frac{1}{3}}}<10^{-7},$
$||a_j-a_{j,\mu}||_{D_{\frac{1}{3}}}<10^{-7},$
$||b_j-b_{j,\mu}||_{D_{\frac{1}{3}}}<10^{-7}$ for $j=1,\ldots,6.$

%******************************************************************

\end{document}